\documentclass[a4, 11pt,reqno]{amsart}

\usepackage{lmodern}
\usepackage[english]{babel}
\usepackage{microtype}
\usepackage{pifont}

\usepackage[active]{srcltx}
\usepackage[utf8]{inputenc}
\usepackage[T1]{fontenc}
\usepackage{lipsum}
\usepackage{upgreek}

\usepackage{amsmath,amssymb,amscd,amsthm,graphicx,enumerate}
\usepackage{mathrsfs}
\usepackage{amssymb,amsthm,a4wide}

\usepackage{nicefrac}

\usepackage{pifont}

\usepackage{color}
\usepackage[usenames,dvipsnames]{xcolor}
\usepackage{xcolor}
\definecolor{ultrablue}{rgb}{0.0,0.0, 1}

\usepackage{hyperref}
\hypersetup{
	linktoc=page,
	colorlinks,
	linkcolor={ultrablue},
	citecolor={ultrablue},
	urlcolor={ultrablue}
}

\hypersetup{breaklinks=true}

\allowdisplaybreaks

\usepackage[all]{xy}
\usepackage{subfig}

\usepackage{enumerate}

\usepackage{soul}
\usepackage{graphicx}

\usepackage{tikz}
\usepackage{tikzpagenodes}
\usepackage{scrlayer-scrpage}

\usepackage{tikz}
\usetikzlibrary{calc,spy}
\usetikzlibrary{shapes.geometric, arrows}
\usepackage{pgfplots}
\usetikzlibrary{intersections, pgfplots.fillbetween}

\makeatletter
\@namedef{subjclassname@2020}{%
	\textup{2020} Mathematics Subject Classification}
\makeatother

\makeatletter 
\def\l@subsection{\@tocline{2}{0pt}{3pc}{6pc}{}}
 \makeatother

\usepackage[a4paper,bindingoffset=0in,%
left=1in,right=1in,top=1.3in,bottom=1.3in,%
footskip=.2in]{geometry}

\theoremstyle{plain}

\newtheorem*{theorem*}{Theorem}

\newtheorem{lemma}{Lemma}[section]
\newtheorem{theorem}[lemma]{Theorem}
\newtheorem{proposition}[lemma]{Proposition}
\newtheorem{corollary}[lemma]{Corollary}

\theoremstyle{definition}
\newtheorem{definition}[lemma]{Definition}
\newtheorem{remark}[lemma]{Remark}
\newtheorem{example}[lemma]{Example}
\newtheorem{anzats}[lemma]{Anzats}

\theoremstyle{remark}

\numberwithin{equation}{section}

\newcommand{\R}{\mathbb R}
\newcommand{\Z}{\mathbb Z}

\newcommand{\Deg}{\mathrm{Deg}}
\newcommand{\D}{\mathrm{D}}

\DeclareMathOperator{\Ric}{Ric}

\newcommand{\N}{\mathcal N}
\newcommand{\K}{\mathcal K}
\newcommand{\M}{\mathcal M}

\DeclareMathOperator{\grad}{grad}

\DeclareMathOperator{\dvol}{dvol}

\newcommand{\bp}{\;\text{\scriptsize $\square$}\;}

\newcount\stylenum  \newdimen\styledim 
\def\varstyle#1{\mathchoice{\stylenum=0 #1}{\stylenum=1 #1}{\stylenum=2 #1}{\stylenum=3 #1}}
\def\mathaxis{\fontdimen22\ifcase\stylenum 
	\textfont\or\textfont\or\scriptfont\or\scriptscriptfont\fi2 }
\def\setstyledim{\styledim=\ifcase\stylenum .1em\or.1em\or.07em\or.05em\fi\relax}

\def\sqdot{\mathbin{\varstyle{\raise\mathaxis\hbox{\setstyledim
				\kern\styledim 
				\vrule width1.2\styledim height.6\styledim depth.6\styledim
				\kern\styledim}}}}

\newcommand{\red}{\textcolor{red}}
\newcommand{\blue}{\textcolor{blue}}

\begin{document}

\title[\footnotesize Bakry-\'Emery Ricci curvature bounds for doubly warped products]{Bakry-\'Emery Ricci curvature bounds for  doubly warped products of weighted spaces}

\address{Zohreh Fathi\\ \newline \phantom{S} Department of Mathematics and Computer Science, Amirkabir University of Technology, 424 hafez Ave., Tehran, Iran}
\email{\href{mailto:z.fathi@aut.ac.ir}{z.fathi@aut.ac.ir}}

\address{Sajjad Lakzian\\  \newline \phantom{S} Department of Mathematical Sciences\\ Isfahan University of  Technology (IUT) \\ Isfahan 8415683111, Iran}

\email{\href{mailto:slakzian@iut.ac.ir}{slakzian@iut.ac.ir}}

\address{School of Mathematics\\
	Institute for Research in Fundamental Sciences (IPM) \\
	P.O. Box 19395-5746\\
	Iran
}
\email{\href{mailto:slakzian@ipm.ir}{slakzian@ipm.ir}}

%

\thanks{*\textit{the corresponding author}}

\maketitle
\begin{center}
\begin{minipage}[t]{0.4\textwidth}
	\centering
\bf	{ \small Zohreh Fathi}\\
\textit{\footnotesize Amirkabir University of Technology}\\
\textit{ \footnotesize Tehran, Iran}\\
	
\end{minipage}
\begin{minipage}[t]{0.4\textwidth}
	\centering
	\bf	{ \small Sajjad Lakzian*}\\
	\textit{\footnotesize  Isfahan University of Technology}\\
	\textit{ \footnotesize Isfahan, Iran}\\
\end{minipage}
\end{center}
\vspace{5mm}

%

\begin{abstract}
\par We introduce a notion of doubly warped product of weighted graphs that is consistent with the doubly warped product in the Riemannian setting. We establish various discrete Bakry-\'Emery Ricci curvature-dimension bounds for such warped products in terms of the curvature of the constituent graphs. This requires deliberate analysis of the quadratic forms involved, prompting the introduction of some crucial notions such as curvature saturation at a vertex. In the spirit of being thorough and to provide a frame of reference, we also introduce the $\left(R_1,R_2\right)$-doubly warped products of smooth measure spaces and establish $\N$-Bakry-\'Emery Ricci curvature (lower) bounds thereof in terms of those of the factors.  At the end of these notes, we present examples and demonstrate applications of warped products with some toy models. 
\end{abstract}
\date{\today}
\tableofcontents

\section{Introduction}
\par A doubly twisted product of two Riemannian manifolds $\left(B^{n_1},g_B\right)$ and $\left( F^{n_2},g_F \right)$ is of the form
\begin{align*}
B {_\alpha\times_\beta} F := \left(  B \times F ,  g:= \alpha^2g_B \oplus \beta^2g_F   \right),
\end{align*}
where
\begin{align*}
\alpha:B\times F \to \R_+ \quad \text{and}  \quad \beta: B\times F  \to \R_+,
\end{align*}
are smooth twisting functions. The product manifold is called a doubly warped product when $\alpha$ and $\beta$ are independent of $B$ and $F$ (resp.). In terms of the square field operator of the doubly twisted product corresponding to the Laplace-Beltrami operator (see~(\ref{eq:Gamma-smooth}) below for the definition), for smooth functions $u,v: B \times F \to \R$ and points $x \in B$ and $p \in F$, setting $u^x(\sqdot) = u \left( x , \sqdot \right)$ and $u^p(\sqdot) = u \left( \sqdot , p \right)$ (and similarly $v^x$ and $v^p$), one observes
\begin{align}\label{eq:Gamma-twisted}
	\Gamma\left( u, v  \right)(x,p)  &= \left<\nabla u(x,p) , \nabla v(x,p)   \right>  \\  &= \alpha^{-2}(x,p)\left< ^B\nabla u(x) , ^B\nabla v(x)   \right>_B  +  \beta^{-2}(x,p)\left< ^F\nabla u(p) , ^F\nabla v(p)    \right>_F  \notag \\ &= \alpha^{-2}(x,p)\; {^B\Gamma} \left( u^p, v^p \right)(x) +  \beta^{-2}(x,p)\;{^F\Gamma} \left( u^x, v^x \right)(p). \notag 
\end{align}
The above will be the relation we will require in the definition of the doubly warped product of graphs. This will be elaborated on later in this section.
\subsubsection*{\bf \textit{Motivation}}
\par The definition of warped product of Riemannian metrics first appeared in~\cite{BO} and has since proven to be an extremely useful tool in geometry inasmuch as selecting suitable warping functions, one can make the product manifold exhibit certain desirable behavior such as globally possess non-positive curvature, have certain spectrum or provide solutions to a given geometric flow; see e.g.~\cite{BO}~\cite{ONeil}~\cite{Ej} and~\cite{IKS}. Also warped product representations of solutions to equations that arise in general relativity are imperative to the subject;~see e.g.~\cite{FGKU} and the references therein. We note that one can also define doubly twisted products of complete geodesic spaces along the lines of~\cite{Ch}.
\par A notion of warped product for weighted graphs has not yet been established. This is partially due to the fact that on the one hand, a graph while seemingly a simple object, could be equipped with vertex weights, edge weights and a variety of distances on it which do not arise from a metric tensor; this causes ambiguity in how to effectively define warped products. On the other hand, one might be curious as to why such a notion should even be studied in the first place. The ``why'' question is really outshined by the ever increasing and widespread interest in generalizing Riemannian geometry concepts to singular and discrete settings; the innate power of the geometry of discrete structures and its links to smooth geometry -- when used in the analysis of big data, network theory and machine learning -- is well established and directly influences people's mundanity. To elaborate, the reason as to why we should study warped products of graphs is twofold.
\subsubsection*{\small \bf \textit{\textsf{From theoretical point of view}}}
From the Mathematical point of view, our definition of a doubly warped product of graphs, in many ways, resembles its Riemannian counterpart; most importantly they share similar form of inner products and similar formulas for curvature. This remarkable feature makes the discrete warped products interesting objects to study. A Mathematical application of the doubly warped or twisted products of graphs would, in analogy with the smooth case, include the use of doubly warped product graphs as local models to define fibered graphs and/or graph submersions which will provide a framework for modeling what we call interplay networks or network of networks in its most general form; the product of network structures are used in matching problems in networks. 
\subsubsection*{\small \bf \textit{\textsf{From application point of view}}}
 In application, our notion of doubly warped product and the curvature bounds we present may be used to model the interplay between complex networks thus has the potential to aid in measuring the robustness of such structures. In literature, a network of networks (e.g.~\cite{CNM}) which is a particular form of an interconnection network, is one that can be modeled by a Cartesian product of weighted graphs; See e.g.~\cite{CNM} for using Cartesian products of graphs in controllability and observability of networks. 
\par A doubly warped product of graphs as we present here, is a more realistic, versatile and sophisticated model for such interaction networks. Indeed, networks (say vertical fibers) that repeat themselves in various geographical areas (geographical areas form the horizontal fibers and the phrase geographical is not literal; for example the horizontal network could consist of various locations in a cancer tissue) can be expressed as a Cartesian product of weighted graphs.
 \par If we let the vertical interactions depend on the area, we get a warped product; if the horizontal interactions depend on the nodes in the vertical fibers, we get a doubly warped product of networks. Of course we can consider multidimensional networks so as to capture as much information about the underlying system as possible. If we let the node interactions to be asymmetric, we get a doubly twisted product of networks (weighted graphs). The model can be further generalized by considering fibered graphs which locally, are doubly twisted products. 
 \par Examples are abundant; if we take the vertical network of airlines and horizontal network of airports, we get a doubly warped multidimensional network describing air travel. Taking vertical network of major cell phone brands and horizontal network of zip codes or cell phone towers, we get a doubly warped product of networks which describes the cell phone communication system. The structure of franchise companies is another example that can be described by doubly warped products. As has recently been evidenced by research works from theoretical network theory to computational biology (e.g.~\cite{Tetal}, \cite{WSJ}, \cite{PMT} and \cite{JL}), different generalizations of the Ricci curvature bounds can be successfully employed to measure the robustness of a given evolutionary or static network; here, robustness (can be quantified via different methods) is generally understood as a gauge to determine the resilience of a given network in maintaining (and/or regaining in short time) its performance in the face of change or malfunction of nodes (perturbations). So knowing the curvature bounds for doubly warped products can be directly applied to measuring the change in robustness of interplay networks. 
 \par In \S\hspace{1.2pt}\ref{sec:toymodels}, we will further expound on some of these ideas and put them to the text in a few toy models. 
\subsubsection*{\bf \textit{Goal and scope}}
\par Our primary purpose in these notes, is to introduce and study a very intuitive notion of doubly warped product of weighted graphs. We have considered the most general case; so our objects will be graphs equipped with vertex and edge weights. There will be no assumptions on the symmetricity of edge weights which makes our construction versatile enough to be used for Markov chains on finite sets as well. The idea behind this definition is to look at the Bakry's square field operator $\Gamma$ on graphs, as the generalization (to the discrete setting) of an inner product acting on gradients of functions. The inner product of gradients, as we outlined in the beginning, is exactly what the square field operator is in the Riemannian setting. As the experts in non-smooth differential geometry attest to, this is a usual way to link classical Riemannian geometry concepts to the ones on metric measure spaces or as in this paper, on graphs.\\
\subsubsection*{\bf \textit{Doubly warped product of weighted graphs}}
Here, a weighted graph is a triple $\left( G,\omega, m   \right)$ where $G$ is countable set of vertices, $\omega$ represents edge weights and $m$, vertex measure (see page 5 for a more precise definition). For two graphs $G$ and $H$, $G \bp H$ means the Cartesian product graph.
\begin{definition}\label{def:graph-dw}
	For two weighted graphs $\left( G_1 , \omega^{G_1}, m^{G_1} \right)$ and $\left( G_2 , \omega^{G_2},  m^{G_2} \right)$, and twisting functions $\alpha , \beta: G_1\bp G_2\to \R_+$, we define
	\begin{align*}
	G_1 \; {_\alpha\diamond_\beta} \; G_2 := \left( G_1 \bp G_2 , \omega, m   \right),
	\end{align*}
	where the edge weights $\omega$ and vertex weights $m$ are given by
	\begin{align*}
	\omega_{ \left( (x,p) , (y,q)  \right)} := \delta_{pq}m^{G_2}(p)\alpha^{-2}(x,p)\omega^{G_1}_{xy} + \delta_{xy}m^{G_1}(x)\beta^{-2}(x,p)\omega^{G_2}_{pq},
	\end{align*}
	and 
	\begin{align*}
	m\left(  (x,p) \right) := m^{G_1}(x)m^{G_2}(p).
	\end{align*}
	Sometimes, we write 
	\begin{align*}
	\omega = m^{G_2}\alpha^{-2}\omega^{G_1}_{xy} \oplus m^{G_1} \beta^{-2}\omega^{G_2}_{pq},
	\end{align*}
	for brevity. When $\alpha$ and $\beta$ are independent of $G_1$ and $G_2$ (resp.), the product graph is called a doubly warped product.
\end{definition}
To demystify the definition of doubly twisted product of graphs, we note that as we will show in Lemma~\ref{lem:D-G-one}, for functions $u,v : G_1\bp G_2 \to \R$ and vertices $x \in G_1$ and $p \in G_2$,
\begin{align*}
\Gamma \left( u,v  \right)(x,p) =\alpha^{-2}(x,p) \Gamma^{G_1}\left( u^p , v^p  \right)(x)  + \beta^{-2}(x,p) \Gamma^{G_2}\left( u^x , u^x  \right)(p),
\end{align*}
which is consistent with the Riemannian version (\ref{eq:Gamma-twisted}). As for the measures, in the Riemannian setting one has
\begin{align*}
\dvol_{\alpha^2g \oplus \beta^2h} = \alpha^{n_1}\beta^{n_2}\dvol_g\dvol_h;
\end{align*}
so, our choice of measure, $m = m^{G_1} m^{G_2}$, is again consistent with the one in the smooth setting since graphs are discrete objects and can meaningfully be considered $0$-dimensional. In these notes, we are only interested in doubly warped products so $\alpha$ and $\beta$ are independent of $B$ and $F$ (resp.)  unless otherwise is specified. 
\subsubsection*{\bf \textit{Weighted doubly warped product of weighted manifolds}}
The Bakry-\'Emery Ricci tensor and the comparison geometry it brings about has gained a lot of steam in recent years. For geometric implications of Bakry-\'Emery Ricci curvature lower bounds, see e.g.~\cite{WW,Lo,MuWa}. In order to have a frame of reference in the Riemannian setting for the reader to compare the discrete results with as well as for their inherent attractiveness as generalizing objects, we will first consider the $\N$-Bakry-\'Emery curvature bounds for general $\left( R_1 , R_2 \right)$-doubly warped products of weighted Riemannian manifolds. These are generalizations of the $N$-warped products of~\cite{CK}. 
\begin{definition}\label{defn:setup}
	For real numbers $R_1$ and $R_2$, the $\left( R_1 , R_2   \right)$-doubly warped product of weighted Riemannian manifolds $\left(B, g_B , e^{-\Phi}\dvol_{g_B} \right)$ and $\left(  F, g_F , e^{-\Psi}\dvol_{g_F} \right)$ is the smooth measure space
	\begin{align*}
	B\;{_\alpha^{R_1}\times_\beta^{R_2}}\;F := \left( B \times F \; ,\;  g := \alpha^2 g_B \oplus \beta^2 g_F \; , \; \alpha^{R_1-n_1}\beta^{R_2-n_2} e^{-\Phi}e^{-\Psi} \dvol_g =: e^{-\chi} \dvol_g \right),
	\end{align*}
	where $\chi = \left(n_1 - R_1 \right)a + \left( n_2 - R_2 \right)b + \Phi + \Psi$, $a := \ln \alpha$ and $b := \ln \beta$.
	Notice that since the warping functions are positive, there is no restriction on the real numbers $R_i$.
\end{definition}
We establish Bakry-\'Emery Ricci curvature lower bounds in terms of those of the underlying factors provided the weights $R_i$ ($i = 1,2$) are positive and warping functions satisfy suitable partial differential inequalities or as we call them, concavity conditions. 
\par Let us briefly recall that, given a complete weighted manifold $\left( M^n, g, e^{-\Phi}\dvol_g \right)$ and $n < \N$, the $\mathcal{N}-$Bakry-\'Emery Ricci tensor is then defined in the interior of $M^n$ as
\begin{align*}
	^{M}\Ric_\Phi^{\mathcal{N}} = \Ric + \nabla^2 \Phi - \left(   \mathcal{N} - n\right)^{-1}\nabla\Phi \otimes \nabla \Phi.
\end{align*}
The case $n = \N$ can only be defined when the logarithmic density function $\Phi$ is a constant function; in the latter case, the above tensor cincides with the Ricci tensor. The other extreme is when $\N = \infty$; in this case, we let $\Ric_\Phi^{\infty} = \Ric + \nabla^2 \Phi$. 
\par It is worth mentioning that since the curvature bounds we derive are pointwise bounds, similar line of reasoning can be used for negative weights $R_i$ to obtain lower bounds however the partial differential inequalities would have to be two sided. In order to avoid redundant arguments, we will only work with positive constants $R_i$. 
\par The $\left( R_1 , R_2   \right)$-doubly warped product can also be defined in the setting of geodesic metric-measure spaces however, the more complicated behavior of geodesics (compared to the geodesics in a warped product), makes it more arduous to obtain weak Ricci curvature bounds via the theory of optimal transport; the Bakry-\'Emery curvature dimension bounds, however, could be obtained by similar calculations as we do for graphs. So at least in $RCD(K,N)$ spaces, where Lott-Sturm-Villani curvature dimension bounds and Bakry-\'Emery coincide~\cite{AGS1, AGS2}, similar calculations will yield curvature-dimension bounds. For curvature bounds of singly warped products of singular or non-Riemannian spaces, see e.g.~\cite{Ch,AB,CK}.
\subsection*{Summary of main results}
Here we will provide a list of our main results. The set up and definitions needed in the statement of these results are explained in \textsection\hspace{1.2pt}\ref{sec:prelim-smooth}.
\par In these notes, an \emph{iso-dimensional} bound refers to a bound in which the synthetic dimensions are the same for the product space and for the constituent factors whereas a \emph{non iso-dimensional} bound refers to the case where the synthetic dimensions are different. 
\subsubsection*{\bf \textit{Smooth setting}}
In this section, we will provide lower bounds for the $\N$-Bakry-\'Emery Ricci tensor in a weighted doubly warped product of weighted manifolds. The upper bounds on this tensor are not treated here and geometrically speaking, the upper bounds are less important than the lower bounds; e.g. see~\cite{Lohkamp} for existence and denseness of Ricci negatively curved Riemannian metrics. Also the lower bounds on the $\N$-Bakry-\'Emery Ricci tensor are very important since they are consistent with the more general framework of curvature-dimension bounds obtained via theory of optimal transportation (the so called Lott-Sturm-Villani curvature bounds). 
\subsubsection*{\small \bf \textit{\textsf{Notations and terminologies:}}}
\par In the geometric quantities that we consider, a superscript on the upper left corner indicates the space in which the quantity is being computed; if there is no spacial superscript present, the quantity is computed in the (doubly warped) product space. For example, $^{B}\nabla^2 b$ means the Hessian of $b$ in $B$ while $\nabla^2 b$ is the Hessian of $b$ in the product space. 
\par In what follows, a dynamic convexity/concavity condition on a quantity merely refers to a  pointwise inequality in terms of the Hessian of that quantity.
\subsubsection*{\small \bf \textit{\textsf{Lower bounds on the Ricci tensor}}}
Consider the doubly warped product manifold $B\;{_\alpha^{R_1}\times_\beta^{R_2}}\;F$ with $\dim B = n_1$ and $\dim F = n_2$. We have the following very versatile bound on the Ricci tensor. The conditions and bounds are all pointwise at a given $(x,p) \in B \times F$. Recall $a := \ln \alpha$ and $b := \ln \beta$. In what follows, $\K_B$ and $\K_F$ denote poitwise lower bounds on the Ricci tensor in $B$ and $F$ respectively; ${^B\Ric}$ and ${^F\Ric} $ denote the Ricci tensors of the said manifolds. 
\par  Also in below, the concavity conditions are tested on the unit tangent bundle; so, $\left( \nabla b  \right)^2$ acting on a unit tangent vector $X$ is $(\nabla_Xb)^2$ and the Hessian $\nabla^2 b$ acting on a unit tangent vector $X$ is to be understood as $\nabla^2b \left(X, X\right)$. 
\begin{theorem}\label{thm:dw-ric-1}
	Suppose ${^B\Ric} \ge \left( n_1 - 1 \right)\K_B\;g_B$ at $x \in B$ and ${^F\Ric} \ge \left( n_2 - 1 \right)\K_F\; g_F $ at $p \in F$. If for given real numbers $K_1, K_2, L_1, L_2 \in \R$, $a$ and $b$ satisfy dynamic concavity conditions
	\begin{align}\label{eq:dw-ric-1}
		&	{^B\nabla}^2 b + \Big( n_2^{-1}\left(n_1 -2\right) + 2 \Big)\left( \nabla b  \right)^2 \\ & \phantom{sajjad}\le \min \left\{-\alpha^{-2}n_2^{-1}K_1 , \left( n_1^{-1}\left(n_2 -2\right) + n_1^{-1}\Big( n_2^{-1}\left(n_1 -2\right) + 2 \Big)\right)\| \nabla b \|_B^2 - n_1^{-1}\alpha^{-2}L_1 \right\}, \notag
	\end{align}
	on $UT_xB$ (unit tangential sphere at $x$ w.r.t. $g_B$) and
	\begin{align}\label{eq:dw-ric-2}
		&	{^F\nabla}^2 a + \Big( n_1^{-1}\left(n_2 -2\right) + 2 \Big)\left( \nabla a  \right)^2 \\ & \phantom{sajjad}\le \min \left\{-\beta^{-2}n_1^{-1}K_2 , \left( n_2^{-1}\left(n_1 -2\right) + n_2^{-1}\Big( n_1^{-1}\left(n_2 -2\right) + 2 \Big)\right)\| \nabla a \|_F^2 - n_2^{-1}\beta^{-2}L_2 \right\}, \notag
	\end{align}
on $UT_pF$  (unit tangential sphere at $p$ w.r.t. $g_F$), then
	\begin{align*}
		\Ric \ge \left( n_1 + n_2 -1   \right)\K g \quad at \quad  (x,p),
	\end{align*}
	where
	\begin{align*}
		\K :=  \frac{\left(n_1 - 1  \right)\K_B + K_1 + L_1}{\left( n_1 + n_2 -1   \right)\alpha^2} \wedge \frac{\left(  n_2 - 1\right) \K_F + K_2 + L_2}{\left( n_1 + n_2 -1   \right)\beta^2}. 
	\end{align*}
\end{theorem}
\begin{remark}
	Notice that the concavity conditions in Theorem~\ref{thm:dw-ric-1}, as well as in the following results, are indeed upper bound conditions on the largest eigenvalue of the quadratic forms that appear on the left hand side.  
\end{remark}
The power in the above theorem is that by tweaking the parameters involved in~\ref{eq:dw-ric-2}, one can get various bounds.
\begin{corollary}
	If ${^B\Ric} \ge \left( n_1 - 1 \right)\K_B\;g_B$, ${^F\Ric} \ge \left( n_2 - 1 \right)\K_F\;g_B $ and $a,b$ satisfy the concavity relations
	\begin{align*}
		{^B\nabla}^2 b  \le - \Big( \frac{n_1 -2}{n_2} + 2 \Big)\left( \nabla b  \right)^2  \quad and \quad {^F\nabla}^2 a  \le -  \Big( \frac{n_2 -2}{n_1} + 2 \Big)\left( \nabla a  \right)^2,
	\end{align*}
	on $UT_xB \oplus UT_pF$, then 
	\begin{align*}
		\Ric \ge \left( n_1 + n_2 -1   \right)\K g \quad at \quad (x,p),
	\end{align*}
	where
	\begin{align*}
		\K :=  \frac{\left(n_1 - 1  \right)\K_B}{\left( n_1 + n_2 -1   \right)\alpha^2} \wedge \frac{\left(n_2 - 1  \right)\K_F}{\left( n_1 + n_2 -1   \right)\beta^{2}}.
	\end{align*}
\end{corollary}
\begin{proof}
	Set $K_1 = K_2 = L_1 = L_2 = 0$ in Theorem~\ref{thm:dw-ric-1}.
\end{proof}
\begin{corollary}
	If ${^B\Ric} \ge \left( n_1 - 1 \right)\K_B\;g_B$, ${^F\Ric} \ge \left( n_2 - 1 \right)\K_F\;g_B $ and $a,b$ satisfy the concavity relations
\begin{align*}
		{^B\nabla}^2 b  + \Big( \frac{n_1 -2}{n_2} + 2 \Big)\left( \nabla b  \right)^2 \le -\rho_1\alpha^{-2} \quad and \quad {^F\nabla}^2 a  \le -  \Big( \frac{n_2 -2}{n_1} + 2 \Big)\left( \nabla a  \right)^2 \le -\rho_2\beta^{-2},
\end{align*}	
on $UT_xB$ and  $UT_pF$ resp. for $\rho_1, \rho_2 \ge 0$, then 
	\begin{align*}
		\Ric \ge \left( n_1 + n_2 -1   \right)\K g \quad at \quad (x,p),
	\end{align*}
	where
	\begin{align*}
		\K :=  \frac{\left(n_1 - 1  \right)\K_B + \left( n_1 + n_2\right)\rho_1}{\left( n_1 + n_2 -1   \right)\alpha^2} \wedge \frac{\left(n_2 - 1  \right)\K_F + \left(n_1 + n_2\right)\rho_2}{\left( n_1 + n_2 -1   \right)\beta^{2}}.
	\end{align*}
\end{corollary}
\begin{proof}
	Set $K_1 = n_2\rho_1$, $L_1  = n_1\rho_1$, $K_2 = n_1\rho_2$ and $L_2 = n_2\rho_2$ in Theorem~\ref{thm:dw-ric-1}.
\end{proof}
\subsubsection*{\small \bf \textit{\textsf{Lower bounds on the $\N$-Bakry-\'Emery Ricci tensor}}}
Consider the doubly weighted doubly warped product $B\;{_\alpha^{R_1}\times_\beta^{R_2}}\;F $ (see Definition~\ref{defn:setup}). In what follows, $\K^{\N}_B$ and $ \K^{\N}_F$ denote pointwise lower bounds on the $\N$-Bakry-\'Emery Ricci tensor of $B$ and $F$ respectively.
\begin{theorem}[Iso-dimensional bounds]\label{thm:main-1-1}
	Let $\N > n_1 + n_2$ ($i = 1,2$). Let $R_1, R_2 \in \R_+$ and $K_1, K_2 \in \R$ be arbitrary. There exists constant $\lambda^{\N}$ only depending on $\N$ and on $n_i,R_i, K_i$ ($i=1,2$) so that if
\begin{align*}
		{^B\Ric}_\Phi^{\N} \ge \left( n_1 - 1   \right)\K^{\N}_B\; g_B\;\; at\;\; x\in B, \quad {^F\Ric}_{\Psi}^{\N} \ge \left( n_2 - 1   \right), \K^{\N}_F\;g_F \;\;at\;\; p \in F
	\end{align*}
	and the dynamic concavity conditions	
	\begin{align}\label{eq:main-con-rel-1-1}
		&R_2{^B\nabla}^2 b  - \lambda^{\N} \left( \nabla \Phi \right)^2  - \lambda^{\N} \left( \nabla b \right)^2 \le \\ &\min\left\{ - \alpha^{-2}K_1, \right. \notag \\ &\left.\phantom{sajjad} - n_1^{-1}\left( \lambda^{\N} + \nicefrac{1}{2}R_2 \right)  \| \nabla \Phi \|_B^2 - n_1^{-1}\Big( \lambda^{\N} + R_2\left(R_2 -2n_2 + \nicefrac{5}{2}   \right)\Big)\|\nabla b\|_B^2 - n_1^{-1}R_2L_1 \right\}, \notag
	\end{align}	
	and
	\begin{align}\label{eq:main-con-rel-1-2}
		&R_1{^F\nabla}^2 a  - \lambda^{\N} \left( \nabla \Psi \right)^2  - \lambda^{\N} \left( \nabla a \right)^2 \le \\ &\min\left\{ - \beta^{-2}K_2, \right. \notag \\ &\left. \phantom{sajjad} - n_2^{-1}\left( \lambda^{\N} + \nicefrac{1}{2} R_1\right)\| \nabla \Psi \|_F^2 - n_2^{-1}\Big( \lambda^{\N} + R_1\left(R_1 -2n_1 + \nicefrac{5}{2}  \right)\Big)\|\nabla a\|_F^2 - n_2^{-1}R_1L_2 \right\}, \notag
	\end{align}	
		hold on $UT_xB$ and $UT_pF$ resp. (here, $UT_x$ and $UT_pF$ are as in Theorem~\ref{thm:dw-ric-1}), then
\begin{align*}
	\Ric_\chi^{\N} \ge (n_1 + n_2 -1 ) \K  g \quad \text{at} \quad (x,p),
\end{align*}
where
\begin{align*}
	\K =  \frac{\left( n_1 - 1  \right)\K^{\N}_B + K_1 + L_1 }{\left( n_1 + n_2 - 1  \right)\alpha^2} \wedge \frac{ \left( n_2 - 1  \right)\K^{\N}_F + K_2 + L_2}{\left( n_1 + n_2 - 1  \right)\beta^2} .
\end{align*}
\end{theorem}
\begin{theorem}[Non iso-dimensional bounds]\label{thm:main-1-2}
	Let $\N_1 > n_1$, $\N_2 > n_2$. For arbitrary $R_1, R_2 \in \R_+$ and $K_1, K_2 \in \R$, there exists constant $\lambda^{\N_1,\N_2}$ only depending on $n_i,R_i,\N_i, K_i$ so that if
	\begin{align*}
			{^B\Ric}_\Phi^{\N_1} \ge \left( n_1 - 1   \right)\K^{\N_1}_B\;g_B\;\; at\;\; x\in B, \quad {^F\Ric}_\Psi^{\N_2} \ge \left( n_2 - 1   \right) \K^{\N_2}_F\; g_F \;\;at\;\; p \in F
		\end{align*}
and the concavity conditions \eqref{eq:main-con-rel-1-1} and \eqref{eq:main-con-rel-1-2}, with $\lambda^{\N}$ replaced by $\lambda^{\N_1,\N_2}$, hold on $UT_xB$ and $UT_pF$ resp., then 
	\begin{align*}
		\Ric_\chi^{\N} \ge (n_1 + n_2 -1 ) \K  g \quad \text{at} \quad (x,p),
	\end{align*}
	where
	\begin{align*}
		\K =  \frac{\left( n_1 - 1  \right)\K_B^{\N_1} + K_1 + L_1}{\left( n_1 + n_2 - 1  \right)\alpha^2} \wedge \frac{ \left( n_2 - 1  \right)\K_F^{\N_2} + K_2 + L_2}{\left( n_1 + n_2 - 1  \right)\beta^2} .
	\end{align*}
\end{theorem}
\begin{corollary}\label{cor:main-1}
Suppose the numbers $n_i,R_i,\N_i, K_i$ are as above. There exists constants $\eta^\N$ only depending on $n_i,R_i,\N_i, K_i$ so that if
		\begin{align*}
		{^B\Ric}_\Phi^{\N} \ge \left( n_1 - 1   \right)\K^{\N}_B\;g_B, \quad {^F\Ric}_\Psi^{\N} \ge \left( n_2 - 1   \right) \K^{\N}_F\;g_F,
		\end{align*}
	and the concavity relations
		\begin{align*}
		R_2 {^B\nabla}^2 b \le  \eta^{\N} \left( \nabla \Phi \right)^2 + \eta^{\N} \left( \nabla b \right)^2 \quad and \quad  R_1{^F\nabla}^2 a \le   \eta^{\N} \left( \nabla \Psi \right)^2  +   \eta^{\N} \left( \nabla a \right)^2,
		\end{align*}
		hold on $UT_xB$ and $ UT_pF$ resp., then
	\begin{align*}
	\Ric_\chi^{\N} \ge \Big( \left( n_1 - 1  \right)\alpha^{-2}\K_B \wedge \left( n_2 - 1  \right)\beta^{-2}\K_F \Big)  g \quad \text{at} \quad (x,p). 
	\end{align*}
\end{corollary}
\begin{proof}
	Set 
	\begin{align*}
		\eta^{\N} := \min\left\{ \lambda^{\N},  - n_1^{-1}\lambda^{\N}, - n_12^{-1}\lambda^{\N}, - n_1^{-1}R_2 \left(R_2 -2n_2 + 2   \right), - n_2^{-1}R_1\left(R_1 -2n_1 + 2 \right)  \right\}.
	\end{align*}
Then by the hypotheses, we have
\begin{align*}
	R_2{^B\nabla}^2 b  \le  \mu^{\N} \left( \nabla \Phi \right)^2  + \mu^{\N} \left( \nabla b \right)^2 \le  \lambda^{\N} \left( \nabla \Phi \right)^2  + \lambda^{\N} \left( \nabla b \right)^2,
\end{align*}	
and
\begin{align*}
	R_2{^B\nabla}^2 b  \le \mu^{\N} \left( \nabla \Phi \right)^2  + \mu^{\N} \left( \nabla b \right)^2 \le   - n_1^{-1}\lambda^{\N}\| \nabla \Phi \|_B^2 - n_1^{-1}\Big( \lambda^{\N} + R_2\left(R_2 -2n_2 + 2   \right)\Big)\|\nabla b\|_B^2,
\end{align*}	
along with similar concavity statements for $a$. 
\par Therefore, the hypotheses of Theorem~\ref{thm:main-1-1} are satisfied in the case where $K_1= K_2= L_1 = L_2 = 0$ and we get the desired conclusion. 
\end{proof}
\begin{corollary}\label{cor:main-2}
	Suppose the numbers $n_i,R_i,\N_i, K_i$ are as above. There exists a constant $\eta^{\N_1,\N_2}$ only depending on $n_i,R_i,\N_i, K_i$ so that if
	\begin{align*}
			{^B\Ric}_\Phi^{\N_1} \ge (n_1 - 1) \K^{\N_1}_B\;g_B , \quad {^F\Ric}_\Psi^{\N_2} \ge (n_2 - 1)\K^{\N_2}_F\;g_F,
		\end{align*}
	and the concavity conditions 
		\begin{align*}
			R_2{^B\nabla}^2 b  \le  \eta^{\N_1,\N_2} \left( \nabla \Phi \right)^2  + \eta^{\N_1,\N_2} \left( \nabla b \right)^2 \quad and \quad  R_1{^F\nabla}^2 a \le  \eta^{\N_1,\N_2} \left( \nabla \Psi \right)^2  +  \eta^{\N_1,\N_2} \left( \nabla a \right)^2,
		\end{align*}
		hold on $UT_xB $ and $UT_pF$,  then
	\begin{align*}
		\Ric_\chi^{\N} \ge \Big( \left( n_1 - 1  \right)\alpha^{-2}\K_B \wedge \left( n_2 - 1  \right)\beta^{-2}\K_F \Big)  g \quad \text{at} \quad (x,p). 
	\end{align*}
\end{corollary}
\begin{proof}
The proof is similar to the proof of Corollary~\ref{cor:main-1} hence, it is omitted. 	
\end{proof}
\begin{remark}
	The constants $\lambda^\N$ and $\lambda^{\N_1,\N_2}$ are indeed lowest eigenvalues of certain quadratic forms that will be obtained in \S\hspace{1.2pt}\ref{sec:NBE-tensor}. For some bounds on these constants and comparison of these constants, see the end of \S\hspace{1.2pt}\ref{sec:NBE-bounds}. 
\end{remark}
\subsubsection*{\bf \textit{Discrete setting}}
\par Our first graph curvature result, is a generalization and sharpening of the structural curvature bounds (curvature bounds in terms of the structure of the graph) of~\cite{LY}. Of course an immediate consequence of the structural curvature bounds is that the best lower bound  $\K_{G,x}(\N)$ at a vertex $x$, is well-defined for the most general discrete Laplacian. 
\par Let us briefly recall that if 
\[
	\Gamma_2 \left(  f \right)(x) \ge \N^{-1}\left( \Delta f(x)\right)^2 + \K\left( \N \right)\Gamma \left(  f \right)(x),
\]
holds for all functions, then, $\K$ is a lower bound for the $\N$-Bakry-\'Emery Ricci curvature at the vertex $x$. Hence, $\K_{G,x}(\N)$ is largest such $\K$. 
See \S\hspace{1.2pt}\ref{sec:disc-curv} for details on the definition of curvature bounds and for more on the terminologies involved. 
\par  Throughout these notes, $\K_{G_1,x}$, $\K_{G_2,p}$ and $\K_{(x,p)}$ denote the best (largest) lower curvature bound i.e. the curvature functions at $x\in G_1$, $p \in G_2$ and $\left(x,p  \right) \in G_1 {_\alpha\diamond_\beta} G_2$ respectively (these are functions of $\N$). Notice that here $\N \in (0,\infty]$. 
\par In the discrete setting, we are interested in both upper and lower bounds on the curvature functions; this is due to the fact that the curvature function is a result of an optimization problem and as a result, there is no formula to compute it in general; therefore finding good upper and lower bounds helps ball parking the curvature function more closely. 
\subsubsection*{\small \bf \textit{\textsf{Structural curvature bounds}}}
With the notations introduced in \S\hspace{1.2pt}\ref{sec:disc-curv}, the following bounds hold. Recall $m$ is a measure on vertices and $\D$ is the degree measure on vertices which is the sum of the edge weights emanating from $x$ divided by the $m$-measure of the vertex $x$, see  \S\hspace{1.2pt}\ref{sec:disc-curv}.
\begin{theorem}[Structural upper bounds]\label{thm:main-3}
	The curvature function at any vertex $x$ in any weighted graph $G$ (with possibly asymmetric edge weights) satisfies 
	\begin{align*}
	\K_{G,x}(\N) &\le \K_{G,x}(\infty) \\ &\le \nicefrac{1}{4}\, m^{-1}_x\D_x \max_{y \sim x} m_y  +  \nicefrac{1}{2}\left( m_x^{-\nicefrac{1}{2}}\D_x^{\nicefrac{1}{2}} - 1  \right) \max_{y\sim x}   m_y  \D_y^{\nicefrac{1}{2}} + \nicefrac{3}{4} \max_{y\sim x} m_y  \D_y,
	\end{align*}
	for all $ \N > 0$.  
\end{theorem}
\begin{theorem}[Structural lower bounds]\label{thm:main-4}
The curvature function at any vertex $x$ in any weighted graph $G$ (with possibly asymmetric edge weights) satisfies 
\begin{align*}
&\K_{G,x}(\N) \ge \\ & \min_{y \sim x} \Big[  -\nicefrac{1}{4}\, \D^2_y + \nicefrac{1}{2}\, \D_y^{\nicefrac{3}{2}} + \left(  m_x^{-\nicefrac{1}{2}}\D_x^{\nicefrac{1}{2}} - \nicefrac{1}{4}   \right)\D_y -  \left(m_x^{-\nicefrac{1}{2}}\D_x^{\nicefrac{1}{2}} \right) \D_y^{\nicefrac{1}{2}}  -  m(y)^{-1} \left(m_x^{-\nicefrac{1}{2}}\D_x^{\nicefrac{1}{2}}\right)  \Big],
\end{align*}
for all $ \N \ge 2$, and
\begin{align*}
\K_{G,x}(\N) \ge  \K_{G,x}(2) - \N^{-1}\left( 2 - \N \right) \D_x\quad   \quad \forall\;\; 0< \N \le 2.
\end{align*}
\end{theorem}
\subsubsection*{\small \bf \textit{\textsf{Upper bounds on the curvature function}}}
\par In order to get more accurate bounds, it is necessary to distinguish between various types of vertices based on the curvature maximizers at those vertices; a curvature maximizer at a vertex $x \in G$ is referred to a function for which 
\begin{align*}
	\Gamma_2(f)(x) = \N^{-1}\left( \Delta f  \right)^2(x) + \K_{G,x} \Gamma (f)(x),
\end{align*}
holds; to this end, we introduce the important notion of curvature saturated vertices.
\begin{definition}[$\N$-curvature saturated vertices]\label{DEF:saturated}
	A vertex $z \in G$ is called 	
	\begin{enumerate}
		\item[\textbf{(i)}] \emph{weakly $\N$-curvature saturated} if there exists $f:G \to \R$ (curvature maximizer) with
		\begin{align*}
			\Gamma_2 \left(  f \right)(z) = \N^{-1}\left( \Delta f(z)\right)^2 + \K_{G,z}\left( \N \right)\Gamma \left(  f \right)(z),
		\end{align*}
		that is harmonic at $z$ i.e. $\Delta f(z) = 0$,
		\item[\textbf{ii)}]  \emph{strongly $\N$-curvature saturated} if all curvature maximizers at $z$ are harmonic at $z$,
		\item[\textbf{iii)}] \emph{$\N$-curvature un-saturated} if all curvature maximizers $f$ at $z$ satisfy $\Delta f(z) \neq 0$.
	\end{enumerate}
\end{definition}
We obtain the following relations between the optimal curvature bounds of the doubly warped product and those of the constituent factors. These relations are reminiscent of what is known in the smooth setting. In what follows we use $\wedge$ and $\vee$ to denote $\min$ and $\max$ resp. 
\begin{theorem}[Non iso-dimensional upper bounds]\label{thm:main-5}
	\begin{align*}
	\K_{(x,p)}\left( \N_1 + \N_2 \right) \le \begin{cases}\Big(\alpha^{-2} \K_{G_1,x}\left( \N_1  \right) + \alpha^2 \mathcal{Q}_1\left( 1,0 \right)  \Big) & \text{if both $x$ and $p$ are } \\ \vee \Big( \beta^{-2} \K_{G_2,p}\left( \N_2  \right)+ \beta^2\mathcal{Q}_2\left( 0,1 \right)  \Big) & \text{ weakly curvature saturated,} \\[5pt]  \alpha^{-2} \K_{G_1,x}\left( \N_1  \right) + \alpha^2 \mathcal{Q}_1\left( 1,0 \right)  & \text{if $x$ is weakly curvature saturated and} \\ + \; 2\N_1^{-1}\N_2\left(  \N_1 + \N_2  \right)^{-1} \D_x & \text{ $p$ is curvature un-saturated,} \\[10pt] \beta^{-2} \K_{G_2,p}\left( \N_2  \right) +  \beta^2\mathcal{Q}_2\left( 0,1 \right)  & \text{$x$ is curvature un-saturated and} \\  +\; 2\N_1\N_2^{-1}\left(  \N_1 + \N_2  \right)^{-1} \D_p & \text{ if $p$ is weakly curvature saturated,} \\[10pt] 
	\Big( \alpha^{-2} \K_{G_1,x}\left( \N_1  \right) + \alpha^2 \mathcal{Q}_1\left( 1,1 \right)  \Big) & \text{if neither $x$ nor $p$ is} \\ \vee \Big( \beta^{-2} \K_{G_2,p}\left( \N_2  \right)+ \beta^2\mathcal{Q}_2\left( 1,1 \right)  \Big) & \text{ strongly curvature saturated}.
	\end{cases}
	\end{align*}
When both $x$ and $p$ are weakly curvature saturated but neither is strongly curvature saturated (this is a subcase of case 1 in the above theorem), we also deduce (this time with minimum instead of max)
\begin{align}
\K_{(x,p)}\left( \N_1 + \N_2 \right) &\le \Big( \alpha^{-2} \K_{G_1,x}\left( \N_1  \right) + \alpha^2 \mathcal{Q}_1\left( 1,0 \right) + 2\N_1^{-1}\N_2\left(  \N_1 + \N_2  \right)^{-1} \D_x \Big) \notag \\ & \wedge \left(  \beta^{-2} \K_{G_2,p}\left( \N_2  \right)+ \beta^2\mathcal{Q}_2\left( 0,1 \right) + 2\N_1\N_2^{-1}\left(  \N_1 + \N_2  \right)^{-1} \D_p \right), \notag 
\end{align}
where $\mathcal{Q}_1$ and $\mathcal{Q}_2$ are (piece-wise) quadratic forms given by 
\begin{align*}
\mathcal{Q}_1\left( c_1,c_2 \right) :=  \nicefrac{1}{2}\, c_1^2 \beta^{-2} \Delta^{G_2} \alpha^{-2} + \left| c_1c_2\right| \Big( \beta^{-2} \D_x + \nicefrac{1}{2}\, \alpha^{-2}\Gamma^{G_1} \left( \beta^{-2}  \right) \Big),
\end{align*}
and
\begin{align*}
\mathcal{Q}_2\left( c_1,c_2 \right) := \nicefrac{1}{2}\, c_2^2 \alpha^{-2} \Delta^{G_1} \beta^{-2} + \left| c_1c_2\right|  \Big( \alpha^{-2} \D_p +\nicefrac{1}{2}\, \beta^{-2}\Gamma^{G_2} \left( \alpha^{-2}  \right) \Big).
\end{align*}
\end{theorem}
Combining the above upper bounds one deduces the following. 
\begin{corollary}[An all inclusive non iso-dimensional upper bound]\label{cor:non-iso-upb}
	The upper bound
	\begin{align}
		\K_{(x,p)}\left( \N_1 + \N_2 \right) &\le \Big( \alpha^{-2} \K_{G_1,x}\left( \N_1  \right) + \alpha^2 \mathcal{Q}_1\left( 1,1 \right) + 2\N_1^{-1}\N_2\left(  \N_1 + \N_2  \right)^{-1} \D_x \Big) \notag \\ & \vee \left(  \beta^{-2} \K_{G_2,p}\left( \N_2  \right)+ \beta^2\mathcal{Q}_2\left( 1,1 \right) + 2\N_1\N_2^{-1}\left(  \N_1 + \N_2  \right)^{-1} \D_p \right), \notag 
	\end{align}
always holds. 
\end{corollary}
The following is an immediate consequence of monotonicity of curvature functions in $\N$.
\begin{corollary}\label{cor:iso-dim-up}
	Letting $\N_1 = \N_2 = \N$, we get
\begin{align}
\K_{(x,p)}\left( \N \right) \le	\K_{(x,p)}\left( 2\N \right) &\le \Big( \alpha^{-2} \K_{G_1,x}\left( \N  \right) + \nicefrac{1}{2} \, \alpha^2 \beta^{-2} \Delta^{G_2} \alpha^{-2} + \N^{-1} \D_x \Big) \notag \\ & \wedge \Big(  \beta^{-2} \K_{G_2,p}\left( \N  \right)+ \nicefrac{1}{2}\, \beta^2 \alpha^{-2} \Delta^{G_1} \beta^{-2} + \N^{-1}\D_p \Big), \notag 
\end{align}	
and similarly, we get the following cases:
\begin{enumerate}
\item if both $x$ and $p$ are weakly curvature saturated, then
\begin{align*}
		\K_{(x,p)}\left( \N \right) \le \K_{(x,p)}\left( 2\N \right) \le \Big(\alpha^{-2} \K_{G_1,x}\left( \N  \right) + \alpha^2 \mathcal{Q}_1\left( 1,0 \right)  \Big)  \vee \Big( \beta^{-2} \K_{G_2,p}\left( \N  \right)+ \beta^2\mathcal{Q}_2\left( 0,1 \right)  \Big);
\end{align*}
\item 	if $x$ is weakly curvature saturated and $p$ is curvature un-saturated, then
	\begin{align*}
		\K_{(x,p)}\left( \N \right) \le \K_{(x,p)}\left( 2\N \right) \le \alpha^{-2} \K_{G_1,x}\left( \N  \right) + \alpha^2 \mathcal{Q}_1\left( 1,0 \right) + \; \N^{-1} \D_x ;
	\end{align*}
\item 	if $p$ is weakly curvature saturated and $x$ is curvature un-saturated, then
	\begin{align*}
	 	\K_{(x,p)}\left( \N \right) \le \K_{(x,p)}\left( 2\N \right) \le	\beta^{-2} \K_{G_2,p}\left( \N  \right) +  \beta^2\mathcal{Q}_2\left( 0,1 \right) +\N^{-1} \D_p;
	\end{align*}
\item 	if neither $x$ nor $p$ is strongly curvature saturated
	\begin{align*}
			\K_{(x,p)}\left( \N \right) \le \K_{(x,p)}\left( 2\N \right) \le \Big( \alpha^{-2} \K_{G_1,x}\left( \N  \right) + \alpha^2 \mathcal{Q}_1\left( 1,1 \right)  \Big)  \vee \Big( \beta^{-2} \K_{G_2,p}\left( \N  \right)+ \beta^2\mathcal{Q}_2\left( 1,1 \right)  \Big).
	\end{align*}
\end{enumerate}
\end{corollary}
Letting $\N \to \infty$ in Corollary~\ref{cor:iso-dim-up}, we get the following.
\begin{corollary}\label{cor:iso-dim-up-infty}
	\begin{align}
		\K_{(x,p)}\left( \infty \right) \le \Big( \alpha^{-2} \K_{G_1,x}\left( \infty \right) + \nicefrac{1}{2} \, \alpha^2 \beta^{-2} \Delta^{G_2} \alpha^{-2}  \Big) \wedge \Big(  \beta^{-2} \K_{G_2,p}\left( \infty \right)+ \nicefrac{1}{2}\, \beta^2 \alpha^{-2} \Delta^{G_1} \beta^{-2} \Big). \notag 
	\end{align}	
\end{corollary}
By Theorem~\ref{thm:main-5}, at vertices where $\alpha$ and $\beta$ are sufficiently convex in the sense of (\ref{eq:g-conv-1}) and (\ref{eq:g-conv-2})  (see also the Remark after Corollary~\ref{cor:one-thirteen}), we get the following neater estimates.
\begin{corollary}
	For any $K_1, K_2 \in \R$ . For every $x \in G_1$ and $p \in G_2$, 
	if
	\begin{align}\label{eq:g-conv-1}
	\alpha^2 \Delta^{G_2} \alpha^{-2} \le - 2\alpha^{-2}\beta^2 K_1 - 2 \alpha^2 \D_x - \beta^2 \Gamma^{G_1} \left( \beta^{-2}  \right),
	\end{align}
	and
	\begin{align}\label{eq:g-conv-2}
	\beta^2 \Delta^{G_2} \alpha^{2}\beta^{-2} \le  - 2\alpha^2 K_2 - 2 \beta^2 \D_p - \alpha^2 \Gamma^{G_2} \left( \alpha^{-2}  \right), 
	\end{align}
hold then
\begin{align*}
	\K_{(x,p)}\left( \N_1 + \N_2 \right) \le \frac{\K_{G_1,x}\left( \N_1\right)  - K_1}{\alpha^2}  \vee \frac{\K_{G_2,p}\left( \N_2\right) - K_2}{\beta^2}.
\end{align*}
\end{corollary}
\begin{proof}
	It is straightforward to see that \eqref{eq:g-conv-1} and \eqref{eq:g-conv-2} imply 
\begin{align*}
	\mathcal{Q}_1\left( 1,1\right) \le -\alpha^{-2}K_1 \quad \text{and} \quad 	\mathcal{Q}_2\left( 1,1\right) \le -\beta^{-2}K_2,
\end{align*}
hence, the conclusion follows from the all cases inclusive non iso-dimensional upper bounds from Corollary~\ref{cor:non-iso-upb}. 
\end{proof}
\begin{theorem}[All inclusive iso-dimensional upper bound]\label{thm:main-6}
	The curvature function of a doubly warped product always satisfies  
\begin{align*}
	\mathcal{K}_{(x,p)}(\N) \le \Big( \alpha^{-2} \mathcal{K}_{G_1,x}(\N) + \nicefrac{1}{2}\, \alpha^{2}\beta^{-2} \Delta^{G_2}  \alpha^{-2}\Big) \wedge \Big(\beta^{-2} \mathcal{K}_{G_2,p}(\N) + \nicefrac{1}{2}\, \beta^2\alpha^{-2} \Delta^{G_1}  \beta^{-2} \Big).
\end{align*}
\end{theorem}
\begin{corollary}\label{cor:one-thirteen}
	For any $K_1,K_2 \in \R$, if
	\begin{align}\label{eq:g-conv-3}
		\alpha^2 \Delta^{G_2} \alpha^{-2} \le - 2\alpha^{-2}\beta^2K_1,
	\end{align}
	and
	\begin{align}\label{eq:g-conv-4}
		\beta^2 \Delta^{G_1} \beta^{-2} \le - 2\beta^{-2}\alpha^2K_2,
	\end{align}
	hold then
\begin{align*}
	\K_{(x,p)}\left( \N \right) \le \frac{\K_{G_1,x}\left( \N\right)  -  K_1}{\alpha^2}  \wedge \frac{\K_{G_2,p}\left( \N\right)  -  K_2}{\beta^2}.
\end{align*}
\end{corollary}
\begin{corollary}
	\eqref{eq:g-conv-3} and \eqref{eq:g-conv-4} imply 
\begin{align*}
	\nicefrac{1}{2}\, \alpha^{2}\beta^{-2} \Delta^{G_2}  \alpha^{-2} \le -\alpha^{-2}K_1 \quad \text{and} \quad  \nicefrac{1}{2}\, \beta^2\alpha^{-2} \Delta^{G_1}  \beta^{-2}  \le -\beta^{-2}K_2
\end{align*}
hence the conclusion follows from Theorem~\ref{thm:main-6}. 
\end{corollary}
\begin{remark}[convexity conditions in discrete setting]\label{rem:convexity}
	Notice in the Riemannian setting where the chain rule is available, we have the identity 
\begin{align*}
	\alpha^2 \Delta^{G_2} \alpha^{-2} = -2 \Delta a - 4 \|\nabla a\|^2.
\end{align*}
	Therefore,  conditions~(\ref{eq:g-conv-1})-(\ref{eq:g-conv-4}) above could be thought of as discrete counterparts of dynamic \emph{convexity} conditions on $a = \ln \alpha$ and $b = \ln \beta$.   By dynamic, we mean vertex dependent. 
\end{remark}
\subsubsection*{\small \bf \textit{\textsf{Lower bounds on the curvature function}}} 
\begin{theorem}[Non iso-dimensional lower bounds]\label{thm:main-lb-1}
	For any $\M > \max \left\{  \N_1, \N_2 \right\}$, the lower bounds 
\begin{enumerate}
		\item  \begin{align*}
		\begin{split}
		&\K_{G, (x,p)} (\M) \ge \\
		&	\Big(   \alpha^{-2} \K_{G_1,x}\left( \N_1  \right)   - \alpha^{2} \beta^{-2}\big(\N_2^{-1} - \M^{-1}  \big)^{-1}  \Gamma^{G_1} \left(  \beta^{-2}\right)  \Big. \\ &  \Big. \phantom{sajjadsajjadsajja} + \nicefrac{1}{2}\, \alpha^{2} \beta^{-2} \Delta^{G_2}\alpha^{-2} - 2\M^{-1}\alpha^{-2}\Deg_i  \Big) \\ &  \wedge  \Big(   \beta^{-2} \K_{G_2,p}\left( \N_2  \right)  -\alpha^{-2} \beta^{2}\big(\N_1^{-1} - \M^{-1}  \big)^{-1}  \Gamma^{G_2} \left(  \alpha^{-2}\right) \Big. \\ & \Big.   \phantom{sajjadsajjadsajja} + \nicefrac{1}{2}\, \alpha^{-2}\beta^{2} \Delta^{G_1}\beta^{-2}  - 2\M^{-1}\beta^{-2}\Deg_2     \Big),
		\end{split}
	\end{align*}
		\item[] and \medskip
	\item[(2)]
		\begin{align*}
			&\K_{G, (x,p)} (2\M) \ge \\
			&  \Big(  \alpha^{-2} \K_{G_1,x}\left( \N_1  \right)   -\alpha^{2} \beta^{-2}\big(\N_2^{-1} - \M^{-1}  \big)^{-1} \Gamma^{G_1} \left(  \beta^{-2}\right)+ \nicefrac{1}{2}\, \alpha^{2} \beta^{-2} \Delta^{G_2}\alpha^{-2} \Big)   \\ & \wedge \Big(  \beta^{-2} \K_{G_2,p}\left( \N_2  \right)  -\alpha^{-2} \beta^{2}\big(\N_1^{-1} - \M^{-1}  \big)^{-1}  \Gamma^{G_2} \left(  \alpha^{-2}\right) + \nicefrac{1}{2}\, \alpha^{-2}\beta^{2} \Delta^{G_1}\beta^{-2}    \Big),
		\end{align*}
	\end{enumerate}
hold on the curvature function. 
\end{theorem}
\begin{corollary}\label{cor:main-lb}
The curvature function satisfies	
	\begin{enumerate}
		\item []
		\begin{align*}
			\begin{split}
				&\K_{G, (x,p)} \left( \N_1 + \N_2 \right)  \ge \\
				&\Big(   \alpha^{-2} \K_{G_1,x}\left( \N_1  \right)   - \alpha^{2} \beta^{-2} \left( \N_1 + \N_2 \right)\N_1^{-1}\N_2 \Gamma^{G_1} \left(  \beta^{-2}\right) \Big. \\ &  \Big. \phantom{sajjadsajjadsajja}  + \nicefrac{1}{2}\, \alpha^{2} \beta^{-2} \Delta^{G_2}\alpha^{-2} - 2 \left( \N_1 + \N_2 \right)^{-1}\alpha^{-2}\Deg_i  \Big) \\ &  \wedge  \Big(  \beta^{-2} \K_{G_2,p}\left( \N_2  \right)  -\alpha^{-2} \beta^{2} \left( \N_1 + \N_2 \right)\N_1\N_2^{-1}  \Gamma^{G_2} \left(  \alpha^{-2}\right) \Big. \\ & \Big.  \phantom{sajjadsajjadsajj}  + \nicefrac{1}{2}\, \alpha^{-2}\beta^{2} \Delta^{G_1}\beta^{-2}  - 2 \left( \N_1 + \N_2 \right) ^{-1}\beta^{-2}\Deg_2    \Big),
			\end{split}
		\end{align*} 
			\item[] and \medskip
			\item []\begin{align*}
		&	\K_{G, (x,p)} \left( 2\left( \N_1 + \N_2 \right) \right) \ge\\
		  &\Big(  \alpha^{-2} \K_{G_1,x}\left( \N_1  \right)   -\alpha^{2} \beta^{-2} \left( \N_1 + \N_2 \right)\N_1^{-1}\N_2  \Gamma^{G_1} \left(  \beta^{-2}\right)+ \nicefrac{1}{2}\, \alpha^{2} \beta^{-2} \Delta^{G_2}\alpha^{-2} \Big)   \\ & \wedge \Big(  \beta^{-2} \K_{G_2,p}\left( \N_2  \right)  -\alpha^{-2} \beta^{2} \left( \N_1 + \N_2 \right)\N_1\N_2^{-1} \Gamma^{G_2} \left(  \alpha^{-2}\right) + \nicefrac{1}{2}\, \alpha^{-2}\beta^{2} \Delta^{G_1}\beta^{-2}    \Big).
		\end{align*}
	\end{enumerate}
In particular,
	\begin{enumerate}
	\item []
	\begin{align}\label{eq:2n-lb}
		\begin{split}
			&\K_{G, (x,p)} \left( 2\N \right)  \ge \\ 
			&\Big(   \alpha^{-2} \K_{G_1,x}\left( \N  \right)   - \alpha^{2} \beta^{-2} \left( 2\N \right) \Gamma^{G_1} \left(  \beta^{-2}\right) + \nicefrac{1}{2}\, \alpha^{2} \beta^{-2} \Delta^{G_2}\alpha^{-2} - 2 \N^{-1}\alpha^{-2}\Deg_i  \Big) \\ &  \wedge  \Big(  \beta^{-2} \K_{G_2,p}\left( \N \right)  -\alpha^{-2} \beta^{2} \left( 2\N \right)  \Gamma^{G_2} \left(  \alpha^{-2}\right) + \nicefrac{1}{2}\, \alpha^{-2}\beta^{2} \Delta^{G_1}\beta^{-2}  - 2 \N^{-1}\beta^{-2}\Deg_2    \Big),
		\end{split}
	\end{align} 
	\item[] and \medskip
	\item []\begin{align*}
		&	\K_{G, (x,p)} \left( 4\N \right) \ge\\
		&\Big(  \alpha^{-2} \K_{G_1,x}\left( \N  \right)   -\alpha^{2} \beta^{-2} \left( 2\N \right)  \Gamma^{G_1} \left(  \beta^{-2}\right)+ \nicefrac{1}{2}\, \alpha^{2} \beta^{-2} \Delta^{G_2}\alpha^{-2} \Big)   \\ & \wedge \Big(  \beta^{-2} \K_{G_2,p}\left( \N \right)  -\alpha^{-2} \beta^{2} \left( 2\N \right) \Gamma^{G_2} \left(  \alpha^{-2}\right) + \nicefrac{1}{2}\, \alpha^{-2}\beta^{2} \Delta^{G_1}\beta^{-2}    \Big).
	\end{align*}
\end{enumerate}
\end{corollary}
\begin{corollary}[Iso-dimensional lower bound]\label{cor:iso-dim-lb}
	\begin{align*}
		\begin{split}
		&	\mathcal{K}_{(x,p)}(\N) \ge  \\
			&\Big(   \alpha^{-2} \K_{G_1,x}\left( \N  \right)   - \alpha^{2} \beta^{-2} \left( 2\N \right) \Gamma^{G_1} \left(  \beta^{-2}\right) + \nicefrac{1}{2}\, \alpha^{2} \beta^{-2} \Delta^{G_2}\alpha^{-2} - 2 \N^{-1}\alpha^{-2}\Deg_i  \Big) \\ &  \wedge  \Big(  \beta^{-2} \K_{G_2,p}\left( \N \right)  -\alpha^{-2} \beta^{2} \left( 2\N \right)  \Gamma^{G_2} \left(  \alpha^{-2}\right) + \nicefrac{1}{2}\, \alpha^{-2}\beta^{2} \Delta^{G_1}\beta^{-2}  - 2 \N^{-1}\beta^{-2}\Deg_2    \Big) \\ & - \N^{-1}\left( \alpha^{-2}\D_x + \beta^{-2}\D_p  \right).
		\end{split} 
	\end{align*} 
\end{corollary}
Keeping $\N_i$ fixed and letting $\M \to \infty$ in either of the above bounds, we get
\begin{theorem}[Lower bounds for dimensionless curvature function]\label{thm:lb-infty}
	\begin{align*}
	\begin{split}
		&\K_{G, (x,p)} (\infty) \ge  \\
		&\Big(   \alpha^{-2} \K_{G_1,x}\left( \N_1  \right)   - \alpha^{2} \beta^{-2}\N_2  \Gamma^{G_1} \left(  \beta^{-2}\right) + \nicefrac{1}{2}\, \alpha^{2} \beta^{-2} \Delta^{G_2}\alpha^{-2} \Big) \\ &  \wedge  \Big(   \beta^{-2} \K_{G_2,p}\left( \N_2  \right)  -\alpha^{-2} \beta^{2}\N_1  \Gamma^{G_2} \left(  \alpha^{-2}\right) + \nicefrac{1}{2}\, \alpha^{-2}\beta^{2} \Delta^{G_1}\beta^{-2}     \Big).
	\end{split}
\end{align*}
\end{theorem}
\begin{corollary}[Special case]\label{cor:lb-infty}
	If $\alpha$ and $\beta$ are constants in balls of radius $1$ around $p$ and $x$ resp. (i.e. they are locally constant at $p$ and $x$) then, 
	\begin{align*}
			\K_{G, (x,p)} (\infty) = \alpha^{-2} \K_{G_1,x}\left( \infty  \right)    \wedge    \beta^{-2} \K_{G_2,p}\left( \infty  \right).
	\end{align*}
\end{corollary}
\begin{proof}
Based on Corollary~\ref{cor:iso-dim-up-infty} and Theorem~\ref{thm:lb-infty}, we have
\begin{align*}
\alpha^{-2} \K_{G_1,x}\left( \N_1  \right)   \wedge    \beta^{-2} \K_{G_2,p}\left( \N_2  \right)   \le	\K_{G, (x,p)} (\infty) \le \alpha^{-2} \K_{G_1,x}\left( \infty  \right)    \wedge    \beta^{-2} \K_{G_2,p}\left( \infty  \right).
\end{align*}
Upon letting $\N_i \to \infty$, we get the desired result. 
\end{proof}
\subsection*{Acknowledgements}
\par \small \textit{SL acknowledges partial support from a ``resident researcher'' grant from IPM (Grant No.1400460424).}\medskip

	\small	\textit{The authors are grateful to Radoslaw Wojciechowski for many insightful conversations about the geometry of graphs and to Florentin M\"unch for helpful correspondence about graph curvature. We are greatly indebted to the anonymous referees of this manuscript, for their time and for providing us with their positive valuable feedback which played an imperative role in improving these notes.}
\normalsize
\section{Preliminaries}
\label{SEC:PRELIM}
\subsection{ Smooth setting}\label{sec:prelim-smooth}
The curvature properties of (doubly) twisted and warped products of Riemannian manifolds have been studied by various authors; e.g.~\cite{FGKU,Un, PR,EJK,EJKS,Ge,FS,DL,Choi}. We start off by discussing the curvature bounds for a generalized doubly warped product of weighted manifolds. 
\par Let $\left( M^n, g, e^{-\Phi}\dvol_g \right)$ be a complete weighted manifold. In the interior of $M$ the corresponding drift Laplacian is defined by
\begin{align*}
\Delta_\Phi = \Delta - \nabla\Phi \cdot \nabla.
\end{align*}
For $\mathcal{N}\ge n$, the $\mathcal{N}-$Bakry-\'Emery Ricci tensor is then given by
\begin{align*}
\Ric_\Phi^{\mathcal{N}} = \Ric + \nabla^2 \Phi - \left(   \mathcal{N} - n\right)^{-1}\nabla\Phi \otimes \nabla \Phi,
\end{align*}
with the conventions $\Ric_\Phi^{\infty} = \Ric + \nabla^2 \Phi$ and $\Ric_\Phi^{n} = \Ric$ (this requires $\Phi$ to be constant). When $\Ric_\Phi^{\mathcal{N}} \ge \K g$, we say the weighted manifold satisfies $BE\left( \K , \mathcal{N}  \right)$ curvature dimension conditions. Considering the square field operator $\Gamma$ defined in~\cite{BE} via
\begin{align}\label{eq:Gamma-smooth}
\Gamma (u,v) := \nicefrac{1}{2}\left( \Delta_\Phi uv - u \Delta_\Phi v  - v \Delta_\Phi u   \right) = \nabla u \cdot \nabla v,
\end{align}
and the iterated $\Gamma_2$ operator given by
\begin{align*}
\Gamma_2 \left(u \right) := \nicefrac{1}{2}\, \Delta_\Phi \Gamma (u) - \Gamma\left( \Delta_\Phi u , u  \right)  = \nicefrac{1}{2}\, \Delta_\Phi \left|  \nabla u  \right|^2 - \nabla u \cdot \nabla \Delta_\Phi u,
\end{align*}
the celebrated Bochner formula can be rewritten as
\begin{align*}
\Gamma_2 \left(f \right) - \Ric^\N_\Phi\left( \nabla f, \nabla f  \right) = \left| \nabla^2 u  \right|^2  +  \left(   \mathcal{N} - n\right)^{-1}\nabla\Phi \otimes \nabla \Phi;
\end{align*}
e.g. \cite[P. 397]{Villani}. The Bochner formula results in
\begin{align*}
\Gamma_2 \left(u \right) - \Ric^\N_\Phi\left( \nabla u, \nabla u  \right) \ge \N^{-1}\left(  \Delta_\Phi u  \right)^2;
\end{align*}
see~\cite{Le}. Therefore, $ \Ric_\Phi^\N \ge \K g$ implies
\begin{align*}
\Gamma_2\left(  u \right) \ge \N^{-1}\left( \Delta_\Phi u \right)^2 + \K\Gamma (u) \quad\quad  \forall u;
\end{align*}
the latter is referred to as $CD(\K,\N)$ curvature dimension condition for the diffusion operator $\Delta_\Phi$; also sometimes called $BE(\K,\N)$ conditions where ``BE'' stands for Bakry-\'Emery. Conversely, if this holds for all smooth functions $u$, by taking curvature maximizers, one deduces $\Ric_\Phi^\N \ge \K g$; see~\cite{St1} for the proof of the above facts in a much more general setting. \emph{We note the difference between the above mentioned $CD(\K,\N)$ curvature dimension conditions and the Lot-Villani-Sturm curvature dimension conditions which are also referred to, in the literature, as  $CD(\K,\N)$ conditions.}
\subsection{ Discrete setting}\label{sec:disc-curv}
In recent years, there has been a substantial interest in curvature of discrete structures, one for the fact that the definitions are simple enough to be programmable and robust enough to determine the geometry. Notions of curvature of graphs started to appear in the literature in as early as the 70's and 80's with~\cite{Sto,DK} and and later on in~\cite{CY,For,Sch}. After Lott, Sturm and Villani's breakthrough in the seminal papers~\cite{St2,St3,LV}, where they developed weak Ricci curvature lower bounds for a broad class of metric spaces, there has been a sudden surge of research in understanding the curvature of discrete structures using methods of optimal transport as in~\cite{BS,Ol,LLY,EM,MW} and using the $\Gamma_2$ calculus methods as had been previously developed in the smooth setting in~\cite{BE,BL}; see e.g.~\cite{LY,JL,CLY,HL,JL,CLP}. Other versions can be found in e.g.~\cite{BHLLMY,Mu,LM}. We also point out to the papers~\cite{EK,GM} that provide some discrete to continuous picture of Wassestein spaces and (dynamic) curvature bounds (super Ricci flows). The literature is too extensive to be covered here so to do justice, we encourage the interested reader to look at the above papers and references therein. 
\par For us, an un-directed weighted graph $G$ is a nonnegative (not-necessarily symmetric) weight function $\omega:\mathbb{Z}^2 \to \R$ satisfying the transition relations $\omega(x,y)=s(x,y)\omega(y,x)$ for $s(x,y) \neq 0$. The vertex and edge sets are respectively identified by
\begin{align*}
V := \{x \in \Z : \exists y \in \Z \; , \;  \omega(x,y)>0 \},
\end{align*}
and
\begin{align*}
E := \{(x,y) \in \Z^2 : \omega(x,y)>0    \}/ (x,y)\sim (y,x).
\end{align*}
\par Finite graphs are given by finitely supported weight functions $\omega$. We write $x\sim_G y$ or $x\sim y$ when there is an edge between $x$ and $y$. We set $\omega_{xy} := \omega(x,y)$ and $\omega_{xy} = s(x,y)\omega_{yx}$ for a nonzero $s(x,y)$. The vertex measure is a function $m:V \to \R_+$. For simplicity, we sometimes write $m_x$ instead of $m(x)$. $G$ will both denote a weighted graph and its vertex set. For any vertex $x$, we set
\begin{align*}
\Deg_G(x) := m(x)^{-1} \sum_{y \sim x} \omega_{xy},
\end{align*}
which will be abbreviated as $\mathrm{D}_x$ when there is no ambiguity in determining the underlying graph. We consider $G$ to be equipped with the most general Laplacian of the form
\begin{align}\label{eq:Laplacian-g}
\Delta f (x):= m(x)^{-1} \sum_{ y \sim x} \left( f(y) - f(x)   \right) \omega_{xy}.
\end{align}
\par The corresponding square field operator $\Gamma$ and the \emph{Ricci form}, $\Gamma_2$ (which is the iterated square field operator) are thus given by
\begin{align}\label{eq:Gamma1-g}
\Gamma \left( u,v \right) =  \nicefrac{1}{2} \left( \Delta (uv) - v\Delta u - u\Delta v \right),
\end{align}
and
\begin{align}\label{eq:Gamma2-g}
\Gamma_2 \left( u,v \right) = \nicefrac{1}{2}\left(  \Delta \Gamma \left( u,v \right) - \Gamma \left( \Delta u , v  \right) - \Gamma \left( u , \Delta v  \right) \right),
\end{align}
respectively. 
\par Analogous to the smooth setting and the Bakry-\'Emery curvature dimension conditions $BE(\K,\N)$, the discrete Bakry-\'Emery curvature dimension conditions $CD\left(  \K,\N \right)$ at a vertex $x \in G$ amounts to the inequality
\begin{align*}
\Gamma_2(f)(x) \ge \N^{-1}\left( \Delta f  \right)^2(x) + \K \Gamma (f)(x) \quad \quad \forall f: G \to \R,
\end{align*}
where $\Gamma(f) := \Gamma(f,f)$ and $\Gamma_2(f) := \Gamma_2 (f,f)$. When this inequality holds globally, we say the graph $G$ satisfies the (discrete) $CD(\K,\mathcal{N})$ curvature dimension conditions. 
\begin{remark}
	\par Before we proceed, we note that the notation $CD\left(  \K,\N \right)$ for discrete Bakry-\'Emery curvature dimension conditions that we use is a well established notation in the context of graph curvature; e.g. see~\cite{CLP}.  This should not be confused with neither the Lott-Sturm-Villani curvature dimension conditions (denoted by the same notation) which are defined on metric measure spaces using optimal transportation, nor with the discrete generalizations defined via discrete optimal transportation. Throughout the rest of this article, a curvature bound for a graph at a vertex means a Bakry-\'Emery curvature dimension bound at that vertex (a number $\K$ such that $CD\left(  \K,\N \right)$ conditions hold). 
\end{remark}
\par For a given $\N$, the best such lower curvature bound at a vertex $x$, will be denoted by $\K_{G,x}(\N)$. It follows from the definitions of these operators that $\Delta$ and $\Gamma$ are linear and $\Gamma_2$ is a quadratic form in terms of the weights $\omega_{xy}$. Hence, setting 
\begin{align*}
G_\lambda := \left( G, \lambda \omega, m  \right),
\end{align*}
we have
\begin{align*}
\K_{G_\lambda,x}(\N) = \lambda \K_{G,x}(\N) \quad \quad \forall \N>0. 
\end{align*}
\par In practice, using the structural curvature dimension bounds obtained in Theorem~\ref{thm:main-3}, one can directly obtain estimates on curvature dimension bounds for doubly twisted products of weighted graphs and networks. However, we take a different approach and instead explore the Ricci form of the doubly warped product to deduce neater bounds in terms of the curvature bounds of the factors. The bounds obtained bear resemblance to the Riemannian curvature bounds. We should mention that the curvature dimension bounds for the un-normalized discrete Laplacian operator in Cartesian products of graphs have been studied in~\cite{LP,CLP}. Important properties of the curvature functions such as their monotonicity and concavity have been discussed in \cite{CLP}. The main difficulty in working with graphs is the lack of chain rule which is hand in hand with the fact that discrete Laplacian is almost never a diffusion operator; so we do not have the chain rule at our disposal.
\section{Proof of main theorems in the smooth setting}
We have introduced a weighted doubly warped product which generalizes the existing notions and is one of the novel aspects of this work. In this section, we calculate its Bakry-\'Emery Ricci curvature bounds in terms of those of the factors.
\par It needs to be emphasized that the curvature estimates are mostly pointwise and they hold wherever the given partial differential inequalities hold at the point therefore, we are not apriori fixing any signs for the constants involved except for the weights $R_1$ and $R_2$ which have to be nonnegative. Obtaining bounds for negative $R_i$ can be done along the same lines but we will not pursue that in these notes. A posteriori, the partial differential inequalities (our so called dynamic concavity relations) will impose restrictions on the size and signs of the constants involved. 
\subsection{ Ricci tensor for weighted doubly warped products}\label{sec:calc-stuff}
Let $\left( B^{n_1}, g_B  \right)$ and $\left( F^{n_2}, g_F  \right)$ be two Riemannian manifolds. Let $\alpha:F \to \R_+$ and $\beta:B \to \R_+$ be smooth positive warping functions. Throughout these notes, we will use the conventions $a := \ln \alpha$ and $b:= \ln \beta$. For the doubly warped product
\begin{align*}
B {_\alpha\times_\beta} F := \left( B\times F , g := \alpha^2g_B \oplus \beta^2 g_F   \right),
\end{align*}
\subsubsection*{\bf \textit{$\nabla$-calculus}}
By the aid of Koszul's formula, it is straightforward task to see that the covariant derivatives are given by
\begin{align*}
\quad \nabla_XY =    {^B\nabla}_XY   - \left< X,Y  \right>\nabla a, \quad   \nabla_VW =  {^F\nabla}_VW   - \left< V,W  \right>\nabla b,
\end{align*}
and
\begin{align*}
\nabla_X V = \nabla_V X = \left(\nabla_X b \right)V + \left( \nabla_V a \right)X;
\end{align*}
see \cite{ONeil} for details. 
\par Now suppose $\Phi: B \to \R$ is a smooth function. We also denote by $\Phi$, the lift of $\Phi$ to a function from $B \times F$ (i.e. the lift $\Phi \circ \mathrm{proj}_B$). We can then compute the Hessian of $\Phi$ by direct calculation. 
\begin{align*}
	\nabla^2 \Phi\left( X+V,X+V \right) &= \nabla_{X+V}\nabla_{X+V}\Phi -   \left(\nabla_{X+V} (X+V) \right)\Phi \\ &= {^B\nabla}_X{^B\nabla}_X  \Phi - d\Phi \left( {^B\nabla}_XX  + 2 \left(\nabla_Va\right) X + 2\left(\nabla_X b \right)V - \|V\|^2 \nabla b \right) \notag \\ &= {^B\nabla}^2 \Phi (X,X) - 2 \nabla_V a \nabla_X \Phi  + \|V\|^2 \left< \nabla \Phi , \nabla b  \right>.
\end{align*}
Similarly for a smooth function $\Psi: F\to \R$ (the lift thereof to the product space $B \times F$ given by $\Psi \circ \mathrm{proj}_F$), we get
\begin{align*}
	\nabla^2 \Psi\left( X+V,X+V \right) &= \nabla_{X+V}\nabla_{X+V}\Psi -   \left(\nabla_{X+V} (X+V) \right)\Psi \\ &= {^F\nabla}_V {^F\nabla}_V \Psi - d\Psi \left( {^F\nabla}_VV  + 2 \left(\nabla_Xb\right) V + 2\left(\nabla_V a \right)X - \|X\|^2 \nabla a \right) \notag \\ &={^F\nabla}^2 \Psi \left( V,V \right) - 2 \nabla_X b \nabla_V \Psi  + \|X\|^2 \left< \nabla \Psi , \nabla a  \right>.
\end{align*}
In particular,
\begin{align*}
	\nabla^2 b \left( X , Y  \right) = {^B\nabla}^2b\left( X,Y   \right), \quad \nabla^2 b \left( V , W  \right) = - \left< V, W\right> \|\nabla b\|^2,
\end{align*}
and similarly
\begin{align*}
	\nabla^2 a \left( V , W  \right) = {^F\nabla}^2a\left( V,W   \right), \quad \nabla^2 a \left( X , Y  \right) = - \left< X , Y\right> \|\nabla a\|^2.
\end{align*}
Now notice that from an orthonormal frame $\{X_i\}^{n_1}_{i=1}$ for $g_B$ and $\{V_j\}^{n_2}_{j=1}$ for $g_F$, one gets the orthonormal frame consisting of $\tilde{X}_i := \frac{1}{\alpha}\{X_i\}^{n_1}_{i=1}$ and $\tilde{V}_i := \frac{1}{\beta}\{V_j\}^{n_2}_{j=1}$ for $g$. Consequently, tracing over the orthonormal frame $\left\{ \tilde{X}, \tilde{V} \right\}$, we get
\begin{align*}
	\Delta b = \mathrm{tr}\nabla^2 b =  \frac{1}{\alpha^2}\,  {^B\Delta} b - n_2\|\nabla b\|^2, \quad 	\Delta a = \frac{1}{\beta^2}\, {^F\Delta} a - n_1\|\nabla a\|^2,
\end{align*}
where
\begin{align*}
	\|\nabla b\|^2 = db\left( \grad(b) \right) = db\left( \frac{1}{\alpha^2}  {^B\grad}(b)\right) = \frac{1}{\alpha^2} \|{^B\nabla} b\|_B^2,
\end{align*}
and similarly, 
\begin{align*}
	\|\nabla a\|^2 = \frac{1}{\beta^2} \|{^F\nabla} a\|_F^2.
\end{align*}

\par Finally, to avoid confusion we note that, by using the orthonormal frames, the norms of $\nabla b$ as a $1$-tensor, are related via
\begin{align*}
	\|\nabla b\|^2 = \sum \left(\nabla_{\tilde{X}_i} b\right)^2 = \alpha^{-2} \sum \left(\nabla_{X_i} b\right)^2 = \alpha^{-2} \|{\nabla} b\|_B^2\quad \text{and}\quad 	\|\nabla a\|^2 = \beta^{-2}\|{\nabla} a\|_F^2
\end{align*}
while obviously the norms of the gradient of $a$ and $b$ (gradient w.r.t. the metric $g$), are related via
\begin{align*}
\| \grad b\|^2 = \alpha^{2} \| \grad b \|_B^2 \quad \text{and}\quad \| \grad a \|^2 = \beta^{2} \| \grad a \|_B^2
\end{align*}
So, a word of caution: in what follows, $\nabla$ denotes the connection (and not the gradient).
\subsubsection*{\bf \textit{The Ricci tensor}} Using the Levi-Civita connection that is computed in above, one calculates the Ricci tensor as follows. 
\begin{proposition}[e.g. \cite{FGKU}]\label{prop:dw-ric}
	The Ricci tensor of the doubly warped product $ B {_\alpha\times_\beta} F$ is given by
	\begin{align}
	\Ric(X+V,Y+W) &= {^B\Ric}(X,Y) +{^F\Ric}(V,W)  \notag \\ & -\left< X,Y  \right> \left(\Delta a + 2 \| \nabla a  \|^2 \right) - \left< V,W  \right> \left(\Delta b + 2 \| \nabla b  \|^2 \right) \notag \\ & - n_2 \left( {^B\nabla}^2b (X,Y)  + \nabla_Xb  \nabla_Yb \right) - n_1 \left({^F\nabla}^2 a (V,W) +  \nabla_Va \nabla_Wa \right) \notag \\ & + (n-2)\nabla_X b \nabla_Wa + (n-2)\nabla_Y b \nabla_Va. \notag
	\end{align}
	where $n_1 = \dim B$, $n_2 = \dim F$ and $n = n_1 + n_2$.
\end{proposition}
\subsection{ Lower bounds for the $\Ric$ tensor}
We start off by deriving pointwise lower Ricci curvature bounds for the doubly warped products under dynamic concavity conditions on the warping functions.
\subsubsection*{\bf \textit{Proof of Theorem~\ref{thm:dw-ric-1}}}
Tracing (\ref{eq:dw-ric-1}) over the frame $\{X_i\}^{n_1}_{i=1}$ (which is orthonormal w.r.t. $g_B$), one obtains
	\begin{align*}
		{^B\Delta}b + \alpha^{2}\Big( n_2^{-1}\left(n_1 -2\right) + 2 \Big) \| \nabla b \|^2 &\le \left(\left(n_2 -2\right) +  \Big( n_2^{-1}\left(n_1 -2\right) + 2 \Big)\right)\| \nabla b \|_B^2 \\ & = \left(n_2 -2\right)  \| \nabla b \|_B^2 + \alpha^{2}\Big( n_2^{-1}\left(n_1 -2\right) + 2 \Big) \| \nabla b \|^2 - \alpha^{-2}L_1,
	\end{align*}
	hence,
	\begin{align}\label{eq:dw-ric-3}
	 \Delta b + 2 \| \nabla b  \|^2 = \alpha^{-2}\left( {^B\Delta}b- \left(n_2 -2 \right) \| \nabla b \|_B^2 \right) \le - \alpha^{-2}L_1.
	\end{align}
	Similarly, tracing (\ref{eq:dw-ric-2}) over the frame $\{V_i\}^{n_2}_{i=1}$), we deduce
	\begin{align}\label{eq:dw-ric-4}
		\Delta a + 2\| \nabla a  \|^2 \le - \beta^{-2}L_2.
	\end{align}
	Therefore, using (\ref{eq:dw-ric-1}) and \eqref{eq:dw-ric-1} in combination with (\ref{eq:dw-ric-3}) and \eqref{eq:dw-ric-4}, we can write
	\begin{align*}
		\Ric(X+V,X+V) &\ge \left( n_1 -1 \right) \frac{\K_B}{\alpha^2} \| X \|^2  + \left( n_2 - 1 \right) \frac{\K_F}{\beta^2} \| V \|^2 \notag \\ & - \| X \|^2 \left(  \Delta a + 2\| \nabla a  \|^2 \right)  - \|V\|^2 \Big(  \Delta b + 2 \| \nabla b  \|^2 \Big) \notag \\ & - n_2 \left( \underset{B}{\nabla}^2 b (X,X)  +  \Big( \frac{n_1 -2}{n_2} + 2 \Big) \left(  \nabla_X b  \right)^2  \right) \notag \\ & - n_1 \left(\underset{F}{\nabla}^2 a (V,V)  +  \Big( \frac{n_2 -2}{n_1} + 2 \Big) \left(  \nabla_V a \right)^2  \right) \notag \\ &\ge \Big( \left( n_1 -1 \right) \frac{\K_B}{\alpha^2} + \frac{K_1}{\alpha^{2}} \Big)\|X\|^2 + \Big( \left( n_2 -1 \right) \frac{\K_F}{\beta^2} +  \frac{K_2}{\beta^{2}} \Big)\|V\|^2 \notag \\ &\ge \min\left\{ \frac{\left( n_1 -1 \right) \K_B + K_1 + L_1}{\alpha^2}, \frac{\left( n_2 -1 \right) \K_F + K_2 + L_2}{\beta^2} \right\} \| X + V \|^2 \\ &= \left( n_1 + n_2 -1  \right) \K \| X + V \|^2; 
	\end{align*}
	notice that in the first inequality,  we have also used the fact that by the Young's inequality, and for any integer $l$, vector fileds $X$ and $Y$ and a functions $f,g$, one has
	\begin{align*}
		l \nabla_X f \nabla_V g \ge - |l| \left(  \left| \nabla_X f \right|^2 + \left| \nabla_V g \right|^2 \right);
	\end{align*}
	the second inequality then follows from (\ref{eq:dw-ric-1}) and (\ref{eq:dw-ric-2}). 
\qed
\subsection{ $\N$-Bakry-\'Emery Ricci tensor for $\left( R_1 , R_2  \right)$-doubly warped products}\label{sec:NBE-tensor}
\begin{proposition}
	Let  $n = n_1 + n_2$ and $\M > n$. Then the $\M$-Bakry-\'Emery Ricci tensor of $B\; {_\alpha^{R_1}\times_\beta^{R_2}} \;F$ is
\begin{align}\label{eq:nice-1}
	\Ric_\upchi^{\mathcal{M}}(X + V) &= {^B\Ric}^\M_\Phi(X) + {^F\Ric}^\M_\Psi(V) \\ & -  \|X\|^2 \Big( \Delta a + \left(R_1 - n_1 + 2  \right)  \| \nabla a \|^2  \Big) \notag \\ & - \|V\|^2 \Big( \Delta b + \left(R_2 - n_2 + 2  \right)  \| \nabla b \|^2  \Big) \notag \\ & - R_1 {^F\nabla}^2_{V,V} a - R_2 {^B\nabla}^2_{X,X} b + \|X\|^2 \left< \nabla \Psi , \nabla a  \right> + \|V\|^2 \left< \nabla \Phi , \nabla b  \right> \notag \\ &+ \left( \M - n \right)^{-1}\mathcal{Q}_{\M}\left( x_1,x_2,x_3,x_4  \right),  \notag
\end{align}
and also in the special case $\N = \N_1 + \N_2$ with $\N_1 > n_1$ and $\N_2 > n_2$, given by
\begin{align*}
	\Ric_\upchi^{\N}(X + V) & = {^B\Ric}^{\N_1}_\Phi(X) + {^F\Ric}^{\N_2}_\Psi(V) \\ & -  \|X\|^2 \Big( \Delta a + \left(R_1 - n_1 + 2  \right)  \| \nabla a \|^2  \Big) \\ & - \|V\|^2 \Big( \Delta b + \left(R_2 - n_2 + 2  \right)  \| \nabla b \|^2  \Big) \\ & - R_1 {^F\nabla}^2_{V,V} a - R_2 {^B\nabla}^2_{X,X} b + \|X\|^2 \left< \nabla \Psi , \nabla a  \right> + \|V\|^2 \left< \nabla \Phi , \nabla b  \right> \\ & + \left( \N - n \right)^{-1}\mathcal{Q}_{\N_1,\N_2}\left( x_1,x_2,x_3,x_4  \right),
\end{align*}
where $(\M - n)^{-1} \mathcal{Q}_{\M}$ is the quadratic form corresponding to symmetric matrix
\begin{align}\label{eq:matrix-1}
	A_{\M} = \frac{1}{d}\begin{pmatrix} 
		-\frac{n_2}{a_1} & -1 & n_2 - c_1 & b_2 \\ 
		-1 &  	- \frac{n_1}{a_2} & b_1   &  n_1 - c_2 \\ 
		n_2 - c_1 &  b_1  &  -b_1^2-n_1d & (n-2)d - b_1b_2 + d(b_1+b_2)\\
		b_2  & n_1 - c_2  & (n-2)d - b_1b_2 + d(b_1+b_2)  &  -b_2^2-n_2d 
	\end{pmatrix},
\end{align}
in which
\begin{align*}
	a_i = \M - n_i  \quad b_i = R_i - n_i \quad c_i = \M - R_i \quad d = \M - n;
\end{align*}
and $(\N - n)^{-1} \mathcal{Q}_{\N_1,\N_2}$ is the quadratic form corresponding to symmetric matrix
\begin{align}\label{eq:matrix-2}
	A_{\N_1,\N_2} = \frac{1}{d}\begin{pmatrix} 
		\frac{a_2}{a_1} & -1 & n_2 - c_1 & b_2 \\ 
		-1 &  	 \frac{a_1}{a_2} & b_1   &  n_1 - c_2 \\ 
		n_2 - c_1 &  b_1  &  -b_1^2-n_1d & (n-2)d - b_1b_2 + d(b_1+b_2)\\
		b_2  & n_1 - c_2  & (n-2)d - b_1b_2 + d(b_1+b_2)  &  -b_2^2-n_2d 
	\end{pmatrix},
\end{align}
in which
\begin{align*}
	a_i = \N_i - n_i  \quad b_i = R_i - n_i \quad c_i = \N_i - R_i \quad d = \ - n.
\end{align*}
\end{proposition}
\begin{proof}
	Recall $\chi = \left(n_1 - R_1 \right)a + \left( n_2 - R_2 \right)b + \Phi + \Psi$, $a := \ln \alpha$ and $b := \ln \beta$. By the $\nabla$-calculus form \S~\ref{sec:calc-stuff}, we get
\begin{align*}
	&\nabla^2 \upchi \left( X+V , X+V  \right) \\&= \nabla^2 \left(  \left( n_2 - R_2 \right)b + \Phi \right) \left( X+V , X+V  \right) + \nabla^2 \left( \left(n_1 - R_1 \right)a + \Psi \right) \left( X+V , X+V  \right) \\ &=  {^B\nabla}^2 \Phi (X,X) + {^F\nabla}^2 \Psi (V,V) \\ &+ \left( n_2 - R_2 \right) {^B\nabla}^2 b (X,X) + \left( n_1 - R_1 \right) {^F\nabla}^2 a (V,V) \\ &+ \|X\|^2\left( \left( n_1 - R_1 \right)\|\nabla a\|^2 + \left< \nabla \Psi, \nabla a \right> \right) + \|V\|^2 \left( \left( n_2 - R_2 \right)\|\nabla b\|^2 + \left< \nabla \Phi, \nabla b \right> \right) \\ & - 2\nabla_V a \nabla_X\Phi  - 2\left( n_2 - R_2 \right) \nabla_V a \nabla_Xb \\ &  - 2\nabla_X b \nabla_V\Psi  - 2\left( n_1 - R_1 \right) \nabla_X b \nabla_Va \\ &= {^B\nabla}^2 \Phi(X,X) - 2 \nabla_V a \nabla_X \Phi  + \|V\|^2 \left< \nabla \Phi , \nabla b  \right> \notag\\ &+ {^F\nabla}^2 \Psi(V,V) - 2 \nabla_X b \nabla_V \Psi + \|X\|^2 \left< \nabla \Psi , \nabla a  \right> \notag \\ & - \left( R_1 - n_1   \right)\left( {^F\nabla}^2 a(V,V)- 2 \nabla_X b \nabla_V a + \|X\|^2 \| \nabla a  \|^2  \right) \notag \\ & - \left( R_2 - n_2 \right)\left( {^B\nabla}^2 b(X,X) - 2 \nabla_X b \nabla_V a + \|V\|^2 \| \nabla b  \|^2   \right), \notag
\end{align*}	
also
\begin{align*}
	\nabla \upchi \otimes \nabla \upchi (X+V,X+V) &= \Big( \left(n_1 - R_1 \right)\nabla_Va + \left( n_2 - R_2 \right)\nabla_Xb + \nabla_X\Phi + \nabla_V\Psi \Big)^2 	\\&= \left( \nabla_X \Phi  \right)^2 + \left( \nabla_V \Psi  \right)^2 + \left( R_1 - n_1   \right)^2 \left( \nabla_V a  \right)^2 + \left( R_2 - n_2  \right)^2 \left( \nabla_X b  \right)^2 \notag \\ & +  2\nabla_X\Phi \nabla_V \Psi -2\left(  R_1 - n_1 \right)\nabla_X\Phi\nabla_Va - 2\left(  R_2 - n_2 \right)\nabla_X\Phi\nabla_Xb \notag \\ & - 2\left(  R_1 - n_1 \right)\nabla_V\Psi\nabla_Va - 2\left(  R_2 - n_2 \right)\nabla_V\Psi\nabla_Xb \notag \\ & + 2\left(  R_1 - n_1 \right)\left(  R_2 - n_2 \right) \nabla_Va \nabla_Xb \notag.
\end{align*}
\par Therefore, by using Proposition~\ref{prop:dw-ric} and the definition of $\M$-Bakry-\'Emery Ricci tensor, for any $\M > n$, we get
\begin{align*}
		&\Ric_\upchi^{\mathcal{M}}(X + V) = \\ & \Ric(X + V)  + \nabla^2 \upchi\left(X + V , X + V\right) - \left(   \mathcal{M} - n\right)^{-1}\nabla_{_{X + V}}\upchi \otimes \nabla_{_{X + V}} \upchi \\ &= \blue{{^B\Ric}(X)} + \red{{^F\Ric}(V)} \notag \\ & - \|X\|^2 \left(  \Delta a  +  2\|\nabla a\|^2 \right)  - \|V\|^2 \left(  \Delta b +  2\|\nabla b\|^2 \right) \notag \\ & -n_2 {^B\nabla}^2 b(X,X) -n_2  \left( \nabla_X b \right)^2 - n_1 {^F\nabla}^2 a(V,V) - n_1 \left( \nabla_V a \right)^2
	+ 2(n-2)\nabla_Xb\nabla_Va \notag \\ & \blue{ +{^B\nabla}^2 \Phi(X,X)} - 2 \nabla_V a \nabla_X \Phi  + \|V\|^2 \left< \nabla \Phi , \nabla b  \right> \notag\\ &\red{+ {^F\nabla}^2 \Psi(V,V)} - 2 \nabla_X b \nabla_V \Psi + \|X\|^2 \left< \nabla \Psi , \nabla a  \right> \notag \\ & - \left( R_1 - n_1   \right)\left( {^F\nabla}^2 a (V,V) - 2 \nabla_X b \nabla_V a + \|X\|^2 \| \nabla a  \|^2  \right) \notag \\ & - \left( R_2 - n_2 \right)\left( {^B\nabla}^2 b(X,X) - 2 \nabla_X b \nabla_V a + \|V\|^2 \| \nabla b  \|^2   \right) \notag \\ & - \frac{1}{\M-n}  \left( \nabla_X \Phi  \right)^2 - \frac{1}{\M-n}  \left( \nabla_V \Psi  \right)^2 + \frac{1}{\M-n_1}  \left( \nabla_X \Phi  \right)^2 + \frac{1}{\M-n_2}  \left( \nabla_V \Psi  \right)^2     \\& \blue{- \frac{1}{\M-n_1}  \left( \nabla_X \Phi  \right)^2} \red{- \frac{1}{\M-n_2}  \left( \nabla_V \Psi  \right)^2}      \\&
	  - \frac{\left( R_1 - n_1   \right)^2}{\M-n}  \left( \nabla_V a  \right)^2 -  \frac{\left( R_2 - n_2  \right)^2}{\M-n}\left( \nabla_X b  \right)^2 \notag \\ & - \frac{2}{\M-n}\nabla_X\Phi\nabla_V\Psi + \frac{2\left( R_1 - n_1  \right)}{\M-n}\nabla_X\Phi \nabla_Va + \frac{2\left( R_2 - n_2  \right)}{\M-n}\nabla_V\Psi \nabla_Xb \notag \\ & + \frac{2\left( R_1 - n_1  \right)}{\M-n} \nabla_V\Psi \nabla_Va + \frac{2\left( R_2 - n_2  \right)}{\M-n} \nabla_X\Phi \nabla_Xb - \frac{2\left( R_1 - n_1 \right)\left( R_2 - n_2 \right)}{\M-n}\nabla_Va\nabla_Xb \notag.
\end{align*}
Therefore, setting $x_1 := \nabla_X \Phi$, $x_2 := \nabla_V \Psi$, $x_3 := \nabla_V a$ and $x_4 := \nabla_X b$, we can simplify as follows
\begin{align*}
	\Ric_\upchi^{\mathcal{M}}(X + V) &= {^B\Ric}^\M_\Phi(X) + {^F\Ric}^\M_\Psi(V) \\ & -  \|X\|^2 \Big( \Delta a + \left(R_1 - n_1 + 2  \right)  \| \nabla a \|^2  \Big) \\ & - \|V\|^2 \Big( \Delta b + \left(R_2 - n_2 + 2  \right)  \| \nabla b \|^2  \Big) \\ & - R_1 {^F\nabla}^2_{V,V} a - R_2 {^B\nabla}^2_{X,X} b + \|X\|^2 \left< \nabla \Psi , \nabla a  \right> + \|V\|^2 \left< \nabla \Phi , \nabla b  \right> \\ &   - \frac{1}{\M-n}  \left( x_1 \right)^2 - \frac{1}{\M-n}  \left( x_2 \right)^2 + \frac{1}{\M-n_1}  \left( x_1  \right)^2 + \frac{1}{\M-n_2}  \left( x_2  \right)^2    \\&
	 -n_2  \left( x_4 \right)^2 - n_1 \left( x_3 \right)^2
	+ 2(n-2)  x_3x_4   - 2  x_1x_3 - 2 x_2x_4 \notag \\ & + 2\left( R_1 - n_1   \right)\left( x_4 x_3  \right)  + 2\left( R_2 - n_2 \right)\left( x_3x_4   \right) \notag \\ &  - \frac{\left( R_1 - n_1   \right)^2}{\M-n}  \left( x_3  \right)^2 -  \frac{\left( R_2 - n_2  \right)^2}{\M-n}\left( x_4 \right)^2 \notag \\ & - \frac{2}{\M-n} x_1 x_2 + \frac{2\left( R_1 - n_1  \right)}{\M-n} x_1 x_3+ \frac{2\left( R_2 - n_2  \right)}{\M-n} x_2 x_4 \notag \\ & + \frac{2\left( R_1 - n_1  \right)}{\M-n} x_2 x_3 + \frac{2\left( R_2 - n_2  \right)}{\M-n} x_1 x_4 - \frac{2\left( R_1 - n_1 \right)\left( R_2 - n_2 \right)}{\M-n} x_3 x_4 \notag  \\ &= {^B\Ric}^\M_\Phi(X) + {^F\Ric}^\M_\Psi(V) \\ & -  \|X\|^2 \Big( \Delta a + \left(R_1 - n_1 + 2  \right)  \| \nabla a \|^2  \Big) \\ & - \|V\|^2 \Big( \Delta b + \left(R_2 - n_2 + 2  \right)  \| \nabla b \|^2  \Big) \\ & - R_1 {^F\nabla}^2_{V,V} a - R_2 {^B\nabla}^2_{X,X} b + \|X\|^2 \left< \nabla \Psi , \nabla a  \right> + \|V\|^2 \left< \nabla \Phi , \nabla b  \right> \\ &+ \left( \M - n \right)^{-1}\mathcal{Q}_{\M}\left( x_1,x_2,x_3,x_4  \right),
\end{align*}
	where
\begin{align*}
	&\mathcal{Q}_{\M}\left( x_1,x_2,x_3,x_4  \right) = \\& \left( -1 + \frac{\M - n}{\M - n_1}   \right)x_1^2 +  \left( -1 + \frac{\M - n}{\M - n_2}   \right) x_2^2\\&- \Big( n_1(\M-n) + \left( R_1 - n_1 \right)^2   \Big) x_3^2  - \Big( n_2(\M-n) +  \left(R_2 - n_2 \right)^2   \Big) x_4^2  \notag \\ & + 2 \left(n_2 + R_1 - \M  \right)x_1x_3 + 2 \left(n_1 + R_2 - \M  \right)x_2x_4 \notag \\ & + 2 \left( R_1 - n_1  \right)x_2x_3 + 2 \left( R_2 - n_2  \right)x_1x_4 -2 x_1x_2 \notag \\ & + 2 \Big( (\M - n)(n-2) - \left(R_1 - n_1  \right) \left(R_2 - n_2  \right) + (\M - n) \left( R_1 - n_1 \right) + (\M - n) \left( R_2 - n_2 \right) \Big) x_3x_4. \notag 
\end{align*}
So, $(\M - n)^{-1} \mathcal{Q}_{\M}\left( \mathbf{x} \right) = \mathbf{x} A_{\M} \mathbf{x}^T$ for the symmetric matrix $A_{\M}$ from \eqref{eq:matrix-1}. 
\par For the special case $\N = \N_1 + \N_2$, using the identities
\begin{align*}
		{^B\Ric}^{\N_1 + \N_2}_\Phi(X) &= {^B\Ric}^{\N_1}_\Phi(X) + \left( \frac{1}{\N_1 - n_1} -\frac{1}{\N_1 + \N_2 - n_1}   \right) \left(x_1\right)^2 \\ &= {^B\Ric}^{\N_1}_\Phi(X) + \frac{\N_2}{\left( \N_1 - n_1\right)\left( \N_1 + \N_2 - n_1\right)}x_1^2,
\end{align*}
and 
\begin{align*}
	{^F\Ric}^{\N_1 + \N_2}_\Psi(V) &= {^F\Ric}^{\N_2}_\Psi(V) + \left( \frac{1}{\N_2 - n_2} -\frac{1}{\N_1 + \N_2 - n_2}   \right) \left(x_2\right)^2 \\ &= {^F\Ric}^{\N_2}_\Psi(V) + \frac{\N_1}{\left( \N_2 - n_2\right)\left( \N_1 + \N_2 - n_2\right)}x_2^2,
\end{align*}
we deduce,
\begin{align*}
	\Ric_\upchi^{\N_1 + \N_2}(X + V) &= {^B\Ric}^{\N_1}_\Phi(X) + {^F\Ric}^{\N_2}_\Psi(V) \\ & -  \|X\|^2 \Big( \Delta a + \left(R_1 - n_1 + 2  \right)  \| \nabla a \|^2  \Big) \\ & - \|V\|^2 \Big( \Delta b + \left(R_2 - n_2 + 2  \right)  \| \nabla b \|^2  \Big) \\ & - R_1 {^F\nabla}^2_{V,V} a - R_2 {^B\nabla}^2_{X,X} b + \|X\|^2 \left< \nabla \Psi , \nabla a  \right> + \|V\|^2 \left< \nabla \Phi , \nabla b  \right> \\ &+ \left( \N_1 + \N_2 - n \right)^{-1}\mathcal{Q}_{\N_1 + \N_2}\left( x_1,x_2,x_3,x_4  \right) \\ & + \frac{\N_2}{\left( \N_1 - n_1\right)\left( \N_1 + \N_2 - n_1\right)}x_1^2  + \frac{\N_1}{\left( \N_2 - n_2\right)\left( \N_1 + \N_2 - n_2\right)}x_2^2 \\& = {^B\Ric}^{\N_1}_\Phi(X) + {^F\Ric}^{\N_2}_\Psi(V) \\ & -  \|X\|^2 \Big( \Delta a + \left(R_1 - n_1 + 2  \right)  \| \nabla a \|^2  \Big) \\ & - \|V\|^2 \Big( \Delta b + \left(R_2 - n_2 + 2  \right)  \| \nabla b \|^2  \Big) \\ & - R_1 {^F\nabla}^2_{V,V} a - R_2 {^B\nabla}^2_{X,X} b + \|X\|^2 \left< \nabla \Psi , \nabla a  \right> + \|V\|^2 \left< \nabla \Phi , \nabla b  \right> \\ & + \left( \N_1 + \N_2 - n \right)^{-1}\mathcal{Q}_{\N_1,\N_2}\left( x_1,x_2,x_3,x_4  \right),
\end{align*}
where
\begin{align*}
	& \mathcal{Q}_{\N_1,\N_2}\left( x_1,x_2,x_3,x_4  \right) = \\ &\frac{\N_2 - n_2}{\N_1 -n_1 }x_1^2 + \frac{\N_1 - n_1}{\N_2 - n_2 }x_2^2  -2 x_1x_2 \notag \\ & - \Big( n_1(\N-n) + \left( R_1 - n_1 \right)^2   \Big) x_3^2  - \Big( n_2(\N-n) +  \left(R_2 - n_2 \right)^2   \Big) x_4^2  \notag \\ & - 2 \left(n_2 + R_1 - \N   \right)x_1x_3 - 2 \left(n_1 + R_2 - \N   \right)x_2x_4 \notag \\ & + 2 \left( R_1 - n_1  \right)x_2x_3 + 2 \left( R_2 - n_2  \right)x_1x_4 \notag \\ & + 2 \Big( (\N-n)(n-2) - \left(R_1 - n_1  \right) \left(R_2 - n_2  \right) + (\M - n) \left( R_1 - n_1 \right) + (\M - n) \left( R_2 - n_2 \right)  \Big) x_3x_4 \notag,
\end{align*}
which means quadratic form $(\N - n)^{-1} \mathcal{Q}_{\N_1,\N_2}$ is  given by the symmetric matrix \eqref{eq:matrix-2}. 
\end{proof}
\subsection{Bakry-\'Emery curvature-dimension bounds for $\left( R_1 , R_2  \right)$-doubly warped products}\label{sec:NBE-bounds}
\subsubsection*{\bf \textit{Proof of Theorem~\ref{thm:main-1-1}}}
Let $\lambda^{\N}$ denote the smallest eigenvalues of $A_{\N}$. Then, from \eqref{eq:nice-1}, we deduce
\begin{align}\label{eq:wdw-ric}
	\Ric_\chi^{\N}(X + V) \ge & \alpha^{-2}\left( n_1 - 1  \right)\K_B \|  X \|^2 + \beta^{-2}\left( n_2 - 1  \right)\K_F \| V \|^2 \notag \\ &-  \|X\|^2 \left( \Delta a + \left(R_1 - n_1 + 2  \right)  \| \nabla a \|^2 \right) - \|V\|^2 \left( \Delta b + \left(R_2 - n_2 + 2  \right)  \| \nabla b \|^2 \right) \notag \\&   - \|X\|^2 \left(  \nicefrac{1}{2}\|\nabla \Psi\|^2 + \nicefrac{1}{2} \|\nabla a\|^2  \right) - \|V\|^2 \left(  \nicefrac{1}{2}\|\nabla \Phi\|^2 + \nicefrac{1}{2} \|\nabla b\|^2  \right) \notag
	\\& - R_2 \Big( {^B\nabla}^2_{X,X} b - \frac{\lambda^{\N}}{R_2} \left( \nabla_X \Phi \right)^2  - \frac{\lambda^{\N}}{R_2} \left( \nabla_X b \right)^2 \Big)  \\ & - R_1 \Big( {^F\nabla}^2_{V,V} a - \frac{\lambda^{\N}}{R_1} \left( \nabla_V \Psi \right)^2  - \frac{\lambda^{\N}}{R_1} \left( \nabla_V a \right)^2 \Big).\notag
\end{align}
\par Now, tracing the concavity condition \eqref{eq:main-con-rel-1-1} over the frame $X_i$ (orthornormal w.r.t. $g_B$), we get
\begin{align*}
	&R_2 {^B\Delta} b  - \alpha^{2}\lambda^{\N}\| \nabla \Phi \|^2  - \alpha^{2}\lambda^{\N} \| \nabla b \|^2 \\ &\le - \left( \lambda^{\N} + \nicefrac{1}{2}R_2 \right)\| \nabla \Phi \|_B^2 - \Big(\lambda^{\N} + R_2\left(R_2 -2n_2 + 2   \right) + \nicefrac{1}{2}R_2 \Big)\|\nabla b\|_B^2  -R_2L_1,
\end{align*}
which upon cancellation of the like terms, dividing by $R_2$ and rearranging terms, becomes
\begin{align*}
	{^B\Delta} b + \nicefrac{1}{2}\|\nabla \Phi\|_B^2 + \nicefrac{1}{2} \|\nabla b\|_B^2 \le - \left(R_2 -2n_2 + 2   \right)\|\nabla b\|_B^2 -L_1.
\end{align*}	
\par Consequently, one has
\begin{align}\label{eq:main-bound-1}
	&\Delta b + \left(R_2 - n_2 + 2  \right)  \| \nabla b \|^2 + \nicefrac{1}{2}\|\nabla \Phi\|^2 + \nicefrac{1}{2} \|\nabla b\|^2\\  &= \frac{1}{\alpha^2} {^B\Delta} b  - n_2\|\nabla b\|^2 + \left(R_2 -n_2 + 2   \right) \|\nabla b\|^2  + \nicefrac{1}{2}\|\nabla \Phi\|^2 + \nicefrac{1}{2} \|\nabla b\|^2 \notag \\ &=  \frac{1}{\alpha^2} {^B\Delta} b  + \left(R_2 - 2n_2 + 2   \right) \|\nabla b\|^2 + \nicefrac{1}{2}\|\nabla \Phi\|^2 + \nicefrac{1}{2} \|\nabla b\|^2 \notag \\ &= \frac{1}{\alpha^2} \Big(   {^B\Delta} b  + \left(R_2 -2n_2 + 2   \right)\|\nabla b\|_B^2 + \nicefrac{1}{2}\|\nabla \Phi\|_B^2 + \nicefrac{1}{2} \|\nabla b\|_B^2 \Big) \le -\alpha^{-2}L_1. \notag
\end{align}
In a similar fashion, by tracing \eqref{eq:main-con-rel-1-2} over the orthonormal frame $\left\{ V_i \right\}$, we deduce
\begin{align*}
	\Delta a + \left(R_1 - n_1 + 2  \right)  \| \nabla a \|^2 + \nicefrac{1}{2}\|\nabla \Psi\|^2 + \nicefrac{1}{2} \|\nabla a\|^2 \le -\beta^{-2}L_2.
\end{align*}
 \par By how we have picked the constant $\lambda^\N$, we have
\begin{align}\label{eq:main-bound-2}
	&- R_1 {^F\nabla}^2_{V,V} a - R_2 {^B\nabla}^2_{X,X} b + \left( \M - n \right)^{-1}\mathcal{Q}_{\M}\left( x_1,x_2,x_3,x_4  \right) \notag \\ &\ge  - R_1 {^F\nabla}^2_{V,V} a - R_2 {^B\nabla}^2_{X,X} b  + \lambda^{\N} \left(  x_1^2 + x_2^2 + x_3^2 + x_4^2 \right) \notag \\ &= - R_2 \Big( {^B\nabla}^2_{X,X} b - \frac{\lambda^{\N}}{R_2} \left( \nabla_X \Phi \right)^2  - \frac{\lambda^{\N}}{R_2} \left( \nabla_X b \right)^2 \Big)  \\ & - R_1 \Big( {^F\nabla}^2_{V,V} a - \frac{\lambda^{\N}}{R_1} \left( \nabla_V \Psi \right)^2  - \frac{\lambda^{\N}}{R_1} \left( \nabla_V a \right)^2 \Big). \notag
\end{align}
\par Finally, using \eqref{eq:wdw-ric} in combination with \eqref{eq:main-bound-1} and \eqref{eq:main-bound-2}, we get
\begin{align*}
	\Ric_\chi^{\N}(X + V) &\ge  \alpha^{-2}\left( n_1 - 1  \right)\K^\N_B \|  X \|^2 + \beta^{-2}\left( n_2 - 1  \right)\K^\N_F \| V \|^2\\ & +\alpha^{-2}K_1\|X\|^2 + \beta^{-2}K_2\|V\|^2  + \alpha^{-2}L_1\|X\|^2 +  \beta^{-2}L_2\|V\|^2\\ &\ge \min\left\{ \frac{\left( n_1 - 1  \right)\K^{\N}_B + K_1 + L_1 }{\alpha^2}, \frac{ \left( n_2 - 1  \right)\K^{\N}_F + K_2 + L_2}{\beta^2}\right\} \|X + V\|,
\end{align*}
which is the desired conclusion.  \qed
\subsubsection*{\bf \textit{Proof of Theorem~\ref{thm:main-1-2}}}
The proof is almost verbatim as in the proof of Theorem~\ref{thm:main-1-1}; the only difference being this time we take $\lambda^{\N_1,\N_2}$ to be the lowest eigenvalue of the quadratic form $A_{\N_1,\N_2}$. 
\qed
\subsubsection*{\bf \textit{Bounds on the (best) constants $\lambda^\M$ and $\lambda^{\N_1, \N_2}$}}
\begin{proposition}
	\begin{align*}
\left| \lambda^{\M} \right| < 7n + 4R_1 + 4R_2 + \frac{2 + 3R_1 + 3R_2 + 6R_1^2 + 6R_2^2 + 3\M + 3n + 8n^2}{\left( \M - n \right)} + \frac{2n}{\left( \M - n \right)^2}.
	\end{align*}	
\end{proposition}
\begin{proof}
	The \emph{operator norm - trace norm} inequality implies
	\begin{align*}
		&|\lambda^\M|^2 \le \left( \sum_{i,j} {A_{\M}}^2_{ij}  \right)^2\\&= \frac{2}{d^2}\Big(  1 + b_1^2 + b_2^2 + \frac{n_2^2}{a_1^2} + \frac{n_1^2}{a_2^2} +  \left( n_2 - c_1  \right)^2 +  \left( n_1 - c_2  \right)^2  \Big. \\ &\Big.  \phantom{sajjadsajjad} +\left( b_1^2 + n_1d  \right)^2 + \left( b_2^2 + n_2d  \right)^2 + \left( (n-2)d - b_1b_2 + d(b_1+b_2) \right)^2\Big),
	\end{align*}
where
\begin{align*}
	a_i = \M - n_i  \quad b_i = R_i - n_i \quad c_i = \M - R_i \quad d = \M - n.
\end{align*}
\par Applying the Young's inequality a few times, we get the following (very crude) estimate
	\begin{align*}
	&|\lambda^\M|^2 < \\& \frac{2}{\left( \M - n \right)^2} \Big( 1 + 2R_1^2 + 2R_2^2  + 4n^2 + \frac{n^2}{\left( \M - n \right)^2} + 2\M^2 + 2R_1^2 + 2\M^2 + 2R_2^2 \Big. \\ &\Big. \phantom{sajjadsajjad} + 8R_1^4 + 8n^4 + 2n^2\left( \M - n\right)^2 +  8R_2^4 + 8n^4 + 2n^2\left( \M - n\right)^2  \Big. \\ & \Big. \phantom{sajjadsajjad} + 2n^2\left( \M - n  \right)^2 + 8R_1^4 + 8R_2^4 + 10n^4 + 2\left( 4R_1^2 + 4R_2^2 + 8n^2 \right)  \left( \M -n  \right)^2 \Big) \\ & = 44n^2 + 16R_1^2 + 16R_2^2 + \frac{2 + 8R_1^2 + 8R_2^2 + 32R_1^4 + 32R_2^4 + 8\M^2 + 8n^2 + 52n^4}{\left( \M - n \right)^2} + \frac{2n^2}{\left( \M - n \right)^4}.
\end{align*}
\par Therefore,
\begin{align*}
	|\lambda^\M| < 7n + 4R_1 + 4R_2 + \frac{2 + 3R_1 + 3R_2 + 6R_1^2 + 6R_2^2 + 3\M + 3n + 8n^2}{\left( \M - n \right)} + \frac{2n}{\left( \M - n \right)^2}.
\end{align*}
\end{proof}
\begin{proposition}
	For $\N = \N_1 + \N_2$, we have
	\begin{align*}
	\lambda^{\N} \le \lambda^{\N_1, \N_2} &\le \lambda^{\N} + \left(  \N - n\right)^{-1} \left( \N_2\left( \N_1 - n_1\right)^{-1} \vee \N_1\left( \N_2 - n_2\right)^{-1} \right).
	\end{align*}
	In particular,
	\begin{align*}
	\lambda^{2\N} \le \lambda^{\N, \N} < \lambda^{\N,\N} + \frac{\N}{(\N - n)^2}.
	\end{align*}
\end{proposition}
\begin{proof}
	\begin{align*}
	A_{\N_1,\N_2} - A_{\N} = \begin{pmatrix}
	 \frac{\N_2}{\left( \N_1 - n_1\right)\left( \N - n_1\right)} & 0 \\
	0 & \frac{\N_1}{\left( \N_2 - n_2\right)\left( \N- n_2\right)}
	\end{pmatrix} \oplus 0_{2\times 2}
	\end{align*}
	which is a \emph{nonnegative definite diagonal} matrix. Therefore, by elementary matrix theory, we have
\begin{align*}
	0 \le \lambda^{\N_1, \N_2} - \lambda^{\N}  &\le \max\left\{ \frac{\N_2}{\left( \N_1 - n_1\right)\left( \N - n_1\right)},  \frac{\N_1}{\left( \N_2 - n_2\right)\left( \N- n_2\right)} \right\} \\ &< \max\left\{ \frac{\N_2}{\left( \N_1 - n_1\right)\left( \N - n\right)},  \frac{\N_1}{\left( \N_2 - n_2\right)\left( \N- n\right)} \right\} \\ &= \left(  \N - n\right)^{-1} \left( \N_2\left( \N_1 - n_1\right)^{-1} \vee \N_1\left( \N_2 - n_2\right)^{-1} \right)
\end{align*}
\end{proof}
\section{Proof of main theorems in the discrete setting}
Recall the Definition~\ref{def:graph-dw} of doubly twisted product of weighted graphs. 
\subsection{ Structural lower bound on the curvature functions}
First we compute the constituent parts of the \textit{Ricci form}, $\Gamma_2$. Recall
\begin{align*}
\Gamma_2 \left(  u \right) := \nicefrac{1}{2} \big( \Delta \Gamma(u) - 2 \Gamma \left(  u , \Delta u \right) \big),
\end{align*}
where $\Delta$ and $\Gamma$ are as in (\ref{eq:Laplacian-g}) and (\ref{eq:Gamma1-g}). So,
\begin{align*}
\Delta u (x) := m_x^{-1} \sum_{y \sim x} \big( u(y) - u(x)  \big) \omega_{xy}, 
\end{align*}
and by now standard calculations, 
\begin{align}
\Gamma \left(u , v \right) (x) &:=  \nicefrac{1}{2} \left( \Delta (uv) - v\Delta u - u\Delta v \right)(x) \notag \\ &=  \left(2m_x\right)^{-1} \sum_{y \sim x} \big( u(y) - v(x)  \big)\big( u(y) - v(x)  \big) \omega_{xy}. \notag
\end{align}
Therefore,
\begin{align}
& 2\Gamma \left(  u , \Delta u \right)(x) \notag \\ &= m_x^{-1} \sum_{y \sim x} \big( u(y) - u(x)  \big)\big( \Delta u(y) - \Delta u(x)  \big) \omega_{xy}  \notag \\ &= m_x^{-1} \sum_{y \sim x} \big( u(y) - u(x)  \big)\left( m_y^{-1} \sum_{z \sim y} \big( u(z) - u(y)  \big) \omega_{yz} -  m_x^{-1} \sum_{w \sim x} \big( u(w) - u(x)  \big)\omega_{xw} \right) \omega_{xy},  \notag
\end{align}
and
\begin{align*}
	\begin{split}
&\Delta \Gamma \left( u  \right)(x)\\ &= \left( 2m_x  \right)^{-1} \sum_{ y \sim x} \left( m_y^{-1} \sum_{z \sim y} \big(  u(z) - u(y)   \big)^2\omega_{yz} -  m_x^{-1} \sum_{w \sim x} \big(  u(w) - u(x)   \big)^2\omega_{xy}\right) \omega_{xy} \notag \\ &=  \left(2m_x\right)^{-1} \left( \sum_{y \sim x} m_y^{-1} \omega_{xy}\sum_{z \sim y} \big( u(y) - u(z)    \big)^2 \omega_{yz}   - m_x^{-1} D_x \sum_{y \sim x} \big( u(x) - u(y)  \big)^2  \omega_{xy}    \right) \notag \\ &= \left( 2m_x  \right)^{-1} \left(  \sum_{y \sim x} m_y^{-1}\omega_{xy}\sum_{y \sim z}  \big( u(y) - u(z)    \big)^2 \omega_{yz}  \right. \\ & \left. \hspace{3cm} - m_x^{-1} D_x \sum_{y \sim x} m_y^{-1}D_y^{-1}\omega_{xy} \sum_{z\sim y}\big( u(x) - u(y)  \big)^2  \omega_{yz}     \right) \notag \\ &=  \left( 2m_x  \right)^{-1} \sum_{y \sim x} m_y^{-1}D_y^{-1}\omega_{xy} \sum_{z \sim y} \left( D_y\big(  u(y) - u(z)  \big)^2 - m_x^{-1}D_x\big( u(x) - u(y) \big)^2  \right) \omega_{yz} \notag \\ &=  \left( 2m_x  \right)^{-1} \sum_{y \sim x} m_y^{-1}D_y^{-1}\omega_{xy} \sum_{z \sim y} \left( m_x^{-\nicefrac{1}{2}}D_x^{\nicefrac{1}{2}}u(x)  -\left(D_y^{\nicefrac{1}{2}} + m_x^{-\nicefrac{1}{2}}D_x^{\nicefrac{1}{2}} \right)u(y) +  D_y^{\nicefrac{1}{2}}  u(z)  \right)^2 \omega_{yz} \notag \\ & - m_x^{-1} \sum_{y \sim x} m_y^{-1}D_y^{-1}\omega_{xy} \sum_{z \sim y} \left( m_x^{-\nicefrac{1}{2}}D_x^{\nicefrac{1}{2}}u(x)  -\left(D_y^{\nicefrac{1}{2}} + m_x^{-\nicefrac{1}{2}}D_x^{\nicefrac{1}{2}} \right)u(y) +  D_y^{\nicefrac{1}{2}}  u(z)  \right)  \notag \\ &\hspace{8.8cm} \boldsymbol{\cdot} \left(  m_x^{-\nicefrac{1}{2}}D_x^{\nicefrac{1}{2}}u(x) - m_x^{-\nicefrac{1}{2}}D_x^{\nicefrac{1}{2}}u(y)  \right)\omega_{yz} . \notag 
\end{split}
\end{align*}
Thus,
\begin{align}\label{eq:strbd-Gam2}
& \Gamma_2 \left(  u \right)(x) \notag \\ &= \left(4m_x\right)^{-1} \left( \sum_{y \sim x} m_y^{-1}D_y^{-1}\omega_{xy} \sum_{z \sim y} \left( m_x^{-\nicefrac{1}{2}}D_x^{\nicefrac{1}{2}}u(x)  -\left(D_y^{\nicefrac{1}{2}} + m_x^{-\nicefrac{1}{2}}D_x^{\nicefrac{1}{2}} \right)u(y) +  D_y^{\nicefrac{1}{2}}  u(z)  \right)^2 \omega_{yz}  \right) \notag  \\ & -  \left(2m_x\right)^{-1} \left( \sum_{y \sim x} m_y^{-1}\omega_{xy} m_x^{-\nicefrac{1}{2}}D_x^{\nicefrac{1}{2}} \big(  u(y) - u(x)   \big)^2 \right)   \\ & +  \left(2m_x\right)^{-1} \left( \sum_{y \sim x} m_y^{-1}\omega_{xy}   D_y^{-\nicefrac{1}{2}} \sum_{z \sim y} \big(  u(y) - u(x)   \big) \big(  u(z) - u(y)   \big) \omega_{yz} \right) \notag \\ & - \left(2m_x\right)^{-1} \left( \sum_{y \sim x} m_y^{-1}\omega_{xy} \sum_{z \sim y} \big(  u(y) - u(x)   \big) \big(  u(z) - u(y)   \big) \omega_{yz} \right) \notag  \\ &+ \nicefrac{1}{2}\left(m_x^{-1}  \sum_{y \sim x} \big(  u(y) - u(x)   \big)  \omega_{xy} \right)^2. \notag 
\end{align}
\subsubsection*{\bf \textit{Proof of Theorem~\ref{thm:main-4}}}
Setting $X = u(y) - u(x)$, $Y = u(z) - u(y)$, $a = m_x^{-\nicefrac{1}{2}}D_x^{\nicefrac{1}{2}}$ and $b = D_y^{\nicefrac{1}{2}} $, $\Gamma_2 \left(  u \right)(x)$ takes the form 
\begin{align*}
&\Gamma_2 \left(  u \right)(x) \\ &= \left(2m_x\right)^{-1} \sum_{y \sim x} m_y^{-1}\omega_{xy} \sum_{z \sim y} \left(\nicefrac{1}{2}\, b^{-2} \left( -aX + bY    \right)^2 + b^{-1}XY - XY - aX^2 \right)\omega_{yz} + \nicefrac{1}{2}\, \big( \Delta u (x)   \big)^2. 
\end{align*}
Applying the identity/inequality,
\begin{align*}
&\nicefrac{1}{2}\, b^{-2}\left( -aX + bY    \right)^2 + b^{-1}XY - XY \\&= \Big(  (2b)^{-1}\left(-2a + b^2 - b\right)X + Y     \Big)^2 +  \Big( a^2b^{-2} -\left( ab^{-1} - \nicefrac{1}{2}\, b + \nicefrac{1}{2}  \right)^2      \Big) X^2 \notag \\ &\ge  \Big( a^2b^{-2} -\left( ab^{-1} - \nicefrac{1}{2}\, b + \nicefrac{1}{2}  \right)^2      \Big) X^2, \notag
\end{align*}
with $X = u(y) - u(x)$, $Y = u(z) - u(y)$, $a = m_x^{-\nicefrac{1}{2}}D_x^{\nicefrac{1}{2}}$ and $b = D_y^{\nicefrac{1}{2}} $ yields
\begin{align*}
& \Gamma_2 \left(  u \right)(x) \notag \\  &\ge \left(2m_x\right)^{-1} \sum_{y \sim x} m_y^{-1}\omega_{xy} \left( m_xD_xD_y^{-1}  - \left( m_x^{-\nicefrac{1}{2}}D_x^{\nicefrac{1}{2}}D_y^{-\nicefrac{1}{2}} - \nicefrac{1}{2}D_y^{\nicefrac{1}{2}} + \nicefrac{1}{2}    \right)^2  \right)\sum_{z \sim y} \big( u(y) - u(x)  \big)^2  \omega_{yz} \notag \\ & - \left(2m_x\right)^{-1} \sum_{y \sim x} \omega_{xy} m_x^{-\nicefrac{1}{2}}m_y^{-1}D_x^{\nicefrac{1}{2}} \big( u(y) - u(x)   \big)^2 + \nicefrac{1}{2}\left(m_x^{-1} \sum_{y \sim x} \big(  u(y) - u(x)   \big)  \omega_{xy} \right)^2 \notag \\ &\ge \min_{y \sim x} \left(   m^{-1}_xD_x  - \left(  m_x^{-\nicefrac{1}{2}}D_x^{\nicefrac{1}{2}} - \nicefrac{1}{2}D_y  + \nicefrac{1}{2} D_y^{\nicefrac{1}{2}}    \right)^2 - m_x^{-\nicefrac{1}{2}}m_y^{-1}D_x^{\nicefrac{1}{2}} \right)\notag  \\ &\boldsymbol{\cdot} \left(2m_x\right)^{-1} \sum_{y \sim x} \big(   u(y) - u(x)      \big)^2 \omega_{xy} \\ &  +\nicefrac{1}{2}\big( \Delta u (x)   \big)^2 \notag \\ &=\nicefrac{1}{2}\big( \Delta u (x)   \big)^2 + \min_{y \sim x} \left(    m^{-1}_xD_x  - \left(  m_x^{-\nicefrac{1}{2}}D_x^{\nicefrac{1}{2}} - \nicefrac{1}{2}D_y  + \nicefrac{1}{2} D_y^{\nicefrac{1}{2}}    \right)^2 - m_x^{-\nicefrac{1}{2}}m_y^{-1}D_x^{\nicefrac{1}{2}} \right)  \Gamma (u)(x).   \notag 
\end{align*}
\qed
\subsection{ Structural upper bound on the curvature functions}
\subsubsection*{\bf \textit{Proof of Theorem~\ref{thm:main-3}}}
Take $u = \delta_x$, then
\begin{align*}
\Delta u (x) = - \D_x \quad \text{and} \quad \Gamma\left( u \right)(x) = \D_x.
\end{align*}
Setting $X$, $Y$, $a$ and $b$ as before, we can compute
\begin{align}\label{eq:strbd-low}
\nicefrac{1}{2}\, b^{-2}\left( -aX + bY    \right)^2 + b^{-1}XY - XY = \begin{cases}  \nicefrac{1}{2}\, a^2b^{-2} + \left(a - 1\right)b^{-1} + \nicefrac{3}{2} &z = x\\[5pt]
\nicefrac{1}{2}\, a^2b^{-2}& z \neq x.
\end{cases}
\end{align}
Using (\ref{eq:strbd-low}) in (\ref{eq:strbd-Gam2}), we deduce
\begin{align}
\Gamma_2 \left(  \delta_x \right)(x) &= \left(2m_x\right)^{-1}  \sum_{y \sim x} \omega_{xy} \Big( \sum_{z \sim y}  \nicefrac{1}{2}\, a^2b^{-2} \omega_{yz}  \Big) + \left(2m_x\right)^{-1}  \sum_{y \sim x} \omega_{xy}\left( \left(a-1\right)b^{-1} + \nicefrac{3}{2}  \right) \omega_{yx} \notag \\  &\le \left(2m_x\right)^{-1} \sum_{y \sim x} \omega_{xy} \left( \nicefrac{1}{2}\, a^2b^{-2} \D_ym_y \right) +  \left(2m_x\right)^{-1}  \sum_{y \sim x} \omega_{xy}\big( \left(a-1\right)b^{-1} + \nicefrac{3}{2}  \big) \D_ym_y \notag \notag \\  &\le \left(\nicefrac{1}{4} \, m^{-1}_x \max_{y \sim x} m_y \D_x +  \nicefrac{1}{2}\, \max_{y\sim x}   m_y  \D_y^{\nicefrac{1}{2}}\left( m_x^{-\nicefrac{1}{2}}\D_x^{\nicefrac{1}{2}} - 1  \right) + \nicefrac{3}{4} \max_{y\sim x} m_y  \D_y \right)  \D_x    \notag 
\end{align}
which means for all $\N>0$,
\begin{align*}
\K_{G,x}(\N) &\le \K_{G,x}(\infty) \\ &\le \nicefrac{1}{4}\,  m^{-1}_x\D_x \max_{y \sim x} m_y  +  \nicefrac{1}{2}\left( m_x^{-\nicefrac{1}{2}}\D_x^{\nicefrac{1}{2}} - 1  \right) \max_{y\sim x}   m_y  \D_y^{\nicefrac{1}{2}} + \nicefrac{3}{4} \max_{y\sim x} m_y  \D_y.
\end{align*}
\qed
\begin{remark}
	By Theorem~\ref{thm:main-3}, one can deduce curvature bounds for the doubly warped and doubly twisted products of weighted graphs. In practice given a twisted product of weighted networks, one can find the above pointwise bounds via a simple code using the relations
	\begin{align*}
	m_{(z,q)}=m_zm_q \quad \text{and} \quad D_{(z,q)} = \alpha^{-2}(q)D_z + \beta^{-2}(z)D_q. 
	\end{align*}
\end{remark}
Below, we establish curvature bounds for doubly warped products by exploiting the algebraic and geometric properties of quadratic forms arising from Bakry-\'Emery curvature dimension conditions. 
\subsection{ Computation of curvature forms}
\begin{lemma}[$\Delta$ and $\Gamma$]\label{lem:D-G-one} Let $\alpha: G_2 \times G_1 \to \R_+$ and $\beta: G_1 \times G_2 \to \R_+$ be twisting functions. Let $u,v : G_1 \bp G_2 \to \R$ be functions and $u^p,u^x,v^p$ and $v^x$ denote the restrictions of $u$ and $v$ to fibers. Then (suppressing the vertices),
	\begin{align*}
	\Delta u = \alpha^{-2} \Delta^{G_1} u^p +  \beta^{-2} \Delta^{G_2} u^x \quad and \quad \Gamma\left(u,v \right) = \alpha^{-2}\Gamma^{G_1} \left( u^p , v^p \right) +  \beta^{-2} \Gamma^{G_2} \left( u^x , v^x \right).
	\end{align*}
	In particular,
	\begin{align*}
	\Delta \left( u_1 \otimes u_2 \right) = u_2\alpha^{-2} \Delta^{G_1} u_1 + u_1\beta^{-2} \Delta^{G_2} u_2, \quad \Delta \left( u_1 \oplus u_2 \right) = \alpha^{-2} \Delta^{G_1} u_1 + \beta^{-2} \Delta^{G_2} u_2,
	\end{align*}
	\begin{align*}
	\Gamma \left( u_1 \otimes u_2  \right) = u_2^2 \alpha^{-2}\Gamma^{G_1} \left(  u_1 \right) + u_1^2\beta^{-2}\Gamma^{G_2} \left( u_2  \right) \;\; \text{and} \;\; \Gamma \left( u_1 \oplus u_2  \right) =  \alpha^{-2} \Gamma^{G_1} \left(  u_1 \right) + \beta^{-2}\Gamma^{G_2} \left( u_2  \right).
	\end{align*}
\end{lemma}
\begin{proof}
	By definition
	\begin{align}
	\Delta u (x,p) &= \frac{1}{m^{G_1}m^{G_2}}\sum_{(x,p)\sim(y,q)} \big(u(y,q) - u(x,p)\big) \left(  \delta_{xy} m^{G_1} \beta^{-2} \omega^{G_2}_{pq} + \delta_{pq} m^{G_2} \alpha^{-2} \omega^{G_1}_{xy}    \right)  \notag \\ &= \frac{1}{m^{G_2}} \sum_{p \sim q}\big( u(x,q) - u(x,p)\big) \beta^{-2}\omega^{G_2}_{pq} + \frac{1}{m^{G_1}}\sum_{x \sim y}\big( u(y,p) - u(x,p)\big) \alpha^{-2}\omega^{G_1}_{xy} \notag \\ &= \alpha^{-2}\Delta^{G_1} u^p(x)  + \beta^{-2}\Delta^{G_2} u^x(p). \notag
	\end{align}
	Similarly,
	\begin{align}
	\Gamma (u,v) (x,p) &= \frac{1}{2m^{G_1}m^{G_2}}\sum_{(x,p)\sim(y,q)} \big( u(y,q) - u(x,p)\big)\big(v(y,q) - v(x,p)\big) \omega_{\left((x,p)(y,q) \right)} \notag \\ &= \frac{1}{2m^{G_2}}\sum_{p \sim q}\big( u(x,q) - u(x,p)\big) \big( v(x,q) - v(x,p)\big) \beta^{-2}\omega^F_{pq} \notag  \\ &+ \frac{1}{2m^{G_1}}\sum_{x \sim y}\big( u(y,p) - u(x,p) \big) \big( v(y,p) - v(x,p)\big) \alpha^{-2}\omega^B_{xy} \notag \\ &= \alpha^{-2} \Gamma^{G_1} \left(u^p,v^p \right)(x)  + \beta^{-2}\Gamma^{G_2} \left(u^x , v^x \right)(p). \notag
	\end{align}
\end{proof}
\begin{lemma}[First formulation for $\Gamma_2$]
	Let $u , v : G_1 \bp G_2 \to \R$. Then,
	\begin{align}\label{eq:Gam2-1}
	\Gamma_2 \left(u , v \right)  &= \alpha^{-4} \Gamma_2^{G_1} \left( u^p , v^p  \right)(x) + \beta^{-4} \Gamma_2^{G_2} \left( u^x , v^x  \right)(p) \notag \\ & + \; \nicefrac{1}{2}\, \alpha^{-2} \; \mathbf{I} + \nicefrac{1}{2} \, \beta^{-2} \; \mathbf{II},
	\end{align}
	where
	\begin{align*}
	\mathbf{I} =   \Delta^{G_1} \big( \beta^{-2} \Gamma^{G_2}\left(  u^{\bullet} , v^{\bullet} \right)(p) \big) - \Gamma^{G_1} \left(  \beta^{-2} \Delta^{G_2} v^{\bullet}(p) , u^p \right) - \Gamma^{G_1} \left(  \beta^{-2} \Delta^{G_2} u^{\bullet}(p) , v^p \right) \notag,
	\end{align*}
	and
	\begin{align*}
	\mathbf{II} =   \Delta^{G_2} \big( \alpha^{-2} \Gamma^{G_1}\left(  u^{\bullet} , v^{\bullet} \right)(x) \big) - \Gamma^{G_2} \left(  \alpha^{-2} \Delta^{G_1} v^{\bullet}(x) , u^x \right) - \Gamma^{G_2} \left(  \alpha^{-2} \Delta^{G_1} u^{\bullet}(x) , v^x \right). \notag
	\end{align*}
	In particular,
	\begin{align}\label{eq:Gam2-1-2}
	\Gamma_2 \left(u\right) &= \alpha^{-4} \Gamma_2^{G_1} \left( u^p \right) + \beta^{-4} \Gamma_2^{G_2} \left( u^x \right)  \notag \\ & + \; \nicefrac{1}{2}\, \alpha^{-2}  \Big(  \Delta^{G_1} \left( \beta^{-2} \Gamma^{G_2}\left(  u^{\bullet} \right)(p) \right) - 2 \Gamma^{G_1} \left(  \beta^{-2} \Delta^{G_2} u^{\bullet}(p) , u^p \right)   \Big) \\ & +  \; \nicefrac{1}{2}\, \beta^{-2} \Big( \Delta^{G_2} \left( \alpha^{-2} \Gamma^{G_1}\left(  u^{\bullet} \right)(x) \right) - 2 \Gamma^{G_2} \left(  \alpha^{-2} \Delta^{G_1} u^{\bullet}(x) , u^x \right)  \Big). \notag
	\end{align}
\end{lemma}
\begin{proof}
	Based on Lemma~\ref{lem:D-G-one}, one observes
\begin{align*}
\left(\Delta u \right)^p = \alpha^{-2}(p) \Delta^{G_1} u^p + \beta^{-2}  \Delta^{G_2} u^\bullet (p)\quad \text{and}\quad  \left(\Delta u \right)^x = \alpha^{-2} \Delta^{G_1} u^\bullet + \beta^{-2}(x)  \Delta^{G_2} u^x,
\end{align*}
as well as 
	\begin{align*}
	\Gamma \left( u,v \right)^p = \alpha^{-2}(p)\Gamma^{G_1} \left( u^p , v^p \right) +  \beta^{-2} \Gamma^{G_2} \left( u^\bullet , v^\bullet \right)(p)
	\end{align*}
and
\begin{align*}
	\Gamma \left( u,v \right)^x = \alpha^{-2}\Gamma^{G_1} \left( u^\bullet , v^\bullet \right)(x) +  \beta^{-2}(x) \Gamma^{G_2} \left( u^x , v^x \right).
\end{align*}
Thus, according to~\eqref{eq:Gamma2-g}, we get
\begin{align*}
	2 \Gamma_2 \left(u , v \right) &= \Delta \Gamma \left(u , v \right)  -  \Gamma \left( \Delta u , v  \right)  -  \Gamma \left( u , \Delta v   \right) \notag \\
	&= \alpha^{-2}\Delta^{G_1}  \big(\Gamma \left( u,v \right)^p\big) + \beta^{-2}\Delta^{G_2} \big( \Gamma \left( u,v\right)^x \big) \notag \\ &  - \; \alpha^{-2}\Gamma^{G_1} \left( (\Delta u)^p , v^p\right) - \beta^{-2}\Gamma^{G_2} \left( (\Delta u)^x , v^x\right) \notag \\ & - \; \alpha^{-2}\Gamma^{G_1} \left(  u^p , (\Delta v)^p\right) - \beta^{-2}\Gamma^{G_2} \left( u^x , (\Delta v)^x \right) \notag \\ &= \phantom{jj} \alpha^{-2}\Delta^{G_1} \Big( \alpha^{-2} \Gamma^{G_1} \left(  u^p , v^p\right) + \beta^{-2} \Gamma^{G_2}\left(  u^{\bullet} , v^{\bullet} \right)(p) \Big) \notag \\ &  + \; \beta^{-2}\Delta^{G_2} \Big( \alpha^{-2}\Gamma^{G_1} \left(  u^{\bullet} , v^{\bullet} \right)(x) + \beta^{-2} \Gamma^{G_2}\left(  u^x , v^x \right) \Big) \notag \\ &- \; \alpha^{-2}\Gamma^{G_1} \Big( \alpha^{-2} \Delta^{G_1} u^p + \beta^{-2} \Delta^{G_2} u^{\bullet}(p), v^p\Big) \notag \\ &- \; \beta^{-2}\Gamma^{G_2} \Big( \alpha^{-2} \Delta^{G_2} u^x + \beta^{-2} \Delta^{G_2} u^{\bullet}(x), v^x\Big) \notag \\ &- \; \alpha^{-2}\Gamma^{G_1} \Big( \alpha^{-2} \Delta^{G_1} v^p + \beta^{-2} \Delta^{G_2} v^{\bullet}(p), u^p\Big) \notag \\ &- \; \beta^{-2}\Gamma^{G_2} \Big( \alpha^{-2} \Delta^{G_2} v^x + \beta^{-2} \Delta^{G_2} v^{\bullet}(x), u^x\Big) \notag \\ &= \blue{\alpha^{-4} \Delta^{G_1} \Gamma^{G_1} \left(  u^p , v^p\right)} + \alpha^{-2} \Delta^{G_1} \big( \beta^{-2} \Gamma^{G_2} \left(  u^{\bullet} , v^{\bullet} \right)(p) \big) \notag \\ &\red{+  \;\beta^{-4} \Delta^{G_2} \Gamma^{G_2} \left(  u^x , v^x\right)} + \beta^{-2} \Delta^{G_2} \big( \alpha^{-2} \Gamma^{G_1} \left(  u^{\bullet} , v^{\bullet} \right)(x) \big) \notag \\ & \blue{- \; \alpha^{-4} \Gamma^{G_1} \left( \Delta^{G_1} u^p , v^p\right)} - \alpha^{-2} \Gamma^{G_1} \left(  \beta^{-2} \Delta^{G_2} u^{\bullet}(p) , v^p \right) \notag \\ & \red{- \; \beta^{-4} \Gamma^{G_2} \left( \Delta^{G_2} u^x , v^x\right)} - \beta^{-2} \Gamma^{G_2} \left(  \alpha^{-2} \Delta^{G_2} u^x , v^x \right)(p) \notag \\ & \blue{- \; \alpha^{-4} \Gamma^{G_1} \left( \Delta^{G_1} v^p , u^p\right)} - \alpha^{-2} \Gamma^{G_1} \left(  \beta^{-2} \Delta^{G_2} v^{\bullet}(p) , u^p \right) \notag \\ & \red{- \; \beta^{-4} \Gamma^{G_2} \left( \Delta^{G_2} v^x , u^x\right)} - \beta^{-2} \Gamma^{G_2} \left(  \alpha^{-2} \Delta^{G_2} v^x , u^x \right)  \notag \\ &= \blue{2 \alpha^{-4} \Gamma_2^{G_1} \left( u^p , v^p  \right)}  \red{+ 2 \beta^{-4} \Gamma_2^{G_2} \left( u^x , v^x  \right)}  + \alpha^{-2} \; \mathbf{I} + \beta^{-2} \; \mathbf{II}. \notag 
	\end{align*}
	Notation: $\bullet$ is used as a dummy variable e.g. $u^\bullet$ denotes the restriction of $u$ to the $\bullet$-fiber.
\end{proof}
\begin{lemma}
	For $U_{c_1,c_2} := c_1f_1 \oplus c_2f_2$ we thus get
	\begin{align*}
	\Gamma_2 \left(U_{c_1,c_2}\right) = c_1^2 \alpha^{-4} \Gamma_2^{G_1} \left( f_1 \right) + c_2^2 \beta^{-4} \Gamma_2^{G_2} \left( f_2 \right) + \mathcal{Q}\left( c_1 , c_2 \right)
	\end{align*}
	where
	\begin{align}\label{eq:def-Q}
	\mathcal{Q}\left( c_1 , c_2 \right) &=  \nicefrac{1}{2}\, c_2^2 \alpha^{-2}   \Gamma^{G_2} \left(  f_2 \right)  \Delta^{G_1}  \beta^{-2}  -  c_1c_2 \alpha^{-2} \Delta^{G_2} f_2 \Gamma^{G_1} \left(  \beta^{-2}  , f_1 \right)  \notag \\ & +  \nicefrac{1}{2} \, c_1^2 \beta^{-2}   \Gamma^B\left(  f_1 \right)  \Delta^{G_2}  \alpha^{-2}  -  c_1c_2 \beta^{-2} \Delta^{G_1} f_1 \Gamma^{G_2} \left(  \alpha^{-2}  , f_2 \right). 
	\end{align}
\end{lemma}
\begin{proof}
	We have
\begin{align*}
U_{c_1,c_2}^p = c_1f_1 + c_2f_2(p), \quad \text{and}, \quad U_{c_1,c_2}^x = c_1f_1(x) + c_2f_2.
\end{align*}
Hence, by (\ref{eq:Gam2-1-2}), we get
\begin{align}
		\Gamma_2 \left(U_{c_1,c_2}\right) &= \alpha^{-4} \Gamma_2^{G_1} \big( c_1f_1 + c_2f_2(p) \big) + \beta^{-4} \Gamma_2^{G_2} \big( c_1f_1(x) + c_2f_2 \big)  \notag \\ & + \; \nicefrac{1}{2}\, \alpha^{-2}  \Big(  \Delta^{G_1} \big( \beta^{-2} \Gamma^{G_2}\left(  U_{c_1,c_2}^{\bullet} \right)(p) \big) - 2 \Gamma^{G_1} \big(  \beta^{-2} \Delta^{G_2} U_{c_1,c_2}^{\bullet}(p) , U_{c_1,c_2}^p \big)   \Big) \notag \\ & +  \; \nicefrac{1}{2}\, \beta^{-2} \Big( \Delta^{G_2} \big( \alpha^{-2} \Gamma^{G_1}\left(  U_{c_1,c_2}^{\bullet} \right)(x) \big) - 2 \Gamma^{G_2} \left(  \alpha^{-2} \Delta^{G_1} U_{c_1,c_2}^{\bullet}(x) , U_{c_1,c_2}^x \right)  \Big). \notag \\ &= c_1^2 \alpha^{-4} \Gamma_2^{G_1} \left( f_1 \right) + c_2^2 \beta^{-4} \Gamma_2^{G_2} \left( f_2 \right) \notag \\ &+  \; \nicefrac{1}{2}\, \alpha^{-2}  \Big(  \Delta^{G_1} \big[ \beta^{-2} \Gamma^{G_2}\left( c_2f_2 \right)(p) \big] - 2 \Gamma^{G_1} \left(  \beta^{-2} \Delta^{G_2} c_2f_2 (p) , c_1f_1 \right)   \Big)  \notag \\ & +  \; \nicefrac{1}{2}\, \beta^{-2} \Big( \Delta^{G_2} \big[ \alpha^{-2} \Gamma^{G_1}\left( c_1f_1 \right)(x) \big] - 2 \Gamma^{G_2} \left(  \alpha^{-2} \Delta^{G_1} c_1f_1 (x) , c_2f_2 \right)  \Big). \notag \\ &= c_1^2 \alpha^{-4} \Gamma_2^{G_1} \left( f_1 \right) + c_2^2 \beta^{-4} \Gamma_2^{G_2} \left( f_2 \right) \notag \\ &
\left. \begin{array}{rr} 
	+ \; \nicefrac{1}{2}\, \alpha^{-2}  \Big(  c_2^2 \Delta^{G_1}\left(\beta^{-2}\right) \Gamma^{G_2}\left(f_2 \right)- 2c_1c_2 \Delta^{G_2}\left(f_2\right) \Gamma^{G_1} \left(  \beta^{-2} , f_1 \right)   \Big) & \\[10pt]  +  \; \nicefrac{1}{2}\, \beta^{-2} \Big( c_1^2 \Delta^{G_2}\left(\alpha^{-2}  \right) \Gamma^{G_1}\left(f_1 \right) - 2 c_1c_2 \Delta^{G_1}\left( f_1\right) \Gamma^{G_2} \left(  \alpha^{-2}  , f_2 \right)  \Big) &
\end{array} \right\} 
	\mathcal{Q}\left( c_1 , c_2 \right). \notag
\end{align}
\end{proof}
\begin{lemma}[Second formulation for $\Gamma_2$]\label{lem:sec-form-Gamma2}
	For the special case $u = u_1 \otimes u_2$ and $v = v_1 \otimes v_2$ where $u_1,v_1:G_1 \to \R$ and $u_2,v_2: G_2 \to \R$, we have
	\begin{align}\label{eq:Gam2-2}
	\Gamma_2 \left(u_1 \otimes u_2 , v_1 \otimes v_2 \right) &=  u_2v_2 \alpha^{-4} \Gamma^{G_1}_2\left(  u_1,v_1 \right) + u_1v_1 \beta^{-4} \Gamma^{G_2}_2\left(  u_2,v_2 \right) \notag \\ & + \nicefrac{1}{2}\, \alpha^{-2} \mathbf{I}  + \nicefrac{1}{2}\, \beta^{-2} \mathbf{II},
	\end{align}	
	where,
	\begin{align*}
	\mathbf{I} := \Gamma^{G_2}(u_2,v_2) \Delta^{G_1} \left( u_1v_1\beta^2 \right) - v_2 \Delta^{G_2}u_2\Gamma^{G_1}\left(u_1\beta^{-2} , v_1  \right) - u_2 \Delta^{G_2}v_2\Gamma^{G_1}\left(v_1\beta^{-2} , u_1  \right)
	\end{align*}
	and
	\begin{align*}
	\mathbf{II} := \Gamma^{G_1}(u_1,v_1) \Delta^{G_2} \left( u_2v_2\alpha^{-2} \right) - v_1 \Delta^{G_1}u_1\Gamma^{G_2}\left(u_2\alpha^{-2} , v_2  \right) - u_1 \Delta^{G_1}v_1\Gamma^{G_2}\left(v_2\alpha^{-2} , u_2  \right).
	\end{align*}
	In particular,
	\begin{align*}
	\Gamma_2 \left(u_1 \otimes u_2 \right) =  u_2^2 \alpha^{-4} \Gamma^{G_1}_2\left(  u_1 \right) + u_1^2 \beta^{-4} \Gamma^{G_2}_2\left(  u_2\right) + \nicefrac{1}{2}\, \alpha^{-2} \; \mathbf{I}  + \nicefrac{1}{2}\, \beta^{-2} \; \mathbf{II},
	\end{align*}
	where,
	\begin{align*}
	\mathbf{I} := \Gamma^{G_2}(u_2) \Delta^{G_1} \left( u_1^2\beta^{-2} \right) - 2 u_2 \Delta^{G_2}u_2\Gamma^{G_1}\left(u_1\beta^{-2} , u_1  \right),
	\end{align*}
	and
	\begin{align*}
	\mathbf{II} := \Gamma^{G_1}(u_1) \Delta^{G_2} \left( u_2^2\alpha^{-2} \right) - 2 u_1 \Delta^{G_1}u_1\Gamma^{G_2}\left(u_2\alpha^{-2} , u_2  \right).
	\end{align*}
\end{lemma}	
\begin{proof}
One clearly has
\begin{align*}
\left(u_1 \otimes u_2\right)^p = u_2(p)u_1, \quad \text{and}\quad \left(u_1 \otimes u_2\right)^x = u_1(x)u_2,
\end{align*}
and similar statements also hold for $v_1\otimes v_2$. 
\par Therefore, by Lemma~\ref{lem:D-G-one}, we get
\begin{align*}
\Delta \left( u_1 \otimes u_2\right) = \alpha^{-2}u_2 \Delta^{G_1} u_1 +  \beta^{-2}u_1 \Delta^{G_2} u_2 = \Delta^{G_1} u_1 \otimes \alpha^{-2}u_2 +  \beta^{-2}u_1 \otimes \Delta^{G_2} u_2,
\end{align*}
and 
	\begin{align}
	\Gamma\left(u_1 \otimes u_2,v_1 \otimes v_2 \right) &= \alpha^{-2}\Gamma^{G_1} \left( \left(u_1 \otimes u_2\right)^p , \left(v_1 \otimes v_2\right)^p \right) +  \beta^{-2} \Gamma^{G_2} \left( \left(u_1 \otimes u_2\right)^x , \left(v_1 \otimes v_2\right)^x \right) \notag \\ &= \alpha^{-2}u_2v_2\Gamma^{G_1} \left(u_1, v_1 \right) + \beta^{-2}u_1v_1\Gamma^{G_2} \left(u_2,v_2 \right) \notag \\ &= \Gamma^{G_1} \left(u_1, v_1 \right)  \otimes \alpha^{-2}u_2v_2 + \beta^{-2}u_1v_1 \otimes  \Gamma^{G_2} \left(u_2,v_2 \right). \notag
	\end{align}
Consequently,
\begin{align*}
	2 \Gamma_2 \left(u_1 \otimes u_2 , v_1 \otimes v_2 \right) &= \Delta \Gamma \left(u_1 \otimes u_2 , v_1 \otimes v_2 \right) \notag \\ & - \;  \Gamma \left( \Delta \left(u_1 \otimes u_2\right) , v_1 \otimes v_2   \right) \\ & - \; \Gamma \left( \left(u_1 \otimes u_2\right) , \Delta \left(v_1 \otimes v_2\right)   \right) \notag \\
	&= \Delta \big( \Gamma^{G_1}\left(u_1 , v_1   \right) \otimes u_2v_2\alpha^{-2} \big) + \Delta \big( u_1v_1\beta^{-2} \otimes \Gamma^{G_2}\left(u_2 , v_2   \right) \big)  \notag \\ & - \; \Gamma \big( \Delta^{G_1} u_1 \otimes u_2\alpha^{-2}, v_1 \otimes v_2   \big) - \Gamma \big( u_1\beta^{-2} \otimes \Delta^{G_2} u_2, v_1 \otimes v_2   \big) \notag \\ & - \Gamma \big( u_1 \otimes u_2   , \Delta^{G_1} v_1 \otimes v_2\alpha^{-2}  \big) - \Gamma \big( u_1 \otimes u_2 ,  v_1\beta^{-2} \otimes \Delta^{G_2} v_2   \big) \notag \\ &= \blue{\alpha^{-4}u_2v_2\Delta^{G_1}\Gamma^{G_1}(u_1,v_1)} + \beta^{-2} \Gamma^{G_1}(u_1,v_1)\Delta^{G_2}(u_2v_2\alpha^{-2}) \notag \\ &  + \;\alpha^{-2}\Gamma^{G_2}(u_2,v_2) \Delta^{G_1} \left( u_1v_1\beta^{-2} \right) \red{+ u_1v_1\beta^{-4} \Delta^{G_2} \Gamma^{G_2} \left( u_2,v_2 \right)} \notag \\ & \blue{- \;  \alpha^{-4} u_2v_2 \Gamma^{G_1} \left(\Delta^{G_1} u_1 , v_1  \right)} -  \beta^{-2}v_1\Delta^{G_1}u_1\Gamma^{G_2}\left(u_2\alpha^{-2},v_2   \right) \notag \\ & - \;  \alpha^{-2} v_2 \Delta^{G_2}u_2\Gamma^{G_1}\left(u_1\beta^{-2} , v_1  \right) \red{- \beta^{-4} u_1v_1 \Gamma^{G_2}\left( \Delta^{G_2} u_2 , v_2  \right)} \notag \\ & \blue{- \; \alpha^{-4} u_2 v_2 \Gamma^{G_1}\left(  u_1 , \Delta^{G_1} v_1 \right)} - \beta^{-2}u_1\Delta^{G_1}v_1\Gamma^{G_2}\left(v_2\alpha^{-2},u_2   \right) \notag \\ & - \;  \alpha^{-2} u_2 \Delta^{G_2}v_2\Gamma^{G_1}\left(v_1\beta^{-2} , u_1  \right) \red{- \beta^{-4} u_1v_1 \Gamma^{G_2}\left( \Delta^{G_2} v_2 , u_2  \right)} \notag  \\ &=   \blue{2\alpha^{-4}u_2v_2\Gamma^{G_1}_2\left(  u_1,v_1 \right)}  \red{+2u_1v_1 \beta^{-4} \Gamma^{G_2}_2\left(  u_2,v_2 \right)} + \alpha^{-2} \mathbf{I} + \beta^{-2}\mathbf{II}.\notag
	\end{align*}
	When $u_i = v_i$, this simplifies to
	\begin{align*}
	2 \Gamma_2 \left(u_1 \otimes u_2 \right)  &= 2\alpha^{-4}u_2^2\Gamma^{G_1}_2\left(  u_1,v_1 \right) +2u_1^2 \beta^{-4} \Gamma^{G_2}_2\left(  u_2,v_2 \right)   \notag \\ &  + \alpha^{-2}\Big( \Gamma^{G_2}(u_2) \Delta^{G_1} \left( u_1^2\beta^{-2} \right) - 2 u_2 \Delta^{G_2}\left(u_2\right)\Gamma^{G_1}\left(u_1\beta^{-2} , u_1  \right)   \Big) \\ & + \beta^{-2}\Big( \Gamma^{G_1}(u_1)\Delta^{G_2}(u_2^2\alpha^{-2}) - 2u_1\Delta^{G_1}\left(u_1\right)\Gamma^{G_2}\left(u_2\alpha^{-2},u_2   \right) \Big)\notag.
\end{align*}	
\par Note that one can also obtain (\ref{eq:Gam2-2}) directly from (\ref{eq:Gam2-1}). 
\end{proof}
\subsection{ Some useful lemmas}
Obtaining upper bounds for the curvature functions is a max-min problem hence we need to use some intersection theoretic facts from algebraic/differential geometry to surfaces which are obtained from curvature maximizers. 
\subsubsection*{\bf \textit{Quadratic surfaces}}
\par Consider the quadratic surface
\begin{align*}
\Sigma: \;\; z= ax^2 + by^2 + cxy.  
\end{align*}
in $\R^3$. By standard surface theory (see e.g. \cite{MP}), the principal curvatures of $\Sigma$ are given by
\begin{align*}
\kappa_{i} = \nicefrac{1}{2}\, \left(a + b \pm \sqrt{\left(   a - b  \right)^2 + c^2}  \right) \quad i = 1,2 \quad \text{and}\quad \kappa_1 \le \kappa_2.
\end{align*}
The principal directions of $\Sigma$ are counterclockwise rotations of the $x$ and $y$ axes by $\theta := \nicefrac{1}{2}\arctan \nicefrac{c}{(a-b)}$ where $\theta \in \left[ 0, \pi \right]$. Here, the direction of $\kappa_{1}$ (which is either $\theta$ or $\theta \pm \frac{\pi}{2}  \in \left[ 0 , \pi \right] $) is called the principal angle of $\Sigma$. The following geometric classification of quadratic surfaces (based on the sign of (Gauss) curvature $\kappa_1\kappa_2$) is standard. 
\begin{enumerate}
	\item Parabolic cylinders:  $\kappa_1\kappa_2 = 0$ and $\kappa_1 + \kappa_2 \neq 0$;
	\item Hyperbolic-paraboloids: $\kappa_{1} \kappa_2 < 0$;
	\item Paraboloids: $\kappa_1\kappa_2 > 0$.
\end{enumerate}
\begin{lemma}
	Surfaces $\Sigma_1$ and $\Sigma_2$ given by
	\begin{align*}
	\Sigma_1: \;\; z_1 = a_1x^2 + b_1y^2 + c_1xy    \quad \text{and} \quad \Sigma_2: \;\; z_2 = a_2x^2 + b_2y^2 + c_2xy,
	\end{align*}
	have at least a line worth of non-trivial intersection if and only if
	\begin{align*}
	- \det \nabla^2 \left( z_2 - z_1  \right) = \left( c_2 - c_1  \right)^2 - 4 \left( a_2 - a_1 \right)\left( b_2 - b_1   \right) \ge 0. 
	\end{align*}
\end{lemma}
\begin{proof}
The proof can be done via direct calculations. However, an intuitive geometric way to see this, is to notice the surface $z_2-z_1$ passes through origin with Gaussian curvature $\det \nabla^2 \left( z_2 - z_1  \right)$. When the Gaussian curvature is non-positive, the surface $z_2-z_1$ (sans origin) near the origin, can not only live on one side of the $x-y$ plane; this can be seen by writing the Taylor series for $z_2-z_1$ and is more or less standard. it follows that $z_2-z_1 = 0$ at some point other than $(0,0)$. Since the intersection of two homogeneous polynomials is a variety, we deduce there exist at least a line worth of zeros. 
\end{proof}
The following rudimentary lemma comes in handy later on. 
\begin{lemma}\label{lem:F-zeros}
	Let $A$ and $B$ be constant numbers and $\M, \N_1, \N_2 >0$. The surface
\begin{align*}
\Sigma : \; z =\left(\N_1^{-1} - \M^{-1}\right)  A^2x^2 + \left(\N_2^{-1} - \M^{-1}\right)   B^2 y^2 - 2\M^{-1} \;  AB xy,
\end{align*}
with both $A$ and $B$ not simultaneously zero, is
\begin{enumerate}
\item A praboloid when $\M > \N_1 + \N_2$;
\item A parabolic cylinder when $\M = \N_1 + \N_2$;
\item A hyperbolic-paraboloid when $\M < \N_1 + \N_2$. 
\end{enumerate}
In particular, 
\begin{enumerate}
\item For any $A$ and $B$, 
\begin{align*}
	\Sigma : \; z =\N_1^{-1}\N_2 \left( \N_1 + \N_2 \right)^{-1}  A^2x^2 + \N_1\N_2^{-1} \left( \N_1 + \N_2 \right)^{-1}  B^2 y^2 - 2 \left( \N_1 + \N_2\right)^{-1}  AB xy,
\end{align*}
is either a parabolic cylinder or the $x-y$ plane. 
\item When $\M \le \N_1 + \N_2$, the variety $z = 0$ has at least dimension 1; i.e. $z=0$ has at least a line worth of solutions. 
\end{enumerate}
\end{lemma}
\begin{proof}
By direct calculation, 
\begin{align*}
	\det \nabla^2 z &= 2\left(\N_1^{-1} - \M^{-1}  \right)  A^2 \boldsymbol{\cdot} 2\left(\N_2^{-1} - \M^{-1}  \right)  B^2 - \left( 2\M^{-1}  AB\right)^2 \notag \\ &= 
	4\M^{-1}\N_1^{-1}\N_2^{-1}\left( \M - \N_1 - \N_2 \right)A^2B^2,
\end{align*}	
so depending on the sign of $\M - \N_1-\N_2$, one gets one of the above cases (the plane is excluded since $A$ and $B$ are not both equal to zero). 
\par In the special case, $\M = \N_1 + \N_2$, we have $ \det \nabla^2 z = 0 $. Furthermore, $z=0$ has exactly one line of zeros when $A$ and $B$ do not vanish simultaneously. If $A,B \neq 0$, the line of zeros is $y =  \frac{\N_2 A}{\N_1 B} x$. If $A \neq 0$ and $B = 0$, $x = 0$ and if $A = 0$ and $B \neq 0$, $y = 0$ are the lines of zeros. 
\par In the case $\M < \N_1 + \N_2$, the variety $z = 0$ is singular and one dimensional in the following cases:
\begin{enumerate}
\item two intersecting lines
\begin{enumerate}
	\item [(1-1)]
	\begin{align*}
		y = \left( \M - \N_2  \right)^{-1} \Big(  \N_2 \pm \big( \M \left( \N_1 + \N_2 - \M  \right)  \big)^{\nicefrac{1}{2}} \Big)  ABx, \quad \M \neq \N_2,
	\end{align*}
	when $A,B \neq 0$ and $\M \neq \N_2$. 
\item [(1-2)]
\begin{align*}
	y = \left( 2\N_1  \right)^{-1}\left( \M - \N_1 \right) AB^{-1}x \quad \text{and}\quad x = 0,
\end{align*}
when $\M = \N_2 $ and $A, B \neq 0$. 
\item [(1-3)]
\begin{align*}
	y = \left( 2\N_1  \right)^{-1}\left( \M - \N_1 \right) AB^{-1}x \quad \text{and}\quad y = 0,
\end{align*}
when $\M = \N_1 $ and $A, B \neq 0$. 
\end{enumerate}
\item a double line
\begin{enumerate}
	\item [(2-1)]
	\begin{align*}
		y = 0 \quad B\neq 0, A=0, \M \neq \N_2,
	\end{align*}
\item [(2-2)]
\begin{align*}
	x = 0 \quad B = 0, A\neq 0, \M \neq \N_1.
\end{align*}
\end{enumerate}
\end{enumerate}
\par The intersection is two dimensional (the trivial $x-y$ plane) if either
\begin{align*}
A = 0, \M = \N_2 \quad \text{or} \quad B = 0, \M = \N_1;
\end{align*}
it goes without saying that the role of $x$ and $y$ is symmetric i.e. writing the intersection as equations of $x$ in terms of $y$ would lead the exact same cases. 
\end{proof}
\subsubsection*{ \bf \textit{Quadratic estimates}}
Obtaining lower bounds for the curvature functions and also various upper bounds, will require estimating certain quadratic forms in terms of the desired quadratic forms which are the square of Laplacian and the gradient. For this , we will need and frequently use the following useful estimates. 
\begin{lemma}\label{lem:est}
	The inequalities 
	\begin{enumerate}
		\item $	\left( \Delta^{G_i} f \right)^2  \le 2\Deg_{G_i} \; \Gamma^{G_i} \left( f  \right)$,
		\item $ \big( \Gamma^{G_i} \left(  f  ,  g \right) \big)^2  \le \Gamma^{G_i} \left( f  \right) \; \Gamma^{G_i} \left( g  \right)$,
		\item $\left| \Delta^{G_1} f \Gamma^{G_2} \left(  g  , h \right) \right| \le \nicefrac{1}{2}\, \Deg_{G_1} \; \Gamma^{G_1} \left( f  \right) + \nicefrac{1}{2}\, \Gamma^{G_2} \left( g  \right) \; \Gamma^{G_2} \left( h  \right)  $,
	\end{enumerate}
	hold on $G_1$ and $G_2$.
\end{lemma}
\begin{proof}
	By Cauchy-Schwarz, for a vertex $z$ in any weighted graph 
	\begin{align}
	\big( \Delta f (z)\big)^2 &= \left( m(z)^{-\nicefrac{1}{2}} \sum_{w \sim z} \big( f(w)   - f(z)   \big) \left( \omega_{zw} \right)^{\nicefrac{1}{2}} m(z)^{-\nicefrac{1}{2}} \omega_{zw}^{\nicefrac{1}{2}} \right)^2 \notag \\ &\le 2\left(  \left( 2m(z) \right)^{-1}\sum_{w \sim z} \big( f(w)   - f(z)   \big)^2  \omega_{zw} \right) \left( m(z)^{-1} \sum_{w \sim z} \omega_{zw} \right) \notag \\ &= 2\Deg(z) \; \Gamma \left( f  \right)(z)\notag,
	\end{align}
	and
	\begin{align}
	\big( \Gamma \left(  f,g \right)(z) \big)^2 &= \nicefrac{1}{4} \left(\sum_{w \sim z} m(z)^{-\nicefrac{1}{2}}  \big( f(w)   - f(z)  \big)\omega_{zw}^{\nicefrac{1}{2}} \sqdot m(z)^{-\nicefrac{1}{2}} \big( g(w)   - g(z)  \big) \omega_{zw}^{\nicefrac{1}{2}} \right)^2 \notag \\ &\le \left( \left(2m(z)\right)^{-1} \sum_{w \sim z} \big( f(w)  - f(z) \big)^2 \omega_{zw} \right)\left( \left(2m(z)\right)^{-1} \sum_{w \sim z}  \big( g(w)   - g(z)  \big)^2 \omega_{zw} \right) \notag \\ &= \Gamma \left( f \right)(z) \; \Gamma \left( g  \right)(z).\notag
	\end{align}
	Applying the Young's inequality,
	\begin{align}
	\left| \Delta^{G_1} f (x)\Gamma^{G_2} \left(  g  , h  \right)(p) \right| &\le  \nicefrac{1}{2}\left( \Delta^{G_1} f(x) \right)^2 + \nicefrac{1}{2}\, \big( \Gamma^{G_2} \left( g  , h \right)(p) \big)^2 \notag \\ &\le \Deg_{G_1}(x) \; \Gamma^{G_1} \left( f  \right)(x) + \nicefrac{1}{2}\, \Gamma^{G_2} \left( g \right)(p) \; \Gamma^{G_2} \left( h  \right)(p). \notag
	\end{align}
\end{proof}
\begin{lemma}
	The quadratic form, $\mathcal{Q}\left( c_1 , c_2  \right)$, given in (\ref{eq:def-Q}), can be bounded as
	\begin{align*}
	\mathcal{Q}\left( c_1, c_2  \right) \le  \mathcal{Q}_1\left(  c_1 , c_2 \right)\Gamma^{G_1}(f_1) +  \mathcal{Q}_2\left(  c_1 , c_2 \right) \Gamma^{G_2}(f_2),
	\end{align*}
	where 
	\begin{align*}
	\mathcal{Q}_1\left(  c_1 , c_2 \right) =  \begin{cases} 	\mathcal{Q}_{11}:= \nicefrac{1}{2}\, c_1^2 \beta^{-2} \Delta^{G_2} \alpha^{-2} + \left|c_1 \right|\left|c_2 \right| \beta^{-2} \D_x; & \Delta^{G_1}f_1 , \Delta^{G_2}f_2 \neq 0 \\ \phantom{sajjad}+\nicefrac{1}{2} \left|c_1 \right|\left|c_2 \right| \alpha^{-2}\Gamma^{G_1} \left( \beta^{-2}  \right) \\[10pt]  \mathcal{Q}_{12}:= \nicefrac{1}{2} c_1^2 \beta^{-2} \Delta^{G_2} \alpha^{-2} +\frac{1}{2} \left|c_1 \right|\left|c_2 \right| \alpha^{-2}\Gamma^{G_1} \left( \beta^{-2}  \right) ; & \Delta^{G_1}f_1 =0 , \Delta^{G_2}f_2 \neq 0 \\[10pt]
\mathcal{Q}_{13}:= \nicefrac{1}{2}\, c_1^2 \beta^{-2} \Delta^{G_2} \alpha^{-2} + \left|c_1 \right|\left|c_2 \right| \beta^{-2} \D_x; & \Delta^{G_1}f_1 \neq 0 , \Delta^{G_2}f_2 = 0 \\[10pt]  \mathcal{Q}_{14}:= \nicefrac{1}{2}\, c_1^2 \beta^{-2} \Delta^{G_2} \alpha^{-2};  & \Delta^{G_1}f_1 , \Delta^{G_2}f_2 = 0 \end{cases} \notag 
	\end{align*}
	and
	\begin{align*}
	\mathcal{Q}_2\left(  c_1 , c_2 \right) =  \begin{cases} \mathcal{Q}_{21}:= \nicefrac{1}{2}\, c_2^2 \alpha^{-2} \Delta^{G_1} \beta^{-2} + \left|c_1 \right|\left|c_2 \right| \alpha^{-2} \D_p; & \Delta^{G_1}f_1 , \Delta^{G_2}f_2 \neq 0 \\ \phantom{sajjad} +\nicefrac{1}{2} \left|c_1 \right|\left|c_2 \right| \beta^{-2}\Gamma^{G_2} \left( \alpha^{-2}  \right) \\[10pt] \mathcal{Q}_{22}:= \nicefrac{1}{2}\, c_2^2 \alpha^{-2} \Delta^{G_1} \beta^{-2} +\nicefrac{1}{2} \left|c_1 \right|\left|c_2 \right| \beta^{-2}\Gamma^{G_2} \left( \alpha^{-2}  \right)  ; & \Delta^{G_1}f_1 \neq 0, \Delta^{G_2}f_2 =0\\[10pt]
	\mathcal{Q}_{23}:= \nicefrac{1}{2}\, c_2^2 \alpha^{-2} \Delta^{G_1} \beta^{-2} + \left|c_1 \right|\left|c_2 \right| \alpha^{-2} \D_p; & \Delta^{G_1}f_1 = 0 , \Delta^{G_2}f_2 \neq 0 \\[10pt]  \mathcal{Q}_{24}:= \nicefrac{1}{2}\, c_2^2 \alpha^{-2} \Delta^{G_1} \beta^{-2};  & \Delta^{G_1}f_1 , \Delta^{G_2}f_2 = 0. \end{cases} \notag 
	\end{align*}	
\par	Furthermore,
	\begin{align*}
	\mathcal{Q}_{1i}\left(1,0\right) = \mathcal{Q}_{1}\left( 1,0  \right) \quad \text{and} \quad \mathcal{Q}_{2i}\left(0,1\right) = \mathcal{Q}_2\left( 0,1  \right).
	\end{align*}
\end{lemma}
\begin{proof}
	Applying the estimates from Lemma~\ref{lem:est}, when $\Delta^{G_1}f_1 , \Delta^{G_2}f_2 \neq 0$,
	\begin{align}
	\mathcal{Q}\left( c_1, c_2  \right) &=  \nicefrac{1}{2}\, c_2^2 \alpha^{-2}(p)   \Gamma^{G_2}\left(  f_2 \right)(p)  \Delta^{G_1}  \beta^{-2}(x)  -  c_1c_2 \alpha^{-2}(p) \Delta^{G_2} f_2 (p)\Gamma^{G_1} \left(  \beta^{-2}  , f_1 \right)(x)  \notag \\ & + \; \nicefrac{1}{2}\, c_1^2 \beta^{-2}(x)   \Gamma^{G_1}\left(  f_1 \right)(x)  \Delta^{G_2}  \alpha^{-2}(p)  - c_1c_2 \beta^{-2}(x) \Delta^{G_1} f_1 (x)\Gamma^{G_2} \left(  \alpha^{-2}  , f_2 \right)(p) \notag  \\ &\le   \Big(  \nicefrac{1}{2}\, c_1^2 \beta^{-2} \Delta^{G_2} \alpha^{-2} + \left|c_1 \right|\left|c_2 \right| \beta^{-2} \D_x +  \nicefrac{1}{2}\, \left|c_1 \right|\left|c_2 \right| \alpha^{-2}\Gamma^{G_1} \left( \beta^{-2}  \right) \Big) \Gamma^{G_1} \left( f_1  \right) \notag \\ & + \Big( \nicefrac{1}{2}\, c_2^2 \alpha^{-2} \Delta^{G_1}  \beta^{-2} + \left|c_1 \right|\left|c_2 \right| \alpha^{-2} \D_p + \nicefrac{1}{2}\, \left|c_1 \right|\left|c_2 \right| \beta^{-2}\Gamma^{G_2} \left( \alpha^{-2}  \right) \Big) \Gamma^{G_2} \left( f_2  \right). \notag
	\end{align}
	The other cases follow similarly. 
\end{proof}
\subsubsection*{\bf \textit{Curvature maximizers}}
\begin{lemma}[Pointwise curvature maximizers]
	For any graph $G$ and any vertex $x$, and any $\N \in (0,\infty]$, one can find a function $f_{\N,x}$ such that the identity
\begin{align*}
\Gamma_2 \left( f_{\N,x} \right)(x) = \N^{-1}\big( \Delta^{G}f_{\N,x}(x) \big)^2 + \K_{G,x}\Gamma^{G} \left( f_{\N,x}  \right)(x),
\end{align*}	
holds true. 
\end{lemma}
\begin{proof}
	The proof is an elementary limiting argument. For any $\epsilon>0$ and based on the optimality of $\K_{G,x}$, one can find \emph{non-trivial} functions $\bar{f}^{\epsilon}_{\N,x}$ such that
	\begin{align*}
 \Gamma_2 \left( \bar{f}_{\N,x} \right)(x) < \N^{-1}\big( \Delta^{G} \bar{f}_{\N,x}(x) \big)^2 + \left( \K_{G,x} + \epsilon \right)\Gamma^{G} \left( \bar{f}_{\N,x}  \right)(x),
	\end{align*}
holds. Set $f^{\epsilon}_{\N,x} := \frac{\bar{f}^{\epsilon}_{\N,x}}{\max |\bar{f}^{\epsilon}_{\N,x}|}$, then 
\begin{align*}
\Gamma_2 \left( f^{\epsilon}_{\N,x} \right)(x) < \N^{-1}\big( \Delta^{G} f^{\epsilon}_{\N,x}(x) \big)^2 + \left( \K_{G,x} + \epsilon \right)\Gamma^{G} \left( f^{\epsilon}_{\N,x}  \right)(x), \quad \| f^{\epsilon}_{\N,x} \|_{\sup} = 1.
\end{align*}
Now by, say a diagonal argument, one can pick a sequence $\epsilon_{j} \downarrow 0$, such that $f^{\epsilon_j}_{\N,x}$ convergence at all the vertices to a function $f_{\N,x}$. This limit function is a curvature maximizer at $x$. In an infinite graph, one can first truncate the $f^{\epsilon}_{\N,x}$ outside a ball of radius $3$ and then do the limiting argument; notice that the $\N$-Bakry-\'Emery curvature at $x$ only sees the graph data in a ball of radius $2$ around $x$.  
\end{proof}
\subsection{ Upper bounds for curvature functions}
\subsubsection*{\bf \textit{Proof of Theorem~\ref{thm:main-5}}}
Let $f_i:G_i \to \R$, $i = 1,2$ be curvature maximizers at $x$ and $p$ respectively i.e.
\begin{align*}
\Gamma_2^{G_i} \left(  f_i \right)\left(z_i \right) = \N_{i}^{-1}\big( \Delta^{G_1} f_i \left( z_i \right) \big)^2 + \K_{G_i,z_i}\left( \N_i  \right)\Gamma^{G_i} \left(  f_i \right)\left(z_i \right) \quad z_1 := x \quad \text{and} \quad z_2 := p.
\end{align*}
\textbf{\small Claim:}	There is a sequence $f_{ij}:G_i \to \R$ and $\epsilon_{ij} \to 0$ with
\begin{align*}
\Gamma_2^{G_i} \left(  f_{ij} \right)\left(z_i \right) = \N_i^{-1}\big( \Delta^{G_1} f_{ij}\left(z_i \right) \big)^2 + \left( \K_{G_i,z_i}\left( \N_1  \right) + \epsilon_{ij} \right)\Gamma^{G_1} \left(  f_{ij} \right)\left(z_i \right) \quad \text{and} \quad \Gamma^{G_i}\left( f_{ij} \right)\left(z_i \right) \neq 0. 
\end{align*}
\textbf{\small Proof of claim:} If 	$\Gamma^{G_i} \left(  f_i \right)(z) \neq 0$, we set $f_{ij} = f_i$ for all $j$ and $\epsilon_j = 0$. If $\Gamma^{G_i} \left(  f_i \right)(z) = 0$, $f_i$ is locally constant at $z$. Set $f_{ij} = f_i + \frac{1}{j}\delta_z$. By the vertex-wise continuity of the curvature dimension inequalities, we can find such sequences $\epsilon_{ij}$. Obviously $\Gamma^{G_i}\left( f_{ij}  \right)\left(z_i \right) \neq 0$ since $f_{ij}$ and $f_i$ can not be locally constant at $z$ simultaneously. Notice that we must have $\epsilon_{ij} \ge 0$ since otherwise the above identity would contradict the definition of the curvature function. 
\scalebox{0.7}{\qed}
\par Without loss of generality, we consider four cases:
\begin{enumerate}
	\item [\textbf{(i)}] \textbf{Neither $x$ nor $p$ is strongly saturated.} 
	Pick curvature maximizers $f_i$ with $\Delta^{G_i} f_i(z_i) \neq 0$. Set $A_j := \alpha^{-2}(p)\Delta^{G_1}f_{1j}(x)$ and $B_j := \beta^{-2}(x)\Delta^{G_2}f_{2j}(p)$. For $j$ large enough, we can assume $A_j, B_j \neq 0$. By Lemma~\ref{lem:F-zeros}, $\mathcal{F}_j\left( c_{1j} , c_{2j}   \right) = 0$ ($\mathcal{F}_j$ is defined using $A_j$ and $B_j$) has a line of zeros. If both $x$ and $p$ are un-saturated, we can, by rescaling, further assume $\Delta^{G_1} f_1(x) = \alpha^2(p)\N_2^{-1}$ and $\Delta^{G_2}f_2(p) = \beta^2(x)\N_1^{-1}$ so the line of zeros satisfies $|c_1| = |c_2|$. Pick the zeros $\left(c_{1j} , c_{2j} \right)$ of $\mathcal{F}_j$ with $\left(|c_{1j}| , |c_{2j}| \right) \to \left( 1 , 1   \right)$ as $j \to \infty$. Then as $j \to \infty$,
	\begin{align*}
	c_{1j}^{-2}\mathcal{Q}_1\left( c_{1j} , c_{2j}   \right) \to \mathcal{Q}_1\left( 1 , 1   \right) =  \nicefrac{1}{2}\,  \beta^{-2} \Delta^{G_2}  \alpha^{-2} +  \beta^{-2} \Deg^{G_1} + \nicefrac{1}{2} \, \alpha^{-2}\Gamma^{G_1} \left( \beta^{-2}  \right),
	\end{align*}
	and
	\begin{align*}
	c_{2j}^{-2}\mathcal{Q}_2\left( c_{1j} , c_{2j}  \right) \to \mathcal{Q}_2\left( 1 , 1  \right)  =   \nicefrac{1}{2} \, \alpha^{-2} \Delta^{G_1}  \beta^{-2} +  \alpha^{-2} \Deg^{G_2} +  \nicefrac{1}{2}\,  \beta^{-2}\Gamma^{G_2} \left( \alpha^{-2}  \right).
	\end{align*}
	Set $ U^j_{c_{1j} , c_{2j} } = c_{1j}f_{1j} \oplus c_{2j}f_{2j}$,
	\begin{align*}
	\Gamma_2 \big(U^j_{c_{1j} , c_{2j} } \big) &= \; c_{1j}^2 \alpha^{-4} \Gamma_2^{G_1} \left( f_{1j} \right) + c_{2j}^2 \beta^{-4} \Gamma_2^{G_2} \left( f_{2j} \right) + \mathcal{Q}\left( c_{1j} , c_{2j} \right)\notag \\ &\le \; c_{1j}^2 \alpha^{-4} \N_1^{-1}\left( \Delta^{G_1}f_{1j} \right)^2 + c_{2j}^2 \beta^{-4} \N_2^{-1}\left( \Delta^{G_2}f_{2j} \right)^2 \notag \\ &+ \;\Big( c_{1j}^2 \alpha^{-4} \big(\K_{G_1,x}\left( \N_1  \right) + \epsilon_{1j} \big) +  \mathcal{Q}_1\left( c_{1j} , c_{2j} \right) \Big)  \Gamma^{G_1} \left(  f_{1j} \right) \notag \\ &+\; \Big( c_{1j}^2 \beta^{-4} \big( \K_{G_2,p}\left( \N_2  \right) + \epsilon_{2j} \big) +  \mathcal{Q}_2\left( c_{1j} , c_{2j} \right) \Big)   \Gamma^{G_2} \left(  f_{2j} \right) \notag \\ &= \left( \N_1 + \N_2  \right)^{-1}\Big( c_{1j} \alpha^{-2} \Delta^{G_1}f_{1j} + c_{2j} \beta^{-2} \Delta^{G_2}f_{2j} \Big)^{2}  \notag \\ &+ \Big(  \alpha^{-2} \left(\K_{G_1,x}\left( \N_1  \right) + \epsilon_{1j} \right) + \alpha^2 c_{1j}^{-2} \mathcal{Q}_1\left( c_{1j} , c_{2j} \right) \Big) \alpha^{-2}c_{1j}^2  \Gamma^{G_1} \left(  f_{1j} \right) \notag \\ &+ \Big( \beta^{-2} \left( \K_{G_2,p}\left( \N_2  \right) + \epsilon_{2j} \right) + \beta^2c_{2j}^{-2} \mathcal{Q}_2\left( c_{1j} , c_{2j} \right) \Big) \beta^{-2}  c_{2j}^2 \Gamma^{G_2} \left(  f_{2j} \right). \notag
	\end{align*}
	Thus,
	\begin{align*}
	\Gamma_2 \big(U^j_{c_{1j} , c_{2j} } \big) \le \big( \N_1 + \N_2  \big)^{-1}\big( \Delta \; U^j_{c_{1j} , c_{2j} } \big)^{2}  + \K_j \; \Gamma \big(U^j_{c_{1j} , c_{2j} } \big),
	\end{align*}
	where
	\begin{align}
	&\K_j =  \notag \\ &\Big(\alpha^{-2} \left(\K_{G_1,x}\left( \N_1  \right) + \epsilon_{1j} \right) + \alpha^2 c_{1j}^{-2} \mathcal{Q}_1\left( c_{1j} , c_{2j} \right)\Big) \vee \Big( \beta^{-2} \left( \K_{G_2,p}\left( \N_2  \right) + \epsilon_{2j} \right) + \beta^2c_{2j}^{-2} \mathcal{Q}_2\left( c_{1j} , c_{2j} \right) \Big). \notag
	\end{align}
	This implies $\K_{(x,p)}\left( \N_1 + \N_2  \right) \le \K_j$. Taking the limit as $j \to \infty$,
	\begin{align*}
	\K_{(x,p)}\left( \N_1 + \N_2  \right) \le \Big(   \alpha^{-2} \K_{G_1,x}\left( \N_1  \right) + \alpha^2 \mathcal{Q}_1\left( 1,1 \right)  \Big) \vee \Big(  \beta^{-2} \K_{G_2,p}\left( \N_2  \right)+ \beta^2\mathcal{Q}_2\left( 1,1 \right)  \Big).
	\end{align*}
	\item [\textbf{(ii)}]\textbf{$x$ is un-saturated and $p$ is weakly saturated or vice versa.} In this case $\Delta^{G_1}f_1 \neq 0$ and $\Delta^{G_2} f_2 = 0$. So we can assume $c_1=0$ and $c_2 = 1$ is a zero of $\mathcal{F}$. Setting $ U^j_{c_{2j}} = c_{2j}f_{2j}$ and using Lemma~\ref{lem:est},
	\relax
	\begin{align*}
	&\Gamma_2 \big(U^j_{c_{2j} } \big) \notag \\ &= c_{2j}^2 \beta^{-4} \Gamma_2^{G_2} \left( f_{2j} \right) + \mathcal{Q}\left( 0 , c_{2j} \right)\notag \\ &\le  c_{2j}^2 \beta^{-4} \N_2^{-1}\big( \Delta^{G_2}f_{2j} \big)^2+ \Big( c_{1j}^2 \beta^{-4} \big( \K_{G_2,p}\left( \N_2  \right) + \epsilon_{2j} \big) +  \mathcal{Q}_{2}\left( 0 , c_{2j} \right) \Big)   \Gamma^{G_2} \left(  f_{2j} \right) \notag \\ &= \left( \N_1 + \N_2  \right)^{-1}\big(  c_{2j} \beta^{-2} \Delta^{G_2}f_{2j} \big)^{2} + \Big( \N_2^{-1} -   \left( \N_1 + \N_2  \right)^{-1}  \Big) \big(  c_{2j} \beta^{-2} \Delta^{G_2}f_{2j} \big)^{2}\notag \\ &+ \Big( \beta^{-2} \left( \K_{G_2,p}\left( \N_2  \right) + \epsilon_{2j} \right) + \nicefrac{1}{2}\, \beta^2 \alpha^{-2} \Delta^{G_1} \beta^{-2} \Big) \beta^{-2} c_{2j}^2 \Gamma^{G_2} \left(  f_{2j} \right). \notag \\ &\le \left( \N_1 + \N_2  \right)^{-1}\left(  c_{2j} \beta^{-2} \Delta^{G_2}f_{2j} \right)^{2} + \K_j\; \beta^{-2} c_{2j}^2 \Gamma^{G_2} \left(  f_{2j} \right) \notag,
	\end{align*}
	where
	\begin{align}
	\K_j =  \beta^{-2} \big( \K_{G_2,p}\left( \N_2  \right) + \epsilon_{2j} \big) + \nicefrac{1}{2}\, \beta^2 \alpha^{-2} \Delta^{G_1} \beta^{-2} + 2\beta^{-2}\N_1\N_2^{-1}\left(  \N_1 + \N_2  \right)^{-1} \D_p. \notag
	\end{align}
	Taking the limit as $j \to \infty$, we deduce
	\begin{align}
	\K_{(x,p)}\left( \N_1 + \N_2 \right) &\le \K_{G_2,p}\left( \N_2  \right) + \nicefrac{1}{2}\, \beta^2 \alpha^{-2} \Delta^{G_1} \beta^{-2} + 2\beta^{-2}\N_1\N_2^{-1}\left(  \N_1 + \N_2  \right)^{-1} \D_p \notag \\ &= \K_{G_2,p}\left( \N_2  \right) + 2\beta^{-2}\N_1\N_2^{-1}\left(  \N_1 + \N_2  \right)^{-1} \D_p + \beta^2\mathcal{Q}_2\left(0,1\right)  \notag.
	\end{align}
	The proof of the other case follows similarly.. 
	\item [\textbf{(iii)}] \textbf{$x$ and $p$ are both weakly saturated.} In this case $\Delta^{G_1}f_1 = \Delta^{G_2} f_2 = 0$ so any $\left(c_1 , c_2 \right)$ solves $\mathcal{F}=0$ therefore, 
	\begin{align}
	\K_{(x,p)}\left( \N_1 + \N_2 \right) &\le \Big( \alpha^{-2} \K_{G_1,x}\left( \N_1  \right)+  \alpha^2  \mathcal{Q}_{14}\left( 1,0  \right) \Big)   \vee \Big( \beta^{-2}\K_{G_2,p}\left( \N_2  \right) + \mathcal{Q}_{24}\left( 0,1  \right)  \Big) \notag \\ &= \Big( \alpha^{-2}\K_{G_1,x}\left( \N_1  \right)+  \alpha^2 \mathcal{Q}_1\left( 1,0  \right)  \Big)  \vee \Big( \beta^{-2}\K_{G_2,p}\left( \N_2  \right) + \beta^2 \mathcal{Q}_2\left( 0,1  \right)  \Big). \notag
	\end{align}
	\item [\textbf{(iv)}] \textbf{$x$ and $p$ are both weakly saturated but neither is strongly saturated} This is a sub case of \textbf{(ii)}. Combining the bounds obtained in \textbf{(ii)}, we deduce
	\begin{align}
	\K_{(x,p)}\left( \N_1 + \N_2 \right) &\le \Big( \alpha^{-2} \K_{G_1,x}\left( \N_1  \right) + \alpha^2 \mathcal{Q}_1\left( 1,0 \right) + 2\alpha^{-2}\N_1^{-1}\N_2\left(  \N_1 + \N_2  \right)^{-1} \D_x \Big) \notag \\ &\wedge \Big( \beta^{-2} \K_{G_2,p}\left( \N_2  \right)+ \beta^2\mathcal{Q}_2\left( 0,1 \right) + 2\beta^{-2}\N_1\N_2^{-1}\left(  \N_1 + \N_2  \right)^{-1} \D_p \Big). \notag 
	\end{align}
\end{enumerate}\qed
\subsubsection*{\bf \textit{Proof of Theorem~\ref{thm:main-6}}}
Using $u_1 \otimes 1$ as a test function in Lemma~\ref{lem:sec-form-Gamma2},
\begin{align}\label{eq:test-func-Gamma2}
\Gamma_2 \left( u_1 \otimes 1  \right)(x,p) = \alpha^{-4}(p) \Gamma_2^{G_1} (u_1)(x) + \nicefrac{1}{2}\, \beta^{-2}(x)\Gamma^{G_1} \left( u_1 \right)(x) \Delta^{G_2} \alpha^{-2} (p), 
\end{align}
\begin{align*}
\Gamma \left(  u_1 \otimes 1  \right) (x,p) = \alpha^{-2}(p)\Gamma^{G_1} \left( u_1  \right)(x) \quad \text{and} \quad \Delta \left(  u_1 \otimes 1  \right) (x,p) = \alpha^{-2}(p)\Delta^{G_1}  u_1 (x).
\end{align*}
Hence, by the definition of $\K_{(x,p)}$ and using (\ref{eq:test-func-Gamma2}) we get
\begin{align*}
\Gamma_2^{G_1} (u_1)(x) \ge \N^{-1}\big( \Delta^{G_1} u_1(x) \big)^2 + \Big( \alpha^2(p)\mathcal{K}_{(x,p)}(\N) - \nicefrac{1}{2}\, \alpha^4(p)\beta^{-2}(x) \Delta^F  \alpha^{-2} (p)  \Big) \Gamma^{G_1} \left( u_1  \right)(x);
\end{align*}
which implies
\begin{align*}
\mathcal{K}_{G_1,x}(\N) \ge \alpha^2(p)\mathcal{K}_{(x,p)}(\N) - \nicefrac{1}{2}\, \alpha^4(p)\beta^{-2}(x) \Delta^{G_2}  \alpha^{-2} (p),
\end{align*}
or
\begin{align*}
\mathcal{K}_{(x,p)}(\N) \le \alpha^{-2}(p) \mathcal{K}_{G_1,x}(\N) + \nicefrac{1}{2}\, \alpha^2(p) \beta^{-2}(x)\Delta^{G_2}  \alpha^{-2} (p);
\end{align*}
similarly,
\begin{align*}
\mathcal{K}_{(x,p)}(\N) \le \beta^{-2}(x) \mathcal{K}_{G_2,x}(\N) + \nicefrac{1}{2}\, \beta^2(x)\alpha^{-2}(p)  \Delta^{G_1}  \beta^{-2} (x).
\end{align*}
\subsection{ Lower bounds for curvature function}
\subsubsection*{\bf \textit{Bounding the curvature form from below}}
As we will see shortly, for lower bounds, we will need to bound
\begin{align}
	\mathcal{Q}\left( 1, 1  \right) &=  \nicefrac{1}{2}\, \alpha^{-2}(p)   \Gamma^{G_2}\left(  f_2 \right)(p)  \Delta^{G_1}  \beta^{-2}(x)  -   \alpha^{-2}(p) \Delta^{G_2} f_2 (p)\Gamma^{G_1} \left(  \beta^{-2}  , f_1 \right)(x)  \notag \\ & + \; \nicefrac{1}{2}\,  \beta^{-2}(x)   \Gamma^{G_1}\left(  f_1 \right)(x)  \Delta^{G_2}  \alpha^{-2}(p)  -  \beta^{-2}(x) \Delta^{G_1} f_1 (x)\Gamma^{G_2} \left(  \alpha^{-2}  , f_2 \right)(p) \notag,
\end{align}
from below. 
\par Using estimates from Lemma~\ref{lem:est} in combination with with the Young's inequality
\begin{align*}
ab \ge - \epsilon^{-2}a^2 - \epsilon^2 b^2,
\end{align*}
we can find lower bounds for the terms involved in $	\mathcal{Q}\left( 1, 1  \right) $ in terms of gradients and Laplacians squared. Indeed,
\begin{align*}
	- \alpha^{-2} \Delta^{G_2} f_2 (p)\Gamma^{G_1} \left(  \beta^{-2}  , f_1 \right)(x) &\ge - \epsilon_2^2 \alpha^{-2} \left( \Delta^{G_2} f_2 (p)  \right)^2 - \epsilon_2^{-2} \alpha^{-2} \left( \Gamma^{G_1} \left(  \beta^{-2}  , f_1 \right)(x)  \right)^2 \\ &\ge - \epsilon_2^2 \alpha^{-2} \left( \Delta^{G_2} f_2 (p)  \right)^2 -\epsilon_2^{-2} \alpha^{-2} \Gamma^{G_1} \left(  \beta^{-2}\right)(x)  \Gamma^{G_1} \left(f_1 \right)(x),
\end{align*}
and
\begin{align*}
	-	\beta^{-2} \Delta^{G_1} f_1 (x)\Gamma^{G_2} \left(  \alpha^{-2}  , f_2 \right)(p) &\ge - \epsilon_1^2 \beta^{-2} \left( \Delta^{G_1} f_1 (x)  \right)^2 - \epsilon_1^{-2} \beta^{-2} \left( \Gamma^{G_2} \left(  \alpha^{-2}  , f_2 \right)(p)  \right)^2 \\ &\ge - \epsilon_1^2 \beta^{-2} \left( \Delta^{G_1} f_1 (x)  \right)^2 - \epsilon_1^{-2} \beta^{-2} \Gamma^{G_2} \left(  \alpha^{-2}\right)(p)  \Gamma^{G_2} \left(f_2 \right)(p),
\end{align*}
hold for $\epsilon_1$ and $\epsilon_2$ to be determined later. 
\par As a result,
\begin{align*}
		\mathcal{Q}\left( 1, 1  \right)  &\ge  \left(-\epsilon^2_1\alpha^{4}\beta^{-2}\right)  \alpha^{-4} \left( \Delta^{G_1} (f_1)\right)^2 + \left(-\epsilon^2_2\alpha^{-2} \beta^{4}\right) \beta^{-4} \left( \Delta^{G_2} (f_2)\right)^2 \\ & + \Big( -\epsilon_2^{-2}\Gamma^{G_1} \left(  \beta^{-2}\right) + \nicefrac{1}{2}\; \alpha^2\beta^{-2} \Delta^{G_2}\alpha^{-2} \Big) \alpha^{-2}\Gamma^{G_1}\left( f_1\right) \\ &+ \Big( -\epsilon_1^{-2} \Gamma^{G_2} \left(  \alpha^{-2}\right)+ \nicefrac{1}{2}\, \alpha^{-2}\beta^{2} \Delta^{G_1}\beta^{-2} \Big) \beta^{-2}\Gamma^{G_2}\left( f_2\right).
\end{align*}
Hence, 
\begin{align*}
	\Gamma_2 \left( f_1 \oplus f_2 \right)(x,p) &= \alpha^{-4} \Gamma_2^{G_1} \left( f_1 \right) + \beta^{-4} \Gamma_2^{G_2} \left( f_2 \right) + \mathcal{Q}\left( 1 , 1 \right)  \\ &\ge \alpha^{-4} \N_1^{-1}\left( \Delta^{G_1}f_1 \right)^2 +  \beta^{-4} \N_2^{-1}\left( \Delta^{G_2}f_2 \right)^2 + \mathcal{Q} \left( 1 , 1 \right)  \\ &+  \alpha^{-4} \K_{G_1,x}\left( \N_1  \right)    \Gamma^{G_1} \left(  f_1 \right) +  \beta^{-4} \K_{G_2,p}\left( \N_2  \right)   \Gamma^{G_2} \left(  f_2 \right) \\ &\ge \left( \N_1^{-1} - \epsilon^2_1\alpha^{4}\beta^{-2} \right) \alpha^{-4} \left( \Delta^{G_1} (f_1)\right)^2 \\ & + \left( \N_2^{-1} -\epsilon^2_2\alpha^{-2}\beta^{4}    \right) \beta^{-4} \left( \Delta^{G_2} (f_2)\right)^2 \\ &+ \left( \alpha^{-2} \K_{G_1,x}\left( \N_1  \right)   -\epsilon_2^{-2} \Gamma^{G_1} \left(  \beta^{-2}\right) + \nicefrac{1}{2}\, \alpha^{2}\beta^{-2} \Delta^{G_2}\alpha^{-2}  \right) \alpha^{-2}\Gamma^{G_1}\left( f_1\right) \\ & + \Big( \beta^{-2} \K_{G_2,p}\left( \N_2  \right)  -\epsilon_1^{-2}\Gamma^{G_2} \left(  \alpha^{-2}\right) + \nicefrac{1}{2}\, \alpha^{-2}\beta^{2} \Delta^{G_1}\beta^{-2}   \Big) \beta^{-2}\Gamma^{G_2}\left( f_2\right).
\end{align*}
\par Now for any $\M > \max \left\{  \N_1, \N_2 \right\}$, pick any pair $\epsilon_{1,\M}$ and $\epsilon_{2,\M}$ via the criteria
\begin{align*}
\N_1^{-1} - \epsilon^2_{1,\M}\alpha^{4}\beta^{-2} \ge \M^{-1} \quad \text{and} \quad  \N_2^{-1} - \epsilon^2_{2,\M}\alpha^{-2}\beta^{4}  \ge \M^{-1};
\end{align*}
i.e. where
\begin{align*}
\epsilon_{1,\M}^2 \le \alpha^2\beta^{-2}\left(\N_1^{-1} - \M^{-1}  \right)\quad  \text{and} \quad \epsilon_{2,\M}^2 \le \alpha^{-2} \beta^2\left(\N_2^{-1} - \M^{-1}  \right).
\end{align*}
For any such set of parameters $\epsilon_i$ (which depend on the synthetic dimensions and the values of the warping functions at the vertices $x \in G_1$ and $p \in G_2$), we get
\begin{align}\label{eq:lb-L-1}
	\Gamma_2 \left( f_1 \oplus f_2 \right)(x,p) &\ge \M^{-1} \alpha^{-4} \left( \Delta^{G_1}f_1 \right)^2 + \M^{-1} \beta^{-4} \left( \Delta^{G_2}f_2 \right)^2 \\ & + \Big( \alpha^{-2} \K_{G_1,x}\left( \N_1  \right)   -\epsilon_{2,\M}^{-2} \Gamma^{G_1} \left(  \beta^{-2}\right) + \nicefrac{1}{2} \, \alpha^{2}\beta^{-2} \Delta^{G_2}\alpha^{-2}  \Big) \alpha^{-2}\Gamma^{G_1}\left( f_1\right) \notag \\ & + \Big( \beta^{-2} \K_{G_2,p}\left( \N_2  \right)  -\epsilon_{1,\M}^{-2} \Gamma^{G_2} \left(  \alpha^{-2}\right) + \nicefrac{1}{2}\, \alpha^{-2}\beta^{2} \Delta^{G_1}\beta^{-2}   \Big) \beta^{-2}\Gamma^{G_2}\left( f_2\right) \notag.
\end{align}
\subsubsection*{\bf \textit{Proof of Theorem~\ref{thm:main-lb-1}}}
\par The Laplacian terms can be treated in two ways.
\par \noindent \textsf{First way:} \emph{Keeping the dimension fixed by decreasing (compromising on) the lower curvature bound:}
\begin{align}\label{eq:lb-L-2}
\M^{-1} \alpha^{-4} \left( \Delta^{G_1}f_1 \right)^2 + \M^{-1} \beta^{-4} \left( \Delta^{G_2}f_2 \right)^2 &= \M^{-1}\big( \alpha^{-2}  \Delta^{G_1}f_1 + \beta^{-2}  \Delta^{G_2}f_2  \big)^2 \notag \\ &- 2\M^{-1}\alpha^{-2}\beta^{-2}\Delta^{G_1}f_1  \Delta^{G_2}f_2 \\ &\ge \M^{-1}\left( \alpha^{-2}  \Delta^{G_1}f_1 + \beta^{-2}  \Delta^{G_2}f_2  \right)^2 \notag \\ &  - \M^{-1}\left( \alpha^{-2}  \Delta^{G_1}f_1   \right)^2 -\M^{-1} \left( \beta^{-2} \Delta^{G_2}f_2 \right)^2 \notag \\ & \ge  \M^{-1}\left( \alpha^{-2}  \Delta^{G_1}f_1 + \beta^{-2}  \Delta^{G_2}f_2  \right)^2 \notag \\ & - 2\M^{-1}\alpha^{-4}\Deg_i\Gamma^{G_1}\left( f_1 \right) - 2\M^{-1}\beta^{-4}\Deg_2\Gamma_2\left(f_2  \right)\notag,
\end{align}
in which, the first inequality follows from Young's inequality and the second inequality follows from the estimates in Lemma~\ref{lem:est}. Now \eqref{eq:lb-L-1} combined with \eqref{eq:lb-L-2} gives
\begin{align*}
	&\Gamma_2 \left( f_1 \oplus f_2 \right)(x,p) \ge  \\ &\M^{-1}\big( \alpha^{-2}  \Delta^{G_1}f_1 + \beta^{-2}  \Delta^{G_2}f_2  \big)^2 \\ & + \Big( \alpha^{-2} \K_{G_1,x}\left( \N_1  \right)   - \epsilon_{2,\M}^{-2} \Gamma^{G_1} \left(  \beta^{-2}\right) + \nicefrac{1}{2}\, \alpha^{2} \beta^{-2} \Delta^{G_2}\alpha^{-2} - 2\M^{-1}\alpha^{-2}\Deg_i \Big) \alpha^{-2}\Gamma^{G_1}\left( f_1\right) \\ & + \Big( \beta^{-2} \K_{G_2,p}\left( \N_2  \right)  -\epsilon_{1,\M}^{-2} \Gamma^{G_2} \left(  \alpha^{-2}\right) + \nicefrac{1}{2}\, \alpha^{-2}\beta^{2} \Delta^{G_1}\beta^{-2}  - 2\M^{-1}\beta^{-2}\Deg_2 \Big) \beta^{-2}\Gamma^{G_2}\left( f_2\right).
\end{align*}
\par In terms of curvature functions, we have thus shown
\begin{align*}
	&\K_{G, (x,p)} (\M) \ge \\
&  \Big(   \alpha^{-2} \K_{G_1,x}\left( \N_1  \right)   - \epsilon_{2,\M}^{-2} \Gamma^{G_1} \left(  \beta^{-2}\right) + \nicefrac{1}{2}\, \alpha^{2} \beta^{-2} \Delta^{G_2}\alpha^{-2} - 2\M^{-1}\alpha^{-2}\Deg_i  \Big)   \\ & \wedge \Big(  \beta^{-2} \K_{G_2,p}\left( \N_2  \right)  -\epsilon_{1,\M}^{-2} \Gamma^{G_2} \left(  \alpha^{-2}\right) + \nicefrac{1}{2}\, \alpha^{-2}\beta^{2} \Delta^{G_1}\beta^{-2}  - 2\M^{-1}\beta^{-2}\Deg_2    \Big),
\end{align*}
holds for any $\M > \max \left\{  \N_1, \N_2 \right\}$. 
\par \noindent  \textsf{Second way:} \emph{Obtaining a larger lower bound by increasing the synthetic dimension:}
\begin{align*}
	& \M^{-1} \alpha^{-4} \left( \Delta^{G_1}f_1 \right)^2 + \M^{-1} \beta^{-4} \left( \Delta^{G_2}f_2 \right)^2 \\ &= \left(2\M\right)^{-1}\big( \alpha^{-2}  \Delta^{G_1}f_1 + \beta^{-2}  \Delta^{G_2}f_2  \big)^2 \\ &+ \left(2\M\right)^{-1} \Big( \left(  \alpha^{-2}  \Delta^{G_1}f_1  \right)^2 + \left(  \beta^{-2}  \Delta^{G_2}f_2  \right)^2 - 2\left( \alpha^{-2}  \Delta^{G_1}f_1  \right)\left(  \beta^{-2}  \Delta^{G_2}f_2 \right)  \Big) \\ &\ge \left(2\M\right)^{-1} \big( \alpha^{-2}  \Delta^{G_1}f_1 + \beta^{-2}  \Delta^{G_2}f_2  \big)^2,
\end{align*}
which upon substitution in \eqref{eq:lb-L-2}, affords us
\begin{align*}
	\Gamma_2 \left( f_1 \oplus f_2 \right)(x,p) &\ge  (2\M)^{-1}\big( \alpha^{-2}  \Delta^{G_1}f_1 + \beta^{-2}  \Delta^{G_2}f_2  \big)^2 \\ & + \Big( \alpha^{-2} \K_{G_1,x}\left( \N_1  \right)   -\epsilon_{2,\M}^{-2} \Gamma^{G_1} \left(  \beta^{-2}\right)+ \nicefrac{1}{2}\, \alpha^{2} \beta^{-2} \Delta^{G_2}\alpha^{-2} \Big) \alpha^{-2}\Gamma^{G_1}\left( f_1\right) \\ & + \Big( \beta^{-2} \K_{G_2,p}\left( \N_2  \right)  - \epsilon_{1,\M}^{-2} \Gamma^{G_2} \left(  \alpha^{-2}\right) + \nicefrac{1}{2}\, \alpha^{-2}\beta^{2} \Delta^{G_1}\beta^{-2} \Big) \beta^{-2}\Gamma^{G_2}\left( f_2\right).
\end{align*}
\par Thus, we have proven
\begin{align*}
	&\K_{G, (x,p)} (2\M) \ge \\
	&  \Big(  \alpha^{-2} \K_{G_1,x}\left( \N_1  \right)   -\epsilon_{2,\M}^{-2} \Gamma^{G_1} \left(  \beta^{-2}\right)+ \nicefrac{1}{2}\, \alpha^{2} \beta^{-2} \Delta^{G_2}\alpha^{-2} \Big)   \\ & \wedge \Big(  \beta^{-2} \K_{G_2,p}\left( \N_2  \right)  - \epsilon_{1,\M}^{-2} \Gamma^{G_2} \left(  \alpha^{-2}\right) + \nicefrac{1}{2}\, \alpha^{-2}\beta^{2} \Delta^{G_1}\beta^{-2}    \Big)
\end{align*}
In the extremal case where 
\begin{align*}
\epsilon_{1,\M} = \alpha\beta^{-1}\big(\N_1^{-1} - \M^{-1}  \big)^{\nicefrac{1}{2}} \quad  \text{and} \quad \epsilon_{2,\M} = \alpha^{-1} \beta\big(\N_2^{-1} - \M^{-1}  \big)^{\nicefrac{1}{2}},
\end{align*}
and for any $\M > \max \left\{  \N_1, \N_2 \right\}$, we therefore get
\begin{enumerate}
\item \medskip
\begin{align*}
	&\K_{G, (x,p)} (\M) \ge \\
	&  \Big(   \alpha^{-2} \K_{G_1,x}\left( \N_1  \right)   - \alpha^{2} \beta^{-2}\big(\N_2^{-1} - \M^{-1}  \big)^{-1}  \Gamma^{G_1} \left(  \beta^{-2}\right) + \nicefrac{1}{2}\, \alpha^{2} \beta^{-2} \Delta^{G_2}\alpha^{-2} - 2\M^{-1}\alpha^{-2}\Deg_i  \Big)   \\ & \wedge \Big(  \beta^{-2} \K_{G_2,p}\left( \N_2  \right)  -\alpha^{-2} \beta^{2}\big(\N_1^{-1} - \M^{-1}  \big)^{-1}  \Gamma^{G_2} \left(  \alpha^{-2}\right) + \nicefrac{1}{2}\, \alpha^{-2}\beta^{2} \Delta^{G_1}\beta^{-2}  - 2\M^{-1}\beta^{-2}\Deg_2    \Big),
\end{align*}
\item[] and  \medskip \medskip
\item[(2)]
 \begin{align*}
	&\K_{G, (x,p)} (2\M) \ge \\
	&  \Big(  \alpha^{-2} \K_{G_1,x}\left( \N_1  \right)   -\alpha^{2} \beta^{-2}\big(\N_2^{-1} - \M^{-1}  \big)^{-1} \Gamma^{G_1} \left(  \beta^{-2}\right)+ \nicefrac{1}{2}\, \alpha^{2} \beta^{-2} \Delta^{G_2}\alpha^{-2} \Big)   \\ & \wedge \Big(  \beta^{-2} \K_{G_2,p}\left( \N_2  \right)  -\alpha^{-2} \beta^{2}\big(\N_1^{-1} - \M^{-1}  \big)^{-1}  \Gamma^{G_2} \left(  \alpha^{-2}\right) + \nicefrac{1}{2}\, \alpha^{-2}\beta^{2} \Delta^{G_1}\beta^{-2}    \Big).
\end{align*}
\end{enumerate} 
\qed
\subsubsection*{\bf \textit{Proof of Corollary~\ref{cor:iso-dim-lb}}}
For any function $u: G_1 \bp G_2 \to \R$, it holds
\begin{eqnarray}
	\Gamma_2(u) &\ge& \frac{\left(\Delta u \right)^2}{2\N} + \mathcal{K}_{(x,p)}(2\N) \Gamma(u) \notag \\ &=& \frac{\left(\Delta u \right)^2}{\N} - \frac{1}{2\N} \left(\Delta u \right)^2 + \mathcal{K}_{(x,p)}(2\N) \Gamma(u) \notag \\ &\ge& \frac{\left(\Delta u \right)^2}{\N} + \Big[ \mathcal{K}_{(x,p)}(2\N) - \N^{-1}\left( \alpha^{-2}\D_x + \beta^{-2}\D_p  \right)  \Big] \Gamma(u) \notag,
\end{eqnarray}
where in the last line we have used Lemma~\ref{lem:est} and the fact that
\begin{align*}
\Deg \left(  (x,p)  \right) = \alpha^{-2}\; \Deg_{G_1}(x) + \beta^{-2}\; \Deg_{G_2}(p). 
\end{align*}
Therefore,
\begin{equation}\label{eq:lb-same-d-2}
	\mathcal{K}_{(x,p)}(\N) \ge \mathcal{K}_{(x,p)}(2\N) - \N^{-1}\left( \alpha^{-2}\D_x + \beta^{-2}\D_p  \right). 
\end{equation}
Combining (\ref{eq:2n-lb}) and (\ref{eq:lb-same-d-2}) it follows
\begin{align*}
	\begin{split}
		\mathcal{K}_{(x,p)}(\N) &\ge 
		\Big(   \alpha^{-2} \K_{G_1,x}\left( \N  \right)   - \alpha^{2} \beta^{-2} \left( 2\N \right) \Gamma^{G_1} \left(  \beta^{-2}\right) + \nicefrac{1}{2}\, \alpha^{2} \beta^{-2} \Delta^{G_2}\alpha^{-2} - 2 \N^{-1}\alpha^{-2}\Deg_i  \Big) \\ &  \wedge  \Big(  \beta^{-2} \K_{G_2,p}\left( \N \right)  -\alpha^{-2} \beta^{2} \left( 2\N \right)  \Gamma^{G_2} \left(  \alpha^{-2}\right) + \nicefrac{1}{2}\, \alpha^{-2}\beta^{2} \Delta^{G_1}\beta^{-2}  - 2 \N^{-1}\beta^{-2}\Deg_2    \Big) \\ & - \N^{-1}\left( \alpha^{-2}\D_x + \beta^{-2}\D_p  \right).
	\end{split} 
\end{align*} \qed
\section{Examples}\label{sec:examples}
Working out examples by hand can become very hard, very fast. So we have considered simple examples where one does not need computational aid or very ad-hoc arguments; also we have used, for the product factors, the type of graphs for which the precise curvature functions are known. 

In this section, we will only invoke the all cases inclusive bounds 
\begin{align*}
	&	\Big(   \alpha^{-2} \K_{G_1,x}\left( \N_1  \right)   - \alpha^{2} \beta^{-2}\N_2  \Gamma^{G_1} \left(  \beta^{-2}\right) + \nicefrac{1}{2}\, \alpha^{2} \beta^{-2} \Delta^{G_2}\alpha^{-2} \Big) \\ &  \wedge  \Big(   \beta^{-2} \K_{G_2,p}\left( \N_2  \right)  -\alpha^{-2} \beta^{2}\N_1  \Gamma^{G_2} \left(  \alpha^{-2}\right) + \nicefrac{1}{2}\, \alpha^{-2}\beta^{2} \Delta^{G_1}\beta^{-2}     \Big) \\[10pt]
	& \le	\mathcal{K}_{(x,p)}(\infty) \le \\[10pt] & \Big( \alpha^{-2} \K_{G_1,x}\left( \infty \right) + \nicefrac{1}{2} \, \alpha^2 \beta^{-2} \Delta^{G_2} \alpha^{-2}  \Big) \wedge \Big(  \beta^{-2} \K_{G_2,p}\left( \infty \right)+ \nicefrac{1}{2}\, \beta^2 \alpha^{-2} \Delta^{G_1} \beta^{-2} \Big) \notag,
\end{align*}
on the curvature function which follows from Corollary~\ref{cor:non-iso-upb} and Theorem~\ref{thm:lb-infty}.
\begin{example}[Warped product of complete graphs with preferred vertices]\label{ex:1}
Consider the complete graphs $K_n$ and $K_m$. Then for the combinatorial Laplacian (i.e. when the vertex measure $m$ and the edge weights $\omega$ are identically equal to $1$) we know
\begin{align}\label{eq:curv-complete}
\K_{K_n,x} = \frac{n+2}{2} - \frac{2(n-1)}{\N}, \quad 0< \N \le \infty;
\end{align}
see~\cite{JL},~\cite{KKRT} for proofs and~\cite{LM} for a different proof when $\N = \infty$. 
\par Consider a preferred vertex $x_0$ in $G_1 := K_n$ and a preferred vertex $p_0$ in $G_2 := K_m$ by considering the warping functions 
\begin{align*}
\alpha(p) := \begin{cases}  \lambda & p = p_0 \\  1 & p \neq p_0  \end{cases} = \lambda\delta_{p_0} + \delta_{G_2\smallsetminus \{p_0\}}, \quad \text{and} \quad \beta(x) := \begin{cases}  \mu & x = x_0 \\ 1  & x \neq x_0  \end{cases} = \lambda\delta_{x_0} + \delta_{G_1\smallsetminus \{x_0\}} .
\end{align*}
\par The doubly warped product $K_n\, {_\alpha\diamond_\beta}\, K_m$ is the graph on $nm$ vertices $(x,p)$ where
\begin{align*}
(x,p) \sim (y,p)\quad \forall x\neq y \quad \text{and} \quad (x,p) \sim (x,q) \quad \forall p \neq q,
\end{align*}
the vertex measure is identically $1$ at all vertices and where
\begin{align*}
\omega_{(x,p),(y,p)} = \begin{cases} \lambda^{-2} & p = p_0 \\ 1 &p \neq p_0    \end{cases}, \quad \text{and}\quad \omega_{(x,p),(x,q)} = \begin{cases} \mu^{-2} & x = x_0 \\ 1 &x \neq x_0    \end{cases}.
\end{align*}
It is also straightforward from the definition of Laplacian and the gradient squared that
\begin{align*}
	\Delta \beta^{-2} (x):= \Delta \left( \mu^{-2}\delta_{x_0} + \delta_{G_1\smallsetminus \{x_0\}} \right) (x) = \begin{cases} (n-1)\left( 1 - \mu^{-2} \right)  & x = x_0 \\ \mu^{-2} - 1  & x \sim x_0 \end{cases},
\end{align*}
\begin{align*}
	\Delta \alpha^{-2} (p):= \Delta \left( \lambda^{-2}\delta_{p_0} + \delta_{G_1\smallsetminus \{p_0\}} \right) (p) = \begin{cases} (m-1)\left( 1 - \lambda^{-2} \right)  & x = p_0 \\ \lambda^{-2} - 1  & p \sim p_0 \end{cases},
\end{align*}
\begin{align*}
	\Gamma^{G_1} \beta^{-2} (x):= \Gamma^{G_1} \left( \mu^{-2}\delta_{x_0} + \delta_{G_1\smallsetminus \{x_0\}} \right) (x) = \begin{cases} \nicefrac{1}{2}\, (n-1)\left( 1 - \mu^{-2} \right)^2  & x = x_0 \\ \nicefrac{1}{2}\,  \left(\mu^{-2} - 1\right)^2  & x \sim x_0 \end{cases},
\end{align*}
\begin{align*}
	\Gamma^{G_2} \alpha^{-2} (p):= \Gamma^{G_2} \left( \lambda^{-2}\delta_{p_0} + \delta_{G_1\smallsetminus \{p_0\}} \right) (p) = \begin{cases} \nicefrac{1}{2}\, (m-1)\left( 1 - \lambda^{-2} \right)^2  & p = p_0 \\ \nicefrac{1}{2}\,  \left(\lambda^{-2} - 1\right)^2  & p \sim p_0 \end{cases}.
\end{align*}
\par There are four possible cases that we discuss in below.
\begin{enumerate}
		\item [\bf (1-1)] $x=x_0$ and $p=p_0$:
\begin{align}\label{eq:curv-1}
	&\Big(   \lambda^{-2} \left( \frac{n+2}{2} - \frac{2(n-1)}{\N_1} \right)  - \lambda^2 \mu^{-2}\N_2  \nicefrac{1}{2}\, (n-1)\left( 1 - \mu^{-2} \right)^2 + \nicefrac{1}{2}\, \lambda^{2} \mu^{-2} (m-1)\left( 1 - \lambda^{-2} \right)  \Big) \notag  \\ & \wedge  \Big(   \mu^{-2} \left( \frac{m+2}{2} - \frac{2(m-1)}{\N_2} \right)    -\lambda^{-2} \mu^{2}\N_1 \nicefrac{1}{2}\, (m-1)\left( 1 - \lambda^{-2} \right)^2 + \nicefrac{1}{2}\, \lambda^{-2}\mu^{2} (n-1)\left( 1 - \mu^{-2} \right)   \Big)\notag  \\[10pt] 	& \le 	\mathcal{K}_{(x_0,p_0)}(\infty)  \le \\[10pt]  &  \Big(  \nicefrac{1}{2}\,  \lambda^{-2} (n+2) + \nicefrac{1}{2}\, \lambda^{2} \mu^{-2} (m-1)\left( 1 - \lambda^{-2} \right)   \Big)\notag  \\ &\wedge \Big(  \nicefrac{1}{2}\,    \mu^{-2}\left( m+2 \right) + \nicefrac{1}{2}\, \lambda^{-2}\mu^{2} (n-1)\left( 1 - \mu^{-2} \right)  \Big) \notag;
\end{align}\medskip
	\item  [\bf (1-2)] $x = x_0$ and $p \sim p_0$:
\begin{align*}
	&\Big(  \left( \frac{n+2}{2} - \frac{2(n-1)}{\N_1} \right)  - \mu^{-2}\N_2  \nicefrac{1}{2}\, (n-1)\left( 1 - \mu^{-2} \right)^2 + \nicefrac{1}{2}\, \mu^{-2}\left( \lambda^{-2} - 1\right)  \Big) \\ & \wedge  \Big(   \mu^{-2} \left( \frac{m+2}{2} - \frac{2(m-1)}{\N_2} \right)    -\mu^{2}\N_1 \nicefrac{1}{2}\,\left( 1 - \lambda^{-2} \right)^2 + \nicefrac{1}{2}\, \mu^{2} (n-1)\left( 1 - \mu^{-2} \right)   \Big) \\[10pt] 	& \le 	\mathcal{K}_{(x_0,p)}(\infty)  \le \\[10pt]  &  \Big(  \nicefrac{1}{2}\, (n+2) + \nicefrac{1}{2}\, \mu^{-2}\left( 1 - \lambda^{-2} \right)   \Big) \wedge \Big(  \nicefrac{1}{2}\,    \mu^{-2}\left( m+2 \right) + \nicefrac{1}{2}\, \mu^{2} (n-1)\left( 1 - \mu^{-2} \right)  \Big) \notag;
\end{align*}
	\item  [\bf (1-3)] $x \sim x_0$ and $p = p_0$:
	\begin{align*}
		&\Big(  \left( \frac{m+2}{2} - \frac{2(m-1)}{\N_2} \right)  - \lambda^{-2}\N_1  \nicefrac{1}{2}\, (m-1)\left( 1 - \lambda^{-2} \right)^2 + \nicefrac{1}{2}\, \lambda^{-2}\left( \mu^{-2} - 1\right)  \Big) \\ & \wedge  \Big(   \lambda^{-2} \left( \frac{n+2}{2} - \frac{2(n-1)}{\N_1} \right)    -\lambda^{2}\N_2 \nicefrac{1}{2}\,\left( 1 - \mu^{-2} \right)^2 + \nicefrac{1}{2}\, \lambda^{2} (m-1)\left( 1 - \lambda^{-2} \right)   \Big) \\[10pt] 	& \le 	\mathcal{K}_{(x,p_0)}(\infty)  \le \\[10pt]  &  \Big(  \nicefrac{1}{2}\, (m+2) + \nicefrac{1}{2}\, \lambda^{-2}\left( 1 - \mu^{-2} \right)   \Big) \wedge \Big(  \nicefrac{1}{2}\,    \lambda^{-2}\left( n+2 \right) + \nicefrac{1}{2}\, \lambda^{2} (m-1)\left( 1 - \lambda^{-2} \right)  \Big) \notag; 
	\end{align*}
	\item[\bf (1-4)] $x \sim x_0$ and $p \sim p_0$:
\begin{align*}
	&\Big( \left( \frac{n+2}{2} - \frac{2(n-1)}{\N_1} \right)  - \N_2  \nicefrac{1}{2}\, \left( 1 - \mu^{-2} \right)^2 + \nicefrac{1}{2}\, \left( 1 - \lambda^{-2} \right)  \Big) \\ & \wedge  \Big(  \left( \frac{m+2}{2} - \frac{2(m-1)}{\N_2} \right)    - \N_1 \nicefrac{1}{2}\, \left( 1 - \lambda^{-2} \right)^2 + \nicefrac{1}{2}\, \left( 1 - \mu^{-2} \right)   \Big) \\[10pt] 	& \le 	\mathcal{K}_{(x,p)}(\infty)  \le \\[10pt]  &  \Big(  \nicefrac{1}{2}\,  (n+2) + \nicefrac{1}{2}\, \left( 1 - \lambda^{-2} \right)   \Big) \wedge \Big(  \nicefrac{1}{2}\,    \left( m+2 \right) + \nicefrac{1}{2}\, \left( 1 - \mu^{-2} \right)  \Big) \notag.
\end{align*}
\end{enumerate}
\par These bounds indicate that as we move towards the vertex $\left( x_0, p_0 \right)$, the curvature decreases in the sense that the lower and upper bounds increase. This should not be very surprising since the edge weights increase as we move towards $\left( x_0, p_0 \right)$ and the curvature is positively correlated with edge weights. 
 \end{example}
\begin{example}[Warped product of a star graph and a complete graph with preferred vertices]\label{ex:2}
	We know in the star graph $S_m$ with the head denoted by $p_0$ and tentacles by $p$, one has
\begin{align*}
\K_{S_m, p_0}\left( \infty \right) = \frac{3-m}{2},\quad  \K_{S_m, p}\left( \infty \right) = \frac{5 - m}{2};
\end{align*}
and
\begin{align*}
	\K_{S_m, p_0}\left( \N\right) = \begin{cases} \frac{3+ m}{2} - \frac{2m}{\N} & 0 < \N \le 2 \\ \frac{3-m}{2} & \N > 2. \end{cases} \quad  \text{and} \quad \K_{S_m, p}\left( \N \right) = \frac{5-m}{2} - \frac{2m}{\N} \quad \forall \N \in (0,\infty];
\end{align*}
see~\cite{JL} for the proof. 
\par Set $G_1 = K_n$ and $G_2 = S_m$ and take the warping functions
\begin{align*}
\alpha(p) := \lambda \delta_{p_0} + \delta_{S_m \smallsetminus \{p_0\}} \quad \text{and}\quad \beta(x) := \mu \delta_{x_0} + \delta_{K_n \smallsetminus \{x_0\}}, \quad \lambda, \mu \ge 1.
\end{align*}
\par We have
\begin{align*}
	\Delta \alpha^{-2} (p):= \Delta \left( \lambda^{-2}\delta_{p_0} + \delta_{G_1\smallsetminus \{p_0\}} \right) (p) = \begin{cases} m \left( 1 - \lambda^{-2} \right)  & x = p_0 \\ \lambda^{-2} - 1  & p \sim p_0  \end{cases},
\end{align*}
and
\begin{align*}
	\Gamma^{G_2} \alpha^{-2} (p):= \Gamma^{G_2} \left( \lambda^{-2}\delta_{p_0} + \delta_{G_1\smallsetminus \{p_0\}} \right) (p) = \begin{cases} \nicefrac{1}{2}\, m \left( 1 - \lambda^{-2} \right)^2  & p = p_0 \\ \nicefrac{1}{2}\,  \left(\lambda^{-2} - 1\right)^2  & p \sim p_0  \end{cases}.
\end{align*}
\par Therefore, if for simplicity, we assume the synthetic dimension $\N_2 > 2$, we get the following bounds hold on the curvature function in the doubly warped product $K_{n} \, {_\alpha\diamond_\beta}\,  S_{m}$. There are again four cases in total. 
\begin{enumerate}
		\item [\bf (2-1)] $x=x_0$ and $p=p_0$:
	\begin{align*}
		&\Big(   \lambda^{-2} \left( \frac{n+2}{2} - \frac{2(n-1)}{\N_1} \right)  - \lambda^2 \mu^{-2}\N_2  \nicefrac{1}{2}\, (n-1)\left( 1 - \mu^{-2} \right)^2 + \nicefrac{1}{2}\, \lambda^{2} \mu^{-2} m \left( 1 - \lambda^{-2} \right)  \Big) \\ & \wedge  \Big(   \mu^{-2} \left( \frac{3-m}{2} \right)    -\lambda^{-2} \mu^{2}\N_1 \nicefrac{1}{2}\, m\left( 1 - \lambda^{-2} \right)^2 + \nicefrac{1}{2}\, \lambda^{-2}\mu^{2} (n-1)\left( 1 - \mu^{-2} \right)   \Big) \\[10pt] 	& \le 	\mathcal{K}_{(x_0,p_0)}(\infty)  \le \\[10pt]  &  \Big(  \nicefrac{1}{2}\,  \lambda^{-2} (n+2) + \nicefrac{1}{2}\, \lambda^{2} \mu^{-2} m\left( 1 - \lambda^{-2} \right)   \Big) \wedge \Big(  \nicefrac{1}{2}\,    \mu^{-2}\left( 3-m \right) + \nicefrac{1}{2}\, \lambda^{-2}\mu^{2} (n-1)\left( 1 - \mu^{-2} \right)  \Big) \notag; 
	\end{align*}
	\medskip
	\item  [\bf (2-2)] $x = x_0$ and $p \sim p_0$:
	\begin{align*}
		&\Big(  \left( \frac{n+2}{2} \right)  - \mu^{-2}\N_2  \nicefrac{1}{2}\, (n-1)\left( 1 - \mu^{-2} \right)^2 + \nicefrac{1}{2}\, \mu^{-2}\left( \lambda^{-2} - 1\right)  \Big) \\ & \wedge  \Big(   \mu^{-2} \left( \frac{5-m}{2} - \frac{2m}{\N_2} \right)    -\mu^{2}\N_1 \nicefrac{1}{2}\,\left( 1 - \lambda^{-2} \right)^2 + \nicefrac{1}{2}\, \mu^{2} (n-1)\left( 1 - \mu^{-2} \right)   \Big) \\[10pt] 	& \le 	\mathcal{K}_{(x_0,p)}(\infty)  \le \\[10pt]  &  \Big(  \nicefrac{1}{2}\, (n+2) + \nicefrac{1}{2}\, \mu^{-2}\left( 1 - \lambda^{-2} \right)   \Big) \wedge \Big(  \nicefrac{1}{2}\,    \mu^{-2}\left( 5-m \right) + \nicefrac{1}{2}\, \mu^{2} (n-1)\left( 1 - \mu^{-2} \right)  \Big) \notag;
	\end{align*}
		\item  [\bf (2-3)] $x \sim x_0$ and $p = p_0$:
	\begin{align*}
		&\Big(  \left( \frac{3-m}{2} \right)  - \lambda^{-2}\N_1  \nicefrac{1}{2}\, m\left( 1 - \lambda^{-2} \right)^2 + \nicefrac{1}{2}\, \lambda^{-2}\left( \mu^{-2} - 1\right)  \Big) \\ & \wedge  \Big(   \lambda^{-2} \left( \frac{n+2}{2} - \frac{2(n-1)}{\N_1} \right)    -\lambda^{2}\N_2 \nicefrac{1}{2}\,\left( 1 - \mu^{-2} \right)^2 + \nicefrac{1}{2}\, \lambda^{2} m\left( 1 - \lambda^{-2} \right)   \Big) \\[10pt] 	& \le 	\mathcal{K}_{(x,p_0)}(\infty)  \le \\[10pt]  &  \Big(  \nicefrac{1}{2}\, (3-m) + \nicefrac{1}{2}\, \lambda^{-2}\left( 1 - \mu^{-2} \right)   \Big) \wedge \Big(  \nicefrac{1}{2}\,    \lambda^{-2}\left( n+2 \right) + \nicefrac{1}{2}\, \lambda^{2} m\left( 1 - \lambda^{-2} \right)  \Big) \notag;
	\end{align*}
		\item[\bf (2-4)] $x \sim x_0$ and $p \sim p_0$:
		\begin{align*}
		&\Big( \left( \frac{n+2}{2} - \frac{2(n-1)}{\N_1} \right)  - \N_2  \nicefrac{1}{2}\, \left( 1 - \mu^{-2} \right)^2 + \nicefrac{1}{2}\, \left( 1 - \lambda^{-2} \right)  \Big) \\ & \wedge  \Big(  \left( \frac{5-m}{2} - \frac{2m}{\N_2} \right)    - \N_1 \nicefrac{1}{2}\, \left( 1 - \lambda^{-2} \right)^2 + \nicefrac{1}{2}\, \left( 1 - \mu^{-2} \right)   \Big) \\[10pt] 	& \le 	\mathcal{K}_{(x,p)}(\infty)  \le \\[10pt]  &  \Big(  \nicefrac{1}{2}\,  (n+2) + \nicefrac{1}{2}\, \left( 1 - \lambda^{-2} \right)   \Big) \wedge \Big(  \nicefrac{1}{2}\,    \left( 5-m\right) + \nicefrac{1}{2}\, \left( 1 - \mu^{-2} \right)  \Big) \notag. 
	\end{align*}
\end{enumerate}
\end{example}
\begin{example}[Warped product of a complete bipartite graph and a complete graph with preferred vertices]
	Based on~\cite{CLP}, for a complete bi-partite graph $\K_{n,m}$ and a vertex $p$ with $\deg (p) = n$, the Bakry-\'Emery curvature function satisfies
\begin{align*}
\K_{K_{n,m}, p} (\infty) = \begin{cases} \frac{m+2 - \left| n - 2m + 2   \right|}{2} & (n,m) \neq (1,1) \\ 2 & o.t. \end{cases},
\end{align*}
and
\begin{align*}
	\K_{K_{n,m}, p} (\N) = \begin{cases} \frac{4 + n - m}{2} - \frac{2n}{\N} & \text{$n=1$ or $n \le 2m-2$}\\ \frac{4 + n - m}{2} - \frac{2n}{\N} &  \text{$n> \max\{1, 2m-2\}$ and  $0 < \N \le \frac{2n}{n-2m+2}$} \\  \frac{3m-n}{2} & \text{$n> \max\{1, 2m-2\}$ and  $ \N > \frac{2n}{n-2m+2}$} \end{cases}.
\end{align*}
Set $G_1 = K_r$ and $G_2 = K_{n,m} $ and the warping functions 
\begin{align*}
\alpha(p) =  \lambda\delta_{p_0} +\delta_{K_{n,m} \smallsetminus \{p_0\}}\quad  (\deg(p_0) = n), \quad \text{and}\quad \beta(x) = \mu \delta_{x_0} + \delta_{K_r \smallsetminus \{x_0\}}.
\end{align*}
where
\begin{align*}
	\Delta \alpha^{-2} (p):= \Delta \left( \lambda^{-2}\delta_{p_0} + \delta_{G_2\smallsetminus \{p_0\}} \right) (p) = \begin{cases} n \left( 1 - \lambda^{-2} \right)  & x = p_0 \\ \lambda^{-2} - 1  & p \sim p_0\;\; (\deg(p) = m)\\ 0 & o.t. \;\;( p\neq p_0, \deg(p) = n) \end{cases},
\end{align*}
and
\begin{align*}
	\Gamma^{G_2} \alpha^{-2} (p):= \Gamma^{G_2} \left( \lambda^{-2}\delta_{p_0} + \delta_{G_2\smallsetminus \{p_0\}} \right) (p) = \begin{cases} \nicefrac{1}{2}\, n \left( 1 - \lambda^{-2} \right)^2  & p = p_0 \\ \nicefrac{1}{2}\,  \left(\lambda^{-2} - 1\right)^2  & p \sim p_0\;\; (i.e. \; \deg(p) = m) \\  0 & o.t. \;\;( p\neq p_0, \deg(p) = n)\end{cases}.
\end{align*}
\par For simplicity and in order to reduce the number of possible cases, we assume $n$ and $m$ and $\N_2$ are chosen so that we have
\begin{align*}
\K_{K_{n,m}, p} (\N_2) =  \frac{4 + n - m}{2} - \frac{2n}{\N} \quad \text{and} \quad \K_{K_{n,m}, p} (\infty) =  \frac{m+2 - \left| n - 2m + 2   \right|}{2};
\end{align*}
the general case is not harder but it involves more cases. For example one way to guarantee the above assumption, is to let $\N_2$ be arbitrary and pick $n$ and $m$ so that $n \le 2m-2 \le 4n - 6$ holds.
\par Then, in the warped product graph $K_{r}\, {_\alpha\diamond_\beta}\, K_{n,m} $,  we have the following six cases. 
\begin{enumerate}
		\item [\bf (3-1)] $x=x_0$ and $p=p_0$:
	\begin{align*}
		&\Big(   \lambda^{-2} \left( \frac{r+2}{2} - \frac{2(r-1)}{\N_1} \right)  - \lambda^2 \mu^{-2}\N_2  \nicefrac{1}{2}\, (r-1)\left( 1 - \mu^{-2} \right)^2 + \nicefrac{1}{2}\, \lambda^{2} \mu^{-2} n \left( 1 - \lambda^{-2} \right)  \Big) \\ & \wedge  \Big(   \mu^{-2} \left(  \frac{4 + n - m}{2} - \frac{2n}{\N_2}  \right)    -\lambda^{-2} \mu^{2}\N_1 \nicefrac{1}{2}\, n \left( 1 - \lambda^{-2} \right)^2  + \nicefrac{1}{2}\, \lambda^{-2}\mu^{2} (r-1)\left( 1 - \mu^{-2} \right)   \Big) \\[10pt] 	& \le 	\mathcal{K}_{(x_0,p_0)}(\infty)  \le \\[10pt]  &  \Big(  \nicefrac{1}{2}\,  \lambda^{-2} (r+2) + \nicefrac{1}{2}\, \lambda^{2} \mu^{-2} n \left( 1 - \lambda^{-2} \right)  \Big) \\ &\wedge \Big(  \nicefrac{1}{2}\,    \mu^{-2}\left( m+2 - \left| n - 2m + 2   \right| \right) + \nicefrac{1}{2}\, \lambda^{-2}\mu^{2} (r-1)\left( 1 - \mu^{-2} \right)  \Big) \notag;
	\end{align*}
	\medskip
	\item  [\bf (3-2)] $x = x_0$ and $p \sim p_0$:
	\begin{align*}
		&\Big(  \left( \frac{r+2}{2} - \frac{2(r-1)}{\N_1} \right)  - \mu^{-2}\N_2  \nicefrac{1}{2}\, (r-1)\left( 1 - \mu^{-2} \right)^2 + \nicefrac{1}{2}\, \mu^{-2}\left( \lambda^{-2} - 1\right)  \Big) \\ & \wedge  \Big(   \mu^{-2} \left(  \frac{4 + n - m}{2} - \frac{2n}{\N_2} \right)    -\mu^{2}\N_1 \nicefrac{1}{2}\,\left( 1 - \lambda^{-2} \right)^2 + \nicefrac{1}{2}\, \mu^{2} (r-1)\left( 1 - \mu^{-2} \right)   \Big) \\[10pt] 	& \le 	\mathcal{K}_{(x_0,p)}(\infty)  \le \\[10pt]  &  \Big(  \nicefrac{1}{2}\, (r+2) + \nicefrac{1}{2}\, \mu^{-2}\left( 1 - \lambda^{-2} \right)   \Big) \wedge \Big(  \nicefrac{1}{2}\,    \mu^{-2}\left( m+2 - \left| n - 2m + 2   \right| \right) + \nicefrac{1}{2}\, \mu^{2} (r-1)\left( 1 - \mu^{-2} \right)  \Big) \notag;
	\end{align*}
	\item[\bf (3-3)] $x = x_0$, $p \neq p_0$ and $p \not \sim p_0$:
		\begin{align*}
				&	\Big(  \left( \frac{r+2}{2} - \frac{2(r-1)}{\N_1} \right)  -  \mu^{-2}\N_2  \nicefrac{1}{2}\, (r-1)\left( 1 - \mu^{-2} \right)^2 \Big) \\ & \wedge  \Big(   \mu^{-2} \left(  \frac{4 + n - m}{2} - \frac{2n}{\N_2} \right)   + \nicefrac{1}{2}\, \mu^{2} (r-1)\left( 1 - \mu^{-2} \right)   \Big) \\[10pt] 
		&	\le \mathcal{K}_{(x_0,p)}(\infty) \le\\[10pt] & \Big(  \nicefrac{1}{2}\, (r+2) \Big) \wedge \Big(  \nicefrac{1}{2}\,    \mu^{-2}\left( m+2 - \left| n - 2m + 2   \right|  \right)+ \nicefrac{1}{2}\,  \mu^{2} (r-1)\left( 1 - \mu^{-2} \right) \Big) \notag;
		\end{align*}
	\item  [\bf (3-4)] $x \sim x_0$ and $p = p_0$:
	\begin{align*}
		&\Big(  \mu^{-2} \left(  \frac{4 + n - m}{2} - \frac{2n}{\N_2} \right)  - \lambda^{-2}\N_1  \nicefrac{1}{2}\, n\left( 1 - \lambda^{-2} \right)^2 + \nicefrac{1}{2}\, \lambda^{-2}\left( \mu^{-2} - 1\right)  \Big) \\ & \wedge  \Big(   \lambda^{-2} \left( \frac{r+2}{2} - \frac{2(r-1)}{\N_1} \right)    -\lambda^{2}\N_2 \nicefrac{1}{2}\,\left( 1 - \mu^{-2} \right)^2 + \nicefrac{1}{2}\, \lambda^{2} n \left( 1 - \lambda^{-2} \right)   \Big) \\[10pt] 	& \le 	\mathcal{K}_{(x,p_0)}(\infty)  \le \\[10pt]  &  \Big(  \nicefrac{1}{2}\, \left( m+2 - \left| n - 2m + 2   \right|  \right) + \nicefrac{1}{2}\, \lambda^{-2}\left( 1 - \mu^{-2} \right)   \Big) \wedge \Big(  \nicefrac{1}{2}\,    \lambda^{-2}\left( r+2 \right) + \nicefrac{1}{2}\, \lambda^{2}n \left( 1 - \lambda^{-2} \right)  \Big) \notag;
	\end{align*}
	\item[\bf (3-5)] $x \sim x_0$, $p \neq p_0$, and $p \not \sim p_0$: 
		\begin{align*}
		&	\Big(   \left( \frac{r+2}{2} - \frac{2(r-1)}{\N_1} \right)    - \N_2 \nicefrac{1}{2}\,  \left(\mu^{-2} - 1\right)^2 \Big) \\ &  \wedge  \Big(  \left(  \frac{4 + n - m}{2} - \frac{2n}{\N_2}  \right)  + \nicefrac{1}{2}\,  \left( \mu^{-2} - 1 \right)     \Big) \\[10pt]
		& \le	\mathcal{K}_{(x,p)}(\infty) \le \\[10pt] & \Big(  \frac{r+2}{2}  \Big) \wedge \Big(  \nicefrac{1}{2}\left( m+2 - \left| n - 2m + 2   \right|  \right) + \nicefrac{1}{2}\, \left( \mu^{-2} - 1 \right) \Big) \notag;
	\end{align*}
	\item[\bf (3-6)] $x \sim x_0$ and $p \sim p_0$:
		\begin{align*}
		&\Big( \left( \frac{r+2}{2} - \frac{2(r-1)}{\N_1} \right)  - \N_2  \nicefrac{1}{2}\, \left( 1 - \mu^{-2} \right)^2 + \nicefrac{1}{2}\, \left( 1 - \lambda^{-2} \right)  \Big) \\ & \wedge  \Big(  \left(  \frac{4 + n - m}{2} - \frac{2n}{\N_2}   \right)    - \N_1 \nicefrac{1}{2}\, \left( 1 - \lambda^{-2} \right)^2 + \nicefrac{1}{2}\, \left( 1 - \mu^{-2} \right)   \Big) \\[10pt] 	& \le 	\mathcal{K}_{(x,p)}(\infty)  \le \\[10pt]  &  \Big(  \nicefrac{1}{2}\,  (r+2) + \nicefrac{1}{2}\, \left( 1 - \lambda^{-2} \right)   \Big) \wedge \Big(  \nicefrac{1}{2}\,    \left(  m+2 - \left| n - 2m + 2   \right|  \right) + \nicefrac{1}{2}\, \left( 1 - \mu^{-2} \right)  \Big) \notag.
	\end{align*}	
\end{enumerate}
\end{example}
\section{Toy models}\label{sec:toymodels}
\par In this closing section, we present toy models which are constructed out of rudimentary graphs and are based on simple real life principles. Of course, these toy models cannot replace working with real data but we believe they can serve to illuminate the main ideas and the philosophy behind our definition of warped product. 
\subsection{ Robustness of networks in relation with curvature}
Examples of complex networks range from micro-biochemical and genetic networks such as cancer transcription network associated to a cancer tissue to world-wide spread infrastructures such as the Internet and social media networks. The theory has found its way to medicine where the effectiveness of drugs on cancer networks can be measured by graph theoretic quantities and to literature and history where the network representing the use of literary and artistic forms in early American writing can be used to study how people used to think in the old days!
\par In almost any network, much like any infrastructural network (such as Internet and/or social networks), the resilience of network in face of attack, malfunction and perturbations is a critical measure of how effective the network is. This resilience is what is generally referred to as the robustness of the network. In other words, robustness of a complex network is referred to the ability of the complex network to retain its functionality under failures or perturbations in the network. Among the various ways to interpret and measure the robustness of a network, percolation theory and fluctuation theorem are the two most mainstreamed ones and both are linked to various type of discrete curvatures. 
\subsubsection{\small \bf \textit{Robustness}}
The robustness $\mathcal{R}$ can be given a precise definition in various ways. The more dominant definition of robustness, in layman's terms, is \emph{the (logarithmic) rate at which fluctuations (from the steady network) in the network converge to their steady state.} This is a very relatable concept from theory of differential equations however now applied to perhaps non-deterministic discrete processes. 
\par The premise for this definition, is that there is a process (Markov, diffusion or stochastic) at work and we are looking at a neighborhood of an attractor; this attractor is what is refereed to as the steady state. So any small enough perturbation from this steady state, has to converge to the steady state as $t \to \infty$. The rate of convergence in probability then can be measured by the aid of Varadhan's large deviation theory; see~\cite{Varad}. Notice that the same concept can be applied to a process obtained from a dynamical system (governed by an ode) instead of a diffusion process. 
\par To be more precise, consider a random bounded perturbation in the network, let $P_{\epsilon}(t)$ denote the probability at which the mean of an energy functional (this is a potential function on the phase space which is set beforehand) deviates from the unperturbed mean at time $t$ by an amount more than $\epsilon$ then the robustness is given by
\begin{align*}
\mathcal{R} := - \lim_{t \to \infty} \left( \frac{1}{t} \ln P_{\epsilon}(t)  \right).
\end{align*}
\par The energy functional acts like a potential function when we have a dynamical system and in the case of a diffusion process, the process is in a sense the gradient flow of the energy functional. Existence of such limit is guaranteed by the theory of large deviations and the number can be found using statistical methods by selecting discrete samples and using the theory of large deviations as is done in~\cite{D}; also see~\cite{Young, Ellis, DGO}.
\par We need to clarify that the process that is used in the definition of robustness is not necessarily a deterministic one and it could as well be a stochastic process and the theory of large deviations still gives a logarithmic rate of convergence. A nice example of computing this rate of convergence for the Ornstein-Uhlenbeck process (and its relation to the Olivier's curvature) can be found in \cite{Ol, Oliv2}. 
\par In an ideal world, this rate of convergence at least for diffusion processes, would be the same as the Gaussian bounds for the heat kernel which again are directly related to the lower bounds on the Ricci curvature (at least in the smooth setting); see e.g. \cite{SV}. 
\par Another possible definition of robustness of a network is, roughly speaking, \emph{the ability of the network to hold it together in the face of vertex or edge malfunction}. This can be quantified via percolation theory. Given an unweighted random network, suppose the nodes could fail to operate with probability $p$. The percolation theory then provides a critical probability $p_c$ at which, the network breaks into many small fragments (clusters). Also there are universal constants that, in combination with this critical probability, will describe many other features of the network; see~\cite{AJB} for the origin of this approach. 
\par Measuring robustness as a number may not be an easy or a fruitful task (in terms of interpreting the obtained number) however it turns out that one can resort to changes in curvature to speculate the change in robustness which is a more useful and desired way of using robustness. Generally speaking, we want our good networks to be more robust and bad networks (like cancer tissues) to be less robust. This can be correlated with the changes in various types of curvature.  
\subsubsection{\small \bf \textit{Positive correlation among robustness, entropy and graph curvatures}}
In~\cite{DGO}, a fluctuation theorem for time dependent evolutionary populatory networks (think of a time dependent weighted graph) has been established. The fluctuation theorem asserts that: ``\emph{The rate of decay of fluctuations to their steady state is positively correlated with Boltzmann's entropy the given evolutionary network''}. The robustness in this context is interpreted as the ability of the evolutionary network to adapt to changes. This phenomenon is expressed as 
\begin{align*}
\Delta \mathcal{R} \Delta\mathrm{Ent} > 0. 
\end{align*}
where $\Delta \mathrm{Ent}$ denotes the change of Boltzmann's entropy of the steady state. The Boltzmann's entropy can be written as 
\begin{align*}
	\mathrm{Ent}\left( \mu \vert m \right) := - \int_X \; \mu \ln \frac{\mathrm{d}\mu}{\mathrm{d}m};
\end{align*}
once we decide what the measures involved are and what integration means. For example in the discrete setting, the integral is a weighted sum. 
\par Recall that for the definition of the robustness, one needs an attractor of a process. See~\cite{DGO, D, DM} for further details on the fluctuation theorem and the definition of entropy. 
\par An intuitive way to think about the fluctuation theorem is that the longest a perturbation takes to converge to the steady state, the more time the process would have to diffuse which means the steady state must have had higher entropy hence, the positive correlation between the robustness and the entropy of the steady state. 
\par On the other hand, in the realm of metric measure spaces equipped with a measure $m$ and a distance $d$, one way to define a lower bound $k$ on weak Ricci curvature is the $(k)$-convexity of the negative of the Boltzmann's entropy
\begin{align*}
\mathrm{Ent}\left( \mu \vert m \right) := - \int_X \; \mu \ln \frac{\mathrm{d}\mu}{\mathrm{d}m},
\end{align*}
along $L^2$-Wasserstein geodesics i.e. loosely speaking as $ \mathrm{Ent}'' \le - \Ric$. These go under the Lott-Sturm-Villani curvature-dimension bounds. To get a dimensional version one requires $(\K, \N)$-convexity of the negative of the Boltzmann entropy; e.g. see~\cite{Villani} for a detailed account. 
\par What Lott-Sturm-Villani bounds and the fluctuation theorem have in common is that both use rates of changes of the Boltzmann entropy (however of different orders). This suggests a correlation between the robustness of a network and the discrete version of Lott-Sturm-Villani curvature-dimension bounds. This correlation in the graph setting has not been rigorously mathematically proven yet all the numerous experimental work thus far is in favor of a positive correlation. We will mention some such experimental works in below.
\par In~\cite{Tetal}, the authors have explored this correlation by using real data based on cancer tissues. They have shown some important networks which are constructed based on cancer tissues posses higher curvature than normal tissues and higher entropy than normal tissues. On the other hand, by the previous works in quantitative biology, one knows the cancer tissues have higher robustness than normal tissues; see~\cite{WBS}. This indicates a positive correlation between curvature and robustness. 
\par In~\cite{SGRetal}, again using real life data, it is evidenced that the networks resulting from different kinds of cancer have different (quantifiably so) curvatures i.e. one can differentiate cancer tissues by their curvature. In the said references, the theory is tested using Lott-Sturm-Villani bounds and using scalar curvature bounds (an average of Olivier's edge curvatures) each intended to measure robustness of specific type of properties of the said networks. 
\par In quantitative finance, data experiments show that the market fragility (a quantity which is negatively correlated with the robustness) is negatively correlated with the Ricci curvature of various kinds; see~\cite{SGT, SGTetal, Netgeom}. 
\par So, at least experimentally, one expects that the positive correlation
\begin{align*}
\Delta \mathcal{R} \Delta \Ric > 0,
\end{align*}
to hold i.e. one expects the changes in robustness and changes in the Ricci curvature (of either of previously mentioned kinds) to be positively correlated as shown via experiments. 
\par What we are concerned with in these notes is the discrete Bakry-\'Emery curvature for weighted graphs. These bounds being more of an algebraic nature, are easier to find as the answer to a semi-definite programming problem; see~\cite[Section 2]{CLP}. This makes the discrete Bakry-\'Emery bounds more apt (here meaning programmable) for network and data science purposes. 
\par We point out that in the Riemannian setting, the Bakry-\'Emery curvature-dimension bounds coincide with the Lott-Sturm-Villani bounds (defined via $(\K, \N)$-convexity of the negative Boltzmann entropy) thanks to a well-behaved heat flow, the Bochner identity and the existence of a chain rule. Of course no such coincidence holds for graphs. The Bakry-\'Emery lower bounds (meaning the best such bound i.e. the curvature functions), the Olivier lower bounds and the lower bounds obtained via $\K$-convexity of the negative of Boltzmann entropy (along with a zoo of other notions of vertex and edge graph curvature), hardly ever coincide for a generic weighted graph; one, for the reason that on graphs, the chain rule for derivatives fails which will imply that the process generated by the discrete Laplacian is not a diffusion process in the standard sense hence the heat flow is not as well-behaved. 
\par Even-though various notions of graph curvature do not coincide, intuitionally speaking, we expect them to be positively correlated. In these notes, we will neither prove any result of this sort nor will we rigorously argue that the discrete Bakry-\'Emery Ricci bound is positively correlated with Robustness. We just note that, from the experiment point of view, there are again data experiments that confirm the positive correlation between the robustness of biological networks and the discrete Bakry-\'Emery curvature-dimension bounds. See~\cite{PRRetal, PMT} for the use of Bakry-\'Emery Ricci curvature bounds in determining the robustness.
\par So for the rest of this section, we work under the anzats 
\begin{anzats}\label{anz:robustness}
	$\Delta \K_{G}\left(\infty \right) \Delta \mathcal{R} > 0$ i.e. Bakey-\'Emery curvature functions and robustness are positively correlated. 
\end{anzats}
The speculations we make on robustness of our toy models only make sense under the above anzats. 
\subsection{ Toy Models: Scenarios modeled on doubly warped products}
Here we explain how doubly warped product of graphs can be used in modeling the interplay between two networks.  We will see how under Anzats~\ref{anz:robustness}, the curvature estimates can be used to speculate the robustness of such interplay networks.
\par The toy models presented here elucidate the general principle behind how doubly warped product can be used to model the interplay of two networks. One of the networks we consider is interpreted as the environment or the network of locations and the other one is another general network that repeats itself in different locations. This is not a far-fetched construction. Indeed, many large scale networks are structured this way and one expects a product structure to occur. We will dub the name $G_1$-$G_2$ interplay networks for the kind of networks we consider in the rest of these notes. 
\subsubsection{\small \bf \textit{Toy Model I: Profession-Community Model (the $P$-$C$ interplay network)}}\label{sec:toy-1}
Let us see how one can look at the interwoven structure of professions in a society and the effect they have on their local communities and vice versa as a doubly warped product of time dependent networks. Of course the principle behind this model can be applied to almost any other setting. 
\par \textsf{\small \bf {The reference $P$-network:}} This is a network of professions as nodes and the interactions amongst them as edge weights. In reality, the $P$-network is an evolutionary network as the interactions between professions vary in time. The edge weights represent the amount of interactions. For example if we consider the accepted norms about profession-related interactions, in a community, the physics professors and the judges are expected to have way less interactions than the policemen and the judges. But all of them are expected to have more or less the same amount of interactions with the barbers in the community. Since one would expect that all professions in a community would have some interactions, we can take the $P$-network to be a complete weighted graph. The parameters in this network can be for example determined by averaging the interactions among various professions in various locations (cities, countries, etc.); see Figure~\ref{fig:CP}.A. 
\par \textsf{\small \bf The reference $C$-network:} The $C$-network serves as a reference network of communities as nodes and their interactions as edge weights based on chosen measurements, norms and standards at the time. This again can be assumed to be a complete weighted network of some given size. The edge weights for example could represent the mutual political influence or interaction they have or the amount of migration in both directions in a \emph{reference ideal society} or to be an average of interactions of all sorts between various locations. Letting the edge weights to depend on time, will provide a more precise evolutionary model that will factor in how cold or warm the interactions of two locations are at different times. One can also add in vertex measures in order to make a more versatile model; see Figure~\ref{fig:CP}.B. 
\par In order to make things simpler to compute, we will set the vertex measures in both models to be $1$. 
\begin{figure}
	\subfloat[Reference P-graph: Graph of four professions with their interactions as edge weights.]
	{
		\begin{tikzpicture}[
			vertex/.style={draw,circle,minimum size=4pt,
				inner sep=2pt},
			weight/.style={font=\scriptsize},]

			\draw[blue , ultra thick] (0,0) node  [vertex, fill=blue, label=above: $P_1:$ Profesor] {};
			\draw[blue , ultra thick] (-2,-4) node  [vertex, fill=blue,label=left:$P_2:$ Priest] {};	
			\draw[blue , ultra thick] (2,-4) node  [vertex, fill=blue,label=right:$P_3:$ Police] {};
			\draw[blue , ultra thick] (0,-8) node  [vertex, fill=blue,label=below:$P_4:$ Postman] {};
			
			\draw[blue, ultra thick] (0,0) --  (-2,-4) [] {};
			\draw[blue, ultra thick, densely dotted] (-2,-4) --  (2,-4) [] {};
			\draw[blue, ultra thick] (2,-4) --  (0,0) [] {};
			\draw[blue, ultra thick] (-2,-4) --  (0,-8) [] {};
			\draw[blue, ultra thick] (0,-8) --  (2,-4) [] {};
			\draw[blue, ultra thick] (0,0) --  (0,-8) [] {};

			
			\draw[black] (-0.8,-1.8) node  [ label=left:$\omega^P_{12}$] {};
			\draw[black] (-0.3,-3.7) node  [ label=left:$\omega^P_{23}$] {};
			\draw[black] (0.8,-1.8) node  [ label=right:$\omega^P_{13}$] {};
			\draw[black] (-0.1,-5.1) node  [ label=right:$\omega^P_{24}$] {};
			\draw[black] (0.6,-6.6) node  [ label=right:$\omega^P_{24}$] {};
			\draw[black] (-0.6,-6.6) node  [ label=left:$\omega^P_{24}$] {};
		\end{tikzpicture}
	}
	\subfloat[Reference C-graph: Graph of four cities (communities) in the upstate New York with their interactions as edge weights. ]
	{
		\begin{tikzpicture}[
			vertex/.style={draw,circle,minimum size=4pt,
				inner sep=2pt},
			weight/.style={font=\scriptsize},]

			\draw[red , ultra thick] (5,0) node  [vertex, fill=red, label=above: $C_2:$ Carmel] {};
			\draw[red , ultra thick] (3,-3) node  [vertex, fill=red,label=right:$C_1:$ Clay] {};	
			\draw[red , ultra thick] (1,1) node  [vertex, fill=red,label=above:$C_3:$ Cairo] {};
			\draw[red , ultra thick] (0,-5) node  [vertex, fill=red,label=below:$C_4:$ Cayuta] {};

			
			\draw[red, ultra thick] (1,1) --  (3,-3) [] {};
			\draw[red, ultra thick] (3,-3) --  (0,-5) [] {};
			\draw[red, ultra thick] (3,-3) --  (5,0) [] {};
				\draw[red, ultra thick] (1,1) --  (3,-3) [] {};
			\draw[densely dotted, red, ultra thick] (5,0) --  (0,-5) [] {};
				\draw[red, ultra thick] (5,0) --  (1,1) [] {};
				\draw[red, ultra thick] (1,1) --  (0,-5) [] {};

			
			\draw[black] (4,-1.4) node  [ label=right:$\omega^C_{12}$] {};
			\draw[black] (3,-1) node  [ label=left:$\omega^C_{13}$] {};
			\draw[black] (1.5,-5) node  [ label=above:$\omega^C_{14}$] {};
			\draw[black] (1.5,-3.4) node  [ label=above:$\omega^C_{24}$] {};
				\draw[black] (3.5,1) node  [ label=left:$\omega^C_{23}$] {};
				\draw[black] (0.7,-1) node  [ label=left:$\omega^C_{34}$] {};

		\end{tikzpicture}
	}
	\caption{P and C networks as weighted graphs.}
	\label{fig:CP}
\end{figure}
\par \textsf{\small \bf The $P$-$C$ interplay network:}
Now we wish to model the global network of professions and their interactions as a doubly warped product $P\, {_\alpha\diamond_\beta} \, C$. To do so, the warping functions are to be chosen in a way that would signify the interplay between the two networks. For a profession $C_i$, $\beta(C_i)$ is to satisfy (be given by)
\begin{align*}
\beta^{-2}(C_i)\;\omega^P_{km} = \textit{interaction between the professions $P_k$ and $P_m$ within the community $C_i$}. 
\end{align*}
so $\beta$ is a measure of how much the interactions between professions in the city $C_i$ deviates from the reference ideal city norms. Similarly we pick the value $\alpha(P_i)$ via
\begin{align*}
	\alpha^{-2}(P_i)\;\omega^C_{km} &=  \textit{ interaction between the members of profession $P_i$ }   \\
	&\phantom{sajj} \textit{between communities $C_k$ and $C_m$}. 
\end{align*}
In the evolutionary version of this model, all these quantities will depend on time. The careful reader might have already noticed that in this model and the way we have picked our functions $\alpha$ and $\beta$, the quantity $\alpha(P_i)$ has to a priori also depend on the edges of the $C$-graph and $\beta(C_j)$ on the edges of the $P$-graph in order to obtain a more accurate model. This in fact is true and it defines what we called a doubly twisted product; see Definition~\ref{def:graph-dw}. We wish to study these type of products in a future but in order to be able to apply the results of the present paper, we will need our functions to only depend on the vertices of the underlying graphs namely, we are working with doubly warped products. 
\begin{figure}{
		\begin{center}
			\begin{tikzpicture}[
				vertex/.style={draw,circle,minimum size=4pt,
					inner sep=2pt},
				weight/.style={font=\scriptsize},
				]
				
				
				\draw[blue , ultra thick] (1,0) node  [vertex, fill=blue, label=above: $\left( P_1 \text{,}\,  C_2 \right)$] {};
				\draw[blue , ultra thick] (0.5,-1) node  [vertex, fill=blue,label=left:] {};	
				\draw[blue , ultra thick] (1.5,-1) node  [vertex, fill=blue,label=right:$\left( P_3 \text{,}\,  C_2 \right)$] {};
				\draw[blue , ultra thick] (1,-2) node  [vertex, fill=blue,label=below:$\left( P_4 \text{,}\,  C_2 \right)$] {};

				\draw[blue , ultra thick] (-3,-3) node  [vertex, fill=blue, label=above: $\left( P_1 \text{,}\,  C_1 \right)$] {};
				\draw[blue , ultra thick] (-3.5,-4) node  [vertex, fill=blue,label=left:$\left( P_2 \text{,}\,  C_1 \right)$] {};	
				\draw[blue , ultra thick] (-2.5,-4) node  [vertex, fill=blue,label=right:] {};
				\draw[blue , ultra thick] (-3,-5) node  [vertex, fill=blue,label=below:$\left( P_4 \text{,}\,  C_1 \right)$] {};

				
				\draw[blue, ultra thick] (1,0) --  (0.5,-1) [] {};
				\draw[blue, ultra thick] (0.5,-1) --  (1.5,-1) [] {};
				\draw[blue, ultra thick] (1.5,-1)  --  (1,0) [] {};
				\draw[blue, ultra thick] (1.5,-1) --  (1,0) [] {};
				\draw[blue, ultra thick] (0.5,-1) --  (1,-2) [] {};
				\draw[blue, ultra thick](1,-2) --  (1.5,-1) [] {};
				\draw[blue, ultra thick] (1,0) --  (1,-2) [] {};

				\draw[densely dotted,blue, ultra thick] (-3.5,-4) --  (-2.5,-4) [] {};
				\draw[densely dotted,blue, ultra thick] (- 3,-3) --  (-2.5,-4) [] {};
				\draw[blue, ultra thick] (-3.5,-4)  --  (-3,-3) [] {};
				\draw[blue, ultra thick] (-3.5,-4) --  (-3,-5) [] {};
				\draw[densely dotted,blue, ultra thick](-3,-5) --  (-2.5,-4) [] {};
				\draw[densely dotted,blue, ultra thick] (-3,-3) --  (-3,-5) [] {};


				\draw[red, ultra thick] (1,0) --  (-3,-3) [] {};
				\draw[red, ultra thick] (0.5,-1) --  (-3.5,-4) [] {};
				\draw[densely dotted, red, ultra thick] (1.5,-1) --   (-2.5,-4) [] {};
				\draw[red, ultra thick] (1,-2) --  (-3,-5) [] {};

				
				\draw[black] (-1,-3.5) node  [ label=right:$ \beta^{-2}(P_4) \omega^C_{12}$] {};
				\draw[black] (-3.1,-3.4) node  [ label=left:$\alpha^{-2}(C_1) \omega^P_{12}$] {};
				\draw[black] (-3.4,-1.4) node  [ label=right:$ \beta^{-2}(P_1) \omega^C_{12}$] {};
				\draw[black] (1.1,-0.3) node  [ label=right:$\alpha^{-2}(C_2) \omega^P_{13}$] {};
			\end{tikzpicture}
	\end{center}}
	\caption{P-C interplay network as a doubly warped product. What is depicted is the interaction between the fibers over the vertices $C_1$ and $C_2$.}
	\label{fig:CP-2}
\end{figure}
\par \textsf{\small \bf Robustness speculations:}\\
Consider the $P$-$C$ interplay network and suppose there is a preferred city (a capital of some sort) which we will assume to be $C_1$ and there is a preferred profession (sort of a pillar for the society), here of course in order to satisfy our fragile ego, we take the profession ``Professorship'' (vertex $P_1$) to be the preferred one. First notice that this is the exact setting of Example~\ref{ex:1} so we already have the needed curvature bounds. Therefore, we can apply the curvature estimates (1-1)-(1-4), by taking $p_0:= P_1$ and $x_0 := C_1$. As we observed in Example~\ref{ex:1}, the closer we get to the vertex $\left(x_0, p_0 \right)$ the lower the curvature becomes (in the sense that upper and lower bounds decrease as we get closer to the preferred vertex). 
\par There are various ways to make the curvature bounds obtained in Example~\ref{ex:1} to increase or decrease and to interpret that change in terms of robustness. Below, we will mention one such observation. 
\par \emph{Observation: ``To have a more robust profession-community network, the mismatch between interactions in the capital city and other places should not be too big; also no one profession should stand out (in terms of interaction with others) to the extreme degree.''}
\par Indeed, suppose the number of cities and professions are fixed. If $\mu^{-2}$ and $\lambda^2$ become large (with the same proportion) then from \eqref{eq:curv-1}, we deduce the upper and lower bounds on $\mathcal{K}_{(x_0,p_0)}(\infty)$ (and also curvature at other vertices) will decrease. This suggests in this case the robustness of the network lessens. 
\par The curvature bounds we have obtained in Example~\ref{ex:1} also enable us to compare two such interplay networks (say from two different countries) in terms of their robustness by just plugging in the parameters in the curvature bounds. 
\subsubsection{\small \bf \textit{Toy Model II: Simple Franchise Model (the $C$-$H$ interplay network)}}
This model is quite in the same spirit of the previous one, however the graph structures used are different. As before there is a $C$-network of communities, cities, etc. (in general geographical locations) modeled with a complete graph and there is a simple hierarchical network, $H$, consisting of a chief and some deputies or a main store and its satellite local branches, etc. 
\par \textsf{\small \bf The $C$-$H$ interplay network:}\\
The simple hierarchical $H$-network is modeled by a star graph. The weights in both graphs are determined as before by averaging the interactions over a large data sample in order to obtain the reference networks; see Figures~\ref{fig:HC} and \ref{fig:HC-2}.
\par Once the reference weights in each network are fixed, one can determine the warping functions $\alpha$ and $\beta$ using the same principle as before; namely, one chooses the values for these warping functions so that
\begin{align*}
\beta^{-2}(C_i)\;\omega^H_{km} = \text{interaction between the branches $H_k$ and $H_m$ within the location $C_i$}. 
\end{align*}
and
\begin{align*}
\alpha^{-2}(H_i)\;\omega^C_{km} = & \; \text{interaction between the branch $H_i$ in location $C_k$} \\  & \;\; \text{and the branch $H_i$ in the location  $C_m$}. 
\end{align*}
More complex franchises can be modeled using a more complex hierarchical network.
\begin{figure}
	{
		\begin{tikzpicture}[
			vertex/.style={draw,circle,minimum size=4pt,
				inner sep=2pt},
			weight/.style={font=\scriptsize}]

			\draw[ ultra thick] (5,0) node  [vertex, fill=black, label=above: $H_2:$ local branch] {};
			\draw[ultra thick] (3,-3) node  [vertex, fill=black,label=right:$H_1:$ local headquarters] {};	
			\draw[ ultra thick] (1,1) node  [vertex, fill=black,label=above:$H_3:$ local branch] {};
			\draw[ ultra thick] (0,-5) node  [vertex, fill=black,label=below:$H_4:$ local branch] {};

			
			\draw[ ultra thick] (1,1) --  (3,-3) [] {};
			\draw[ultra thick] (3,-3) --  (0,-5) [] {};
			\draw[ ultra thick] (3,-3) --  (5,0) [] {};

			
			\draw[black] (4,-1.4) node  [ label=right:$\omega^H_{12}$] {};
			\draw[black] (2.2,-1.3) node  [ label=left:$\omega^H_{23}$] {};
			\draw[black] (1.2,-4.1) node  [ label=above:$\omega^H_{13}$] {};
			
		\end{tikzpicture}
	}
	\caption{The hierarchical network $H$ which depicts a simple chain of commands.}
	\label{fig:HC}
\end{figure}
\begin{figure}
	{
		\begin{tikzpicture}[
			vertex/.style={draw,circle,minimum size=4pt,
				inner sep=2pt},
			weight/.style={font=\scriptsize},]

			\draw[ ultra thick] (0,0) node  [vertex, fill=black, label=left:  \scalebox{0.8}{$C_{_1}$}] {};
			\draw[ ultra thick] (-2,-4) node  [vertex, fill=black,label=below: \scalebox{0.8}{$C_{_2}$}] {};	
			\draw[ultra thick] (2,-4) node  [vertex, fill=black,label=below:  \scalebox{0.8}{$C_{_3}$}] {};
			\draw[ultra thick] (0,-8) node  [vertex, fill=black,label=left:  \scalebox{0.8}{$C_{_4}$}] {};
			
			\draw[ultra thick] (0,0) --  (-2,-4) [] {};
			\draw[ultra thick, densely dotted] (-2,-4) --  (2,-4) [] {};
			\draw[ultra thick] (2,-4) --  (0,0) [] {};
			\draw[ultra thick] (-2,-4) --  (0,-8) [] {};
			\draw[ultra thick] (0,-8) --  (2,-4) [] {};
			\draw[ultra thick] (0,0) --  (0,-8) [] {};

			 \begin{scope}[ shift={(0,0)}, scale = 1]
				\draw [ultra thick] (0,0) -- (0,1);
				\draw [ultra thick] (0,0) -- (0.7,1);
			\draw [ultra thick] (0,0) -- (-0.7,1);
			\draw[ ultra thick] (0,1) node  [vertex, fill=black] {};
			\draw[ ultra thick] (-0.7,1) node  [vertex, fill=black] {};	
			\draw[ultra thick] (0.7,1) node  [vertex, fill=black] {};
				
			\end{scope}

			 \begin{scope}[ shift={(-2,-4)}, scale = 1, rotate=90]
				\draw [ultra thick] (0,0) -- (0,1);
				\draw [ultra thick] (0,0) -- (0.7,1);
				\draw [ultra thick] (0,0) -- (-0.7,1);
				\draw[ ultra thick] (0,1) node  [vertex, fill=black] {};
				\draw[ ultra thick] (-0.7,1) node  [vertex, fill=black] {};	
				\draw[ultra thick] (0.7,1) node  [vertex, fill=black] {};
				
			\end{scope}

			\begin{scope}[ shift={(2,-4)}, scale = 1, rotate=-90]
				\draw [ultra thick] (0,0) -- (0,1);
				\draw [ultra thick] (0,0) -- (0.7,1);
				\draw [ultra thick] (0,0) -- (-0.7,1);
				\draw[ ultra thick] (0,1) node  [vertex, fill=black] {};
				\draw[ ultra thick] (-0.7,1) node  [vertex, fill=black] {};	
				\draw[ultra thick] (0.7,1) node  [vertex, fill=black] {};
				
			\end{scope}

			\begin{scope}[ shift={(0,-8)}, scale = 1, rotate=180]
			\draw [ultra thick] (0,0) -- (0,1);
			\draw [ultra thick] (0,0) -- (0.7,1);
			\draw [ultra thick] (0,0) -- (-0.7,1);
			\draw[ ultra thick] (0,1) node  [vertex, fill=black] {};
			\draw[ ultra thick] (-0.7,1) node  [vertex, fill=black] {};	
			\draw[ultra thick] (0.7,1) node  [vertex, fill=black] {};
			
		\end{scope}

		\draw [very thin] (0,1) -- (3,-4);
		\draw [very thin] (0.7,1) -- (3,-4.7);
		\draw [very thin] (-0.7,1) -- (3,-3.3);

			\begin{scope}[ xscale=-1]
				
				\draw [very thin] (0,1) -- (3,-4);
				\draw [very thin] (0.7,1) -- (3,-4.7);
				\draw [very thin] (-0.7,1) -- (3,-3.3);
				
		   \end{scope}

		\begin{scope}[ yscale=-1, shift={(0,8)}]
		
		\draw [very thin] (0,1) -- (3,-4);
		\draw [very thin] (0.7,1) -- (3,-4.7);
		\draw [very thin] (-0.7,1) -- (3,-3.3);
		
	\end{scope}

	\begin{scope}[ yscale=-1, xscale=-1, shift={(0,8)}]
	
	\draw [very thin] (0,1) -- (3,-4);
	\draw [very thin] (0.7,1) -- (3,-4.7);
	\draw [very thin] (-0.7,1) -- (3,-3.3);
	
\end{scope}
			
		\end{tikzpicture}
	}
	\caption{ The product $C {_\alpha\diamond_\beta} H$ as the whole franchise.}
	\label{fig:HC-2}
\end{figure}
\par \textsf{\small \bf Robustness speculations:}\\
For this model, we can apply the curvature bounds obtained in Example~\ref{ex:2}. So again under the assumption that robustness is positively correlated with the Bakry-\'Emery curvature function (Anzats~\ref{anz:robustness}), we point out the following observations.
\begin{enumerate}
\item Robustness decreases as the number of deputies (or the number of local branches under the local headquarter; the number $m$) increases. However 
\item For small number of deputies ($m \le 3$), the robustness is highest near the headquarters at the capital (the vertex $\left(p_0,x_0\right)$) however for larger number of deputies in the hierarchy ($m \ge 5$), the robustness is the the lowest at the main headquarters. 
\end{enumerate}
\subsubsection{\small \bf \textit{Toy Model III: Model of rival campaigns}}
\phantom{sajjad}
\par \textsf{\small \bf Description of the model:}\\
Suppose there exist various groups which are campaigning for an election and are in competition with each other. Assume that these campaigns have to use the same advertisement resources such as national TV channels, newspapers, etc. In addition, if we assume these advertisement outlets are also in competition with each other (namely, there is no official interaction among them as a confidentiality agreement would require the employees to act), we can locally (in one area) model this network by a complete bi-partite graph. Here, the graph is bi-partite since no two campaigns interact and no two media interact but all the campaigns interact with all the media. 
\par As before, multiplying this bi-partite graph by a complete graph of cities (campaigns and media have branches in all cities), we will get the whole campaigning network. 
\par \textsf{\small \bf Robustness speculations:}\\
With a reasoning similar to what we did in \S\hspace{1.2pt}\ref{sec:toy-1} we can observe that in such an election model, ``\emph{if there are one or two campaigns that are singled out from the others to the extreme degree or if there exist monarchical media (which outshine other media), then the robustness of the election network is less than when there is an even political or media play field.}'' 
\subsubsection{\small \bf \textit{Toy Model IV: Modeling quarantine in a pandemic}}
\phantom{sajjad}
\par \textsf{\small \bf Description of the model:}
\par A simple way of modeling quarantine (here, by quarantine, we mean imposing interaction policies) in a pandemic such as the Covid-19 pandemic is to divide the environment into various geographical blocks and keep track of the number of virus positive, virus negative, vaccine administered, at risk, etc. agents in each geographical block. This way, the whole system can be modeled on a doubly warped product of $K_m$ where $m$ is the number of categories that we put people into and $K_n$ for $n$ a large number which signifies the number of geographical blocks. The warping functions in $\K_m {_\alpha\diamond_\beta} K_n$ are again determined by the amount of interactions between the people in various locations and the interaction of people in different categories in one location. 
\begin{figure}{
		\begin{center}
			\begin{tikzpicture}[
				vertex/.style={draw,circle,minimum size=4pt,
					inner sep=2pt},
				weight/.style={font=\scriptsize},
				]
				
				
				\draw[ blue, ultra thick] (1,0) node  [vertex, fill=blue, label=above: virus negative in location $C_2$ ] {};
				\draw[ultra thick] (0.5,-1) node  [vertex, fill=black,label=left:] {};	
				\draw[ ultra thick] (1.5,-1) node  [vertex, fill=black,label=right:] {};
				\draw[red, ultra thick] (1,-2) node  [vertex, fill=red, label=right: virus positive in location $C_2$] {};

				\draw[blue, ultra thick] (-3,-3) node  [vertex, fill=blue, label=left: virus negative in location $C_1$] {};
				\draw[ ultra thick] (-3.5,-4) node  [vertex, fill=black,label=left: ] {};	
				\draw[ ultra thick] (-2.5,-4) node  [vertex, fill=black,label=right:] {};
				\draw[red, ultra thick] (-3,-5) node  [vertex, fill= red, label=below: virus positive in location $C_1$] {};

				
				\draw[ ultra thick] (1,0) --  (0.5,-1) [] {};
				\draw[ ultra thick] (0.5,-1) --  (1.5,-1) [] {};
				\draw[ ultra thick] (1.5,-1)  --  (1,0) [] {};
				\draw[ultra thick] (1.5,-1) --  (1,0) [] {};
				\draw[ultra thick] (0.5,-1) --  (1,-2) [] {};
				\draw[ultra thick](1,-2) --  (1.5,-1) [] {};
				\draw[ ultra thick] (1,0) --  (1,-2) [] {};

				\draw[densely dotted, ultra thick] (-3.5,-4) --  (-2.5,-4) [] {};
				\draw[densely dotted, ultra thick] (- 3,-3) --  (-2.5,-4) [] {};
				\draw[ ultra thick] (-3.5,-4)  --  (-3,-3) [] {};
				\draw[ ultra thick] (-3.5,-4) --  (-3,-5) [] {};
				\draw[densely dotted, ultra thick](-3,-5) --  (-2.5,-4) [] {};
				\draw[densely dotted, ultra thick] (-3,-3) --  (-3,-5) [] {};


				\draw[blue, ultra thick] (1,0) --  (-3,-3) [] {};
				\draw[ ultra thick] (0.5,-1) --  (-3.5,-4) [] {};
				\draw[densely dotted, ultra thick] (1.5,-1) --   (-2.5,-4) [] {};
				\draw[red, ultra thick] (1,-2) --  (-3,-5) [] {};

%
			\end{tikzpicture}
	\end{center}}
	\caption{The model of interactions between two locations in a quarantine.}
	\label{fig:pandemic}
\end{figure}
\par Notice that a doubly warped product model works very well for this scenario since as an agent moves from one geographical block to the other, the type of that agent does not change; that is why in this model there is no edge between two different categories in two different geographical blocks; see Figure~\ref{fig:pandemic}. 
\par Of course, this is an evolutionary (time dependent) model. In such a quarantine model, one wishes to reach a balance between reducing the interactions and robustness since lessening interactions too much would decrease the curvature and as a result the robustness (since curvature is positively correlated with edge weights). 
\par \textsf{\small \bf Robustness speculations:}\\
Based on the curvature reasoning similar to what is done in Toy Models I and III, one observes that under the Anzats~\ref{anz:robustness}
(which is substantiated experimentally), \emph{the robustness of this quarantine network decreases as the mismatch between interactions in the capital (as compared to other locations) is too big and if there is one category which is favored above the other categories (in terms of leniency on interaction policies) to the extreme degree.} 
\end{document}